\documentclass[leqno,12pt]{amsart}

\newlength{\textlarg}

\usepackage{amssymb, ulem}
\usepackage{euscript}
\usepackage{mathrsfs} %%% Nice caligraphic fonts
\usepackage[latin1]{inputenc}
\usepackage{amsmath}
\usepackage{amsthm}
\usepackage{amsfonts}
\usepackage{mathabx}%%%Allows to use ``\widecheck''
\usepackage{enumerate}

 \usepackage{graphicx}% Include figure files
 \usepackage{float} %Pour fixer l'endroit du dessin avec [H]
 
 \usepackage[colorlinks = true, linkcolor = blue, urlcolor  = blue, citecolor = blue, anchorcolor = blue]{hyperref} %%Allows to give details of referred items

 \usepackage{pgf,tikz}
\usepackage{mathrsfs}
\usetikzlibrary{arrows}
\usetikzlibrary{shapes}
\usetikzlibrary{backgrounds}
\evensidemargin\oddsidemargin
\def\thesection{\arabic{section}}
\renewcommand{\theequation}{\thesection.\arabic{equation}}

\newtheorem{theorem}{Theorem}[section]
\newtheorem{lemma}[theorem]{Lemma}
\newtheorem{proposition}[theorem]{Proposition}
\newtheorem{corollary}[theorem]{Corollary}

\newtheorem{notation}[theorem]{Notation}
\newtheorem{definition}[theorem]{Definition}
\theoremstyle{definition}   % text in "Roman"
\newtheorem{remark}[theorem]{Remark}

\newcommand{\eqnsection}{
\renewcommand{\theequation}{\thesection.\arabic{equation}}
    \makeatletter
    \csname  @addtoreset\endcsname{equation}{section}
    \makeatother}
\eqnsection
\def\r{{\mathbb R}}
\def\e{{\mathbb E}}
\def\p{{\mathbb P}}
\def\d{\mathrm{d}}
\begin{document}

%%%%%%%%%%%%%% Beginning of the text

 %\vglue30pt

\centerline{\large\bf   
The stochastic Jacobi flow}

\bigskip
\bigskip

 \centerline{by}

\medskip

 \centerline{Elie A\"{i}d\'ekon\footnote{\scriptsize SMS, Fudan University, China, {\tt aidekon@fudan.edu.cn} partially supported by NSFC grant QXH1411004.},   Yueyun Hu\footnote{\scriptsize LAGA, Universit\'e Sorbonne Paris Nord, France, {\tt yueyun@math.univ-paris13.fr}}, 
and Zhan Shi\footnote{\scriptsize  AMSS,
Chinese Academy of Sciences,  China, {\tt shizhan@amss.ac.cn}}}

\bigskip
\bigskip

\bigskip
%\centerline{\scriptsize This version:   01/06/2023 }
\bigskip

{\leftskip=2truecm \rightskip=2truecm \baselineskip=15pt \small

\noindent{\slshape\bfseries Summary.} The problem of conditioning on the occupation field was investigated for the Brownian motion in 1998 independently  by Aldous \cite{aldous}  and Warren and Yor  \cite{warren-yor} and recently for the loop soup at intensity $1/2$ by Werner \cite{wernermarkov}, Sabot and Tarr\`es \cite{st16}, and Lupu, Sabot and Tarr\`es \cite{lst19}.  We consider this problem in the case of the Brownian loop soup on the real line, and show that it is connected with a flow version of Jacobi processes, 
 called Jacobi flow. We give a pathwise construction of this flow simultaneously for all parameters by means of a common Brownian motion,
 via the perturbed reflecting Brownian motion. The Jacobi flow  is related to Fleming--Viot processes, as established by Bertoin and Le Gall \cite{bertoin-legall-II} and Dawson and Li \cite{dawson-li12}. This relation allows us to interpret  Perkins' disintegration theorem between Feller continuous state branching-processes and Fleming--Viot processes as a decomposition of Gaussian measures.  Our approach gives a unified framework for the problems of disintegrating on the real line. The connection with Bass--Burdzy flows which was drawn in Warren \cite{warren} and Lupu, Sabot and Tarr\`es \cite{lst20} is shown to be valid in the general case.

\medskip


\medskip
 

} %%%%%% End of narrower

\bigskip
\bigskip

\section{Introduction}

A Fleming--Viot process is a measure-valued branching process which models the evolution of a population of constant size $1$.  In \cite{bertoin-legall-I}, Bertoin and Le Gall construct a generalized version of the process through a flow of bridges, which are in correspondance with exchangeable coalescents. We can understand this flow as follows (it is  the flow $\widehat B$ in the notation of \cite{bertoin-legall-I}). One represents the population at any time as the interval $[0,1]$, each point of the interval (the label) representing a particle. Only a finite number of particles at time $0$ will beget descendants at time $t>0$. Partitioning the population at time $t$ into families will yield a finite partition of the interval $[0,1]$, the $i$-th interval from the left representing the descendants at time $t$ of the ancestor at time $0$ with the $i$-th  lowest label.  One can then construct at each time $t$ a piecewise constant bridge from $0$ to $1$. The bridge jumps at the label of a particle at time $0$ which has descendants at time $t$, the size of the jump being the size of its offspring, and stays constant otherwise. One can proceed similarly between times $s$ and $t$ (and not only $0$ and $t$), and obtain a collection of bridges which naturally satisfies a flow property. Incorporating immigration in the model translates into adding an extra jump at $0$ to the bridge, see \cite{foucart12}. In the case where the coalescents are so-called $\Lambda$-coalescents, Dawson and Li \cite{dawson-li12}  consider the flow ${\mathcal Y}$ such that the bridges of \cite{bertoin-legall-I}, \cite{foucart12} are the maps  $v \in [0,1]\to {\mathcal Y}_{s,t}(v)$, with ${\mathcal Y}_{s,t}(0)$ representing the size of the population descending from immigrants which arrived at rate $\delta$ between time $s$ and $t$. They construct it as a solution of a stochastic differential equation (SDE) driven by a white noise and a Poisson random measure.

We  consider the flow ${\mathcal Y}$ in the case where the SDE appearing in \cite{dawson-li12} is only driven by the white noise.  The flow lines of ${\mathcal Y}$ are  the solutions of  the SDE
\begin{equation}\label{eq:jacobi}
 \d  Y_t = 2\sqrt{ Y_t(1-  Y_t)} \d \gamma_t + \left(\delta (1-Y_t) - \delta' Y_t \right)\d t, 
 \end{equation}

\noindent where $\gamma$ denotes a one-dimensional Brownian motion, and $\delta, \delta'\in \r$.  In the case $\delta=\delta'=0$, the dual of ${\mathcal Y}$ corresponds to the classical Kingman coalescent.
Following \cite{warren-yor} we call the solution $(Y_t)$ of \eqref{eq:jacobi} a ${\rm Jacobi}(\delta,\delta')$ process.  These processes arise as the ratio of Gamma processes, see Proposition 8 in \cite{warren-yor} in the case of nonnegative parameters, and Theorem 4 in \cite{pal}. A ${\rm Jacobi}(\delta,\delta')$ process is called in  population genetics  a Wright--Fisher  diffusion with mutation rates $(\delta,\delta')$.

The Jacobi flow ${\mathcal Y}$ is the flow version of the Jacobi processes ${\rm Jacobi}(\delta,\delta')$, see Definitions \ref{d:jacobi} and \ref{d:jacobi-bis}. %for the ${\rm Jacobi}(\delta,\delta')$ flow with $\delta, \delta'\in \r$.
    
One of the goals of this paper is to give a pathwise construction, from a one-dimensional two-sided Brownian motion, of ${\rm Jacobi}(\delta,\delta')$ flows  simultaneously for all $\delta, \delta' \in \r$.  This is done by means of the two-sided  perturbed reflecting Brownian motion defined as follows.  Let $(B_t, \, t\ge 0)$ and $(B_t', \, t\ge 0)$ be two independent standard one-dimensional Brownian motions. Denote by ${\mathfrak L}$ and ${\mathfrak L}'$ their associated local time processes at position zero. The two-sided Brownian motion $(B_t,\,t\in \r)$ is defined as
$$
B_t := \left\{ 
\begin{array}{ll}
B'_{-t} & t\le 0, \\
B_t & t\ge 0. 
\end{array}
\right.
$$

\noindent We also let ${\mathfrak L}_{t}:= - {\mathfrak L}_{-t}'$ for $t\le 0$. For $\mu>0$, we define the two-sided  perturbed reflecting Brownian motion (PRBM) or $\mu$-process  by
\begin{equation} \label{def:muprocess}
X_t :=
\left\{
\begin{array}{ll}
 |B_{-t}'| + \mu {\mathfrak L}'_{-t}, &\hbox{ if } t\le 0,\\
  |B_t| - \mu {\mathfrak L}_t, &\hbox{ if } t\ge 0,
 \end{array}
 \right.
\end{equation}

\noindent which can simply be written as
$X_t := |B_t| - \mu {\mathfrak L}_t,\, t\in \r$.   [The case $\mu=1$ is special: on the positive half-line, it is distributed as a Brownian motion while in the negative half-line, it is a time-reversed three-dimensional Bessel process.] For general properties of PRBM, see \cite{legall-yor, yor92} and the references therein.  

The promised pathwise construction of ${\rm Jacobi}(\delta,\delta')$ flows is as follows: 

\bigskip
{\noindent\bf Theorem I (Theorem \ref{t:perkins}).} {\it  Let $\delta, \delta' \in \r$. Let ${\mathcal Y}$ be defined as in \eqref{def:Y2}. Then ${\mathcal Y}$ is a ${\rm Jacobi}(\delta,\delta')$ flow. 
}

\bigskip

We will see that the flow ${\mathcal Y}$  constructed in \eqref{def:Y2} is a measurable function of the PRBM $(X_t)$. Furthermore, ${\mathcal Y}$ satisfies an SDE driven by a specific martingale measure, see \eqref{eq:Y1}. This last result can be interpreted as  a ``flow version'' of the well-known Perkins' disintegration theorem where we disintegrate, with respect to a flow line,   an associated flow of squared Bessel processes; see Section \ref{s:BESQflow} for the definition and construction of this associated flow of squared Bessel processes. In particular, Theorem \ref{t:perkins} shows that Perkins' disintegration theorem between Feller continuous state branching-processes and Fleming--Viot processes amounts to a decomposition of Gaussian measures.

Another goal of the paper is to describe the ``contour function'' of Jacobi flows. We show that this contour function is a version of the PRBM conditioned on its occupation field. As such, it gives a unified framework for related problems  on disintegrating a standard Brownian motion \cite{aldous, warren-yor, warren, berestycki, GKW} and inverting the Ray--Knight identity of the Gaussian free field on the real line \cite{lst20}. The connection with Bass--Burdzy flows which was drawn in these papers is shown to be valid in the general case. It provides a construction  of the  PRBM, or equivalently (\cite{lupu18}) of the one-dimensional loop soup, conditioned on its occupation field. The corresponding result for the loop soup at intensity $\frac12$ on a metric graph is given in \cite{elie}.

The problem of conditioning a Brownian motion  on its occupation field has been treated by Aldous \cite{aldous} and Warren and Yor \cite{warren-yor}, via different approaches. Aldous used  the tree structure of the Brownian excursion to show that the genealogy of the conditioned Brownian motion is a time-changed Kingman coalescent (which is the dual of a Fleming--Viot process). Later, Berestycki and Berestycki \cite{berestycki}  found an analogous result with excursion theory.  Warren and Yor \cite{warren-yor} solved the question by constructing the conditioned Brownian motion in time (rather that in space, as in the case of Aldous \cite{aldous}). If $|B|$ is a reflecting Brownian motion, they define a process $Z$, called Brownian burglar, by
$$
Z_t := \int_0^{|B|_{A_t}} {\d r \over L^{|B|}(\tau_1^0(|B|),r)}, \qquad
A_t:= \inf\left\{ s >0\,:\, \int_0^s {\d r \over L^{|B|}(\tau_1^0(|B|),|B_r|)^2} > t\right\} 
$$

\noindent where $L^{|B|}$ and $\tau(|B|)$ are naturally the local time and the inverse local time of $|B|$. They showed that $Z$ is independent of the occupation field $(L^{|B|}(\tau_1^0(|B|),r),\, r\ge 0)$. In these works, the link with the Fleming--Viot process has been suggested but not made explicit.  In \cite{warren-yor}, the authors  state that ``{\it the results [...] can be seen as describing a contour process for the Fleming--Viot process}''. Recently, Gufler, Kersting and Wakolbinger \cite{GKW} gave a rigorous connection between these two models by constructing the Brownian excursion via an enriched version of the lookdown process, hence giving a precise meaning to the statement of \cite{warren-yor}.

In our case, we will recover the connection with Fleming--Viot processes by constructing  $Z$  through our version of Perkins' disintegration theorem. The idea, which originates in a paper of T\'oth and Werner \cite{toth-werner}, is to define the burglar via its local time flow.  Actually, this approach will allow us to solve the analogous problem for the whole class of PRBM, giving rise to a family of ``burglars''.

To construct the burglars, we define  
$$
L(t, \, r) 
:=
\lim_{\varepsilon\to 0} \frac{1}{\varepsilon} \int_{-\infty}^t  {\bf 1}_{\{ r \le X_s \le r+\varepsilon\} } \, \mathrm{d} s,
\qquad
t\in \r , \, r\in \r,
$$

\noindent as the local time of the continuous semimartingale $X$ at time $t$ and position $r$; we work with a bicontinuous version of local times $L(\cdot, \cdot)$, as in (\cite{revuz-yor}, Theorem VI.1.7).  Let
\begin{equation}
    \label{def-tau}
    \tau^r_t:=\inf\{ s\in \r\,:\, L(s,r)> t \}  
\end{equation}

\noindent be the inverse local time of $X$. Let
$$
t^*_0:=\sup\left\{t\le 0\,:\, L(\tau_1^0,X_t)=0\right\}, \qquad 
t^*:=\inf\left\{t\ge 0\,:\, L(\tau_1^0,X_t)=0 \right\},
$$ 

\noindent with the convention $\inf\emptyset=+\infty$ and $\sup\emptyset =-\infty$. Define for $t\in \r$,
\begin{equation}\label{eq:Z10}
Z^{(1)}_t := \int_0^{X_{A^{(1)}_t}} {\d r \over L(\tau_1^0,r)}, \qquad
A^{(1)}_t:= \inf\left\{ s  \in \r \,:\, \int_0^s {\d r \over L(\tau_1^0,X_s)^2} > t\right\}.
\end{equation}

\noindent  The process $Z^{(1)}$ is the burglar associated with the process $X$ between times $t_0^*$ and $t^*$. Similarly, let for $t\ge 0$,
\begin{equation}\label{eq:Z20}
Z^{(2)}_t := \int_0^{X_{A^{(2)}_t}} {\d r \over L(\tau_1^0,r)}, \qquad
A^{(2)}_t:= \inf\left\{ s  \ge 0 \,:\, \int_{\tau_1^0-s}^{\tau_1^0} {\d r \over L(\tau_1^0,X_s)^2} > t\right\}.
\end{equation}

\noindent The burglar $Z^{(2)}$ is associated with the process $X$ between times $t^*$ and $\tau_1^0$. A precise statement of the following result  will be given in Section \ref{s:burglar}.

\bigskip
{\noindent\bf Theorem II (Theorem \ref{p:loctime}).} {\it 
The processes $Z^{(1)}$ and $Z^{(2)}$ are independent of $(L(\tau_1^0,r),\, r\in \r)$. The local time flows of the processes $Z^{(1)}$ and $Z^{(2)}$ are the left part and the right part of a flow which is a {\rm Jacobi($\delta,0$)} flow in the positive time-axis and a {\rm Jacobi($\delta,2$)} flow in the negative time-axis. 
}

The above result  gives  a description of the conditional law of a PRBM given its occupation field up to $\tau_0^1$ (and from $t_0^*$). It suffices to take the processes $Z^{(1)}$ and $Z^{(2)}$ independently of the occupation field, and invert the transformations \eqref{eq:Z10} and \eqref{eq:Z20}.    We can also obtain a burglar by disintegrating the positive part of a PRBM with respect to its occupation field, see Section \ref{s:FV}. This burglar can be interpreted as the contour function of the ${\rm Jacobi}(\delta,0)$ flow. Each of its excursion away from $0$ being associated with a continuous tree in the manner of Aldous \cite{aldous} and Duquesne and Le Gall \cite{DuLG02}, the burglar can be viewed as the contour function of a Fleming--Viot forest.

\bigskip

We present now an interesting  connection with the Bass--Burdzy flow. This flow was introduced by Bass and Burdzy in \cite{bassburdzy}. The Bass--Burdzy flow with parameters $(\beta_1,\beta_2)$, for $\beta_1,\, \beta_2\in \r$, is the collection of processes $({\mathcal R}_t(x),\,t\ge 0,\,x\in \r)$ which are the strong solutions of 
$$
{\mathcal R}_t(x) = x +  \gamma_t + \beta_1\int_0^t {\bf 1}_{\{{\mathcal R}_s(x)< 0\}} \d s +\beta_2\int_0^t {\bf 1}_{\{{\mathcal R}_s(x)> 0\}} \d s
$$

\noindent where $\gamma$ is a standard Brownian motion common for all processes ${\mathcal R}_{\cdot}(x)$. In words, ${\mathcal R}_t(x)$ is a Brownian motion with drift $\beta_1$ when negative and $\beta_2$ when positive. Following \cite{hw00}, for $t\ge 0$, denote by ${\mathcal R}_t^{-1}(0)$ the real $x$ such that ${\mathcal R}_t(x)=0$. When $\beta_1=0$ and $\beta_2=1$, the process ${\mathcal R}_t^{-1}(0)$ is shown in Warren \cite{warren} to be a time-change of the Brownian burglar of Warren and Yor \cite{warren-yor}. In a recent paper \cite{lst20}, Lupu, Sabot and Tarr\`es showed that in the case $\beta_1=-\frac12, \beta_2=\frac12$, the process ${\mathcal R}_t^{-1}(0)$ is the scaling limit of
a self-interacting process involved in the inversion of the Ray--Knight identity, see the next paragraph. We show that in general, one can recover all the burglars $Z^{(2)}$ associated with the PRBM.

\bigskip

{\noindent\bf Theorem III (Theorem \ref{p:bassburdzy} and Proposition \ref{p:parameters}).} {\it 
Taking $\beta_1 ={\delta\over 2}-1$ and $\beta_2={\delta\over 2}$, the process ${\mathcal R}_t^{-1}(0)$ is a time-change of the process $Z^{(2)}$ defined in \eqref{eq:Z20}.
}

\bigskip

To prove the theorem, we use a new approach based on a renewal argument. The process $Z^{(2)}$ has a kind of Markov property. If one appropriately scales the process after time $t$, the scaled process will be independent of the past and with always the same distribution, see Theorem \ref{t:markov}. We use this property to show that the images of a point on the real line by these transformations form actually a L\'evy process up to hitting $0$, and then deduce that it must be a flow line of a Bass--Burdzy flow for some parameters. Finally, we use some random variable whose distribution is easily identified to compute the parameters. This renewal property breaks down when looking at $Z^{(1)}$. Roughly speaking, one needs to remember where the infimum was before doing the scaling. Still,  we believe that a similar connection should hold with some perturbed Bass--Burdzy flow.

\bigskip

Yet another, non-trivial, application of our construction concerns the inversion of the Ray--Knight identity on the line and its connection with loop soups.  Consider the PRBM $X$. By excursion theory, the excursions of $X$ above its infimum process forms a Poisson point process of Brownian excursions rooted at various points of the real line. By exploring the real line in the upwards direction, one can consider the excursions that one encounters along the way as Brownian loops rooted at their minimum. As shown by Lupu \cite{lupu18}, the collection of loops has the law of a Brownian loop soup on the real line with intensity ${1\over \mu}$. %See Figure \ref{f:muprocess}. 
As a result, conditioning a PRBM on its occupation field may be seen  as conditioning  a loop soup on its occupation field.

A similar problem was studied recently by Lupu, Sabot and Tarr\`es \cite{lst19} in the case $\mu=2$, equivalently $\delta=1$, which is the case of the loop soup at intensity ${1\over 2}$, related to the Gaussian free field (the local time is ${\rm BESQ}^1$, hence the square of the Gaussian Free Field, which is simply the Brownian motion on the real line). In their setting (reformulated in terms of loop soups), the authors add a Brownian motion up to a fixed local time  at zero to a loop soup on $\r\backslash\{0\}$. Conditionally on the occupation field, they manage to reconstruct the Brownian motion path up to a random time. The reconstruction process is actually in terms of the burglar $Z^{(2)}$, stopped when it reaches local time $1$ at some position (which happens in a finite time almost surely).

We solve the problem for any intensity of the loop soup. Specifically, take a loop soup in the positive half-line, and add Brownian excursions up to local time $1$. We get the positive part of a PRBM. A space-time transformation of this PRBM gives a burglar $Z^+$ similar to $Z^{(2)}$ in \eqref{eq:Z20}, which we have interpreted as the contour function of the ${\rm Jacobi}(\delta,0)$ flow; see \eqref{eq:Z+} for the definition of $Z^+$. One then gets a reconstruction of the Brownian excursions added to the loop soup in terms of the Jacobi flow, or in terms of a Bass--Burdzy flow driven by a reflecting Brownian motion with drift. More precisely, we have

\bigskip
{\noindent\bf Theorem IV (Theorems \ref{t:RK1} and \ref{t:FV}).} {\it Let $Z^+$ be as in \eqref{eq:Z+} and the occupation field $f$ as in \eqref{def:f}. Then the process $Z^+$ is independent of $f$ and is the burglar associated with the positive part of the PRBM. }

\bigskip

Finally, we may wonder whether it is possible to find the law of a loop soup conditioned on its occupation field on a cable graph introduced by Lupu in \cite{lupu16}. In general, the loop soup loses its Markovian properties. But in the case of intensity ${1\over 2}$, the relation of the loop soup with the Gaussian free field indicates that a certain Markov property should hold. It is the topic of \cite{elie}, where the description of the conditioned loop soup uses the framework of this paper. In the discrete case, the analogous problem was answered by Werner \cite{wernermarkov} (description that we could  label in ``space'', via the link with the random current model) and Sabot, Tarr\`es \cite{st16}, Lupu, Sabot and Tarr\`es \cite{lst19} (description in ``time'', via a self-interacting process).

\bigskip

The paper is organized as follows. Flows of BESQ processes are studied in Section \ref{s:BESQflow}. Section \ref{s:transform} collects all the transformations that will be needed to state the disintegration theorems and construct the burglars. Section \ref{s:perkins} contains the disintegration result stated for the ${\rm BESQ}^{\delta,\delta'}$ flows and as a corollary the construction of the Jacobi flows from the PRBM. Theorem II on the conditioning of a PRBM  is proved in Section \ref{s:burglar}. In Section \ref{s:bassburdzy}, we will state the Markov property for $Z^{(2)}$ and prove the link with the Bass--Burdzy flow given in Theorem III. Section \ref{s:FV} connects the contour function of the ${\rm Jacobi}(\delta,0)$ flow to the Bass--Burdzy flow driven by a reflecting Brownian motion with drift.  Section \ref{s:perfectflow} studies properties of the BESQ and Jacobi flows with emphasis on perfect flow properties, and Section  \ref{s:girsanov} contains a Girsanov theorem for Jacobi flows.

\bigskip

\noindent {\bf Acknowledgements}. We thank  Hui He, Zenghu Li, Titus Lupu, Emmanuel Schertzer, Anton Wakolbinger for helpful discussions. The research of EA was partially supported by NSFC grant QXH1411004.

\section{The BESQ flows}
\label{s:BESQflow}

We give the definition of the BESQ flows,  then embed them in the PRBM.  For any bounded Borel function $g:\r_+\times \r \to \r$ with compact support,  let 
\begin{equation}\label{def:W}
W(g):= \int_{-\infty}^\infty g\big(L(t, X_t), X_t\big)  {\rm sgn}(B_t) \d B_t.
\end{equation}

\noindent As shown in \cite{eyzRK},  the  stochastic integral is well defined and $W$ is a white noise on $\r_+\times \r$.  Define for any $x\ge r$ in $\r$ and any $a\ge 0$,
\begin{equation}\label{def:S}
  S_{r,x}(a) := L(\tau_a^r, x ), \qquad S^*_{r,x}(a):= L(\tau_a^{-r}, -x ).
\end{equation}

\begin{theorem} [\cite{eyzRK}, Theorem 5.1]  \label{t:RK}
Let $W$ be the white noise defined via \eqref{def:W}. For any $r\in \r$ and any $a\ge 0$, $S_{r,r+\cdot}(a)$ and $S^*_{r,r+\cdot}(a)$ are the pathwise unique solutions, which are strong\footnote{By strong solution, we mean that for instance $S_{r,r+h}(a)$ is measurable with respect to $\sigma\{W(\cdot, [r, s]): r\le s \le r+h\}$. This terminology will be used elsewhere with the same remark.}, of the following SDEs:
\begin{eqnarray}
\label{eq:S1}
S_{r,r+h}(a) &=& 
a + 2 \int_r^{r+h}  W([0,S_{r,s}(a)], \d s) + \frac{2}{\mu} h, \, h\ge 0, \\
\label{eq:S2}
S^*_{r,r+h}(a) &=& 
a - 2 \int_{r}^{r+h}  W^*([0, S^*_{r,s}(a)], \d s) + (2-\frac{2}{\mu}) h, \, h\in [0, T_0(S^*_{r,r+\cdot}(a))],
\end{eqnarray}

\noindent   where  $W^*$ is the image of $W$ under the map $(a, s)\mapsto (a, -s)$ and $T_0(S^*_{r,r+\cdot}(a)):=\inf\{h\ge 0\,:\, S^*_{r,r+h}(a) =0\}$.
\end{theorem}

{\bf Remark}. If $I_t:=\inf_{s\le t} X_s$ denotes the infimum process of $X$, we can rewrite $T_0(S^*_{r,r+\cdot}(a))$ in \eqref{eq:S2} as $- I_{\tau_a^{-r}}-r$.

\bigskip

Ray--Knight theorems are usually statements on marginal distributions. Recall that the squared Bessel process of dimension $\delta \in \r$ started at $x\ge 0$, denoted ${\rm BESQ}^\delta_x$,  is the pathwise unique solution of
$$
S_t = x+2\int_0^t \sqrt{|S_s|} \, \d \gamma_s + \delta t,\qquad t\ge 0,
$$

\noindent where, as before,
$\gamma$ is a standard Brownian motion. The ${\rm BESQ}^\delta$ hits zero at a positive time if and only if $\delta< 2$. It is absorbed at $0$ when $\delta= 0$ and is  reflecting
at $0$ when $\delta\in (0,2)$. When $\delta<0$, after hitting $0$, it behaves as   a ${\rm BESQ}^{-\delta}_0$ in the negative half-line, see e.g.\ \cite{gy03}. See Le Gall and Yor \cite{legall-yor} and Yor (\cite{yor92}, Chapter 9) for references on Ray--Knight theorems.

\bigskip

Ray--Knight theorems show that the flows of squared Bessel processes are embedded in the PRBM. The setting of Dawson and Li \cite{dawson-li12} includes the construction of such flows. We impose some further regularity conditions in order to give it the structure of a flow in the sense of \cite{arratia, toth-werner}.

\begin{definition}\label{d:BESQflow}
Let $\delta >0$.  We call ${\rm BESQ}^\delta$ flow (or non-killed ${\rm BESQ}^\delta$ flow) a collection ${\mathcal S}$ of continuous processes $({\mathcal S}_{r,x}(a),\, x\ge r)_{r\in \r,a\ge 0}$ such that: 
\begin{enumerate}[1)]
\item
for each $(r,a)\in \r\times \r_+$, the process $({\mathcal S}_{r,x}(a),\, x\ge r)$ is almost surely the strong solution of the following SDE  
\begin{equation}\label{def:BESQflow}
{\mathcal S}_{r,x}(a) = 
a + 2 \int_r^x  {\mathcal W}([0,{\mathcal S}_{r,s}(a)], \d s) + 
\delta (x-r) 
\end{equation}

\noindent where ${\mathcal W}$ is a white noise on $\r_+\times \r$.

\item Almost surely,
\begin{enumerate}[(i)]
\item for all $r\in \r$ and $a\ge 0$, $ {\mathcal S}_{r,r}(a)=a$,% and $\lim_{r'\uparrow r} {\mathcal S}_{r',r}(a)=a$, 
\item for all $r\le x$, $a\mapsto  {\mathcal S}_{r,x}(a)$ is c\`adl\`ag,
\item for all $r'\le r$ and all $a',a\ge 0$, if $ {\mathcal S}_{r',r}(a')> a$ (resp. $ {\mathcal S}_{r',r}(a')< a$), then $ {\mathcal S}_{r',x}(a')\ge {\mathcal S}_{r,x}(a)$ (resp. $ {\mathcal S}_{r',x}(a')\le {\mathcal S}_{r,x}(a)$) for all $x\ge r$.
\end{enumerate} 
\end{enumerate}
\end{definition}

\bigskip

\begin{definition}\label{def:BESQflowgeneral}
For $-\infty<\delta< 2$, we call killed ${\rm BESQ}^\delta$ flow the flow solution of \eqref{def:BESQflow}, where the process is absorbed when hitting $0$, and which satisfies the same regularity conditions.

For short, we call general ${\rm BESQ}^\delta$ flow a killed ${\rm BESQ}^\delta$ when $\delta\le 0$, a  (non-killed) ${\rm BESQ}^\delta$ flow when $\delta\ge 2$ and either a killed or non-killed ${\rm BESQ}^\delta$ flow when $\delta\in (0,2)$. 
\end{definition}

Dawson and Li \cite{dawson-li12} showed (in a more general setting) that equation \eqref{def:BESQflow} possesses a pathwise unique solution. We can show that the same is true in the case $\delta<0$, with arguments similar to the ones of the proof of Theorem \ref{t:perkins}.

\bigskip

We check that the BESQ flows are naturally embedded in the PRBM $X$, which will give the existence of these flows for free.  In agreement with \eqref{def:S}, we define  $S_{r,x}(a) = L(\tau_a^r, x )$ for any $r,x\in \r$. We call the collection of processes ${\mathcal L}_X:= (S_{r,x}(a),\, -\infty < r,x< \infty,\, a\ge 0)$ the {\it local time flow} of $X$. The flow $S:=(S_{r,x}(a),\, x\ge r)_{r\in\r,a\ge 0}$ is called the {\it forward local time flow of $X$}, while $S^*:=(S_{-r,-x}(a),\, x\ge r)_{r\in\r,a\ge 0}$ is called the {\it backward local time flow of $X$}. Both flows are dual as we will see in equation \eqref{eq:Sdual}.

\begin{proposition}\label{p:Sembedded}
Let $\mu>0$. 

The forward flow $(S_{r,x}(a),\, x\ge r)_{r\in \r,a\ge 0}$ is a ${\rm BESQ}^{2/\mu}$ flow.

If $\mu\in(0,1]$,  the backward flow  $(S_{-r,-x}(a),\, x\ge r)_{r\in\r,a\ge 0}$ is a killed ${\rm BESQ}^{2-2/\mu}$ flow.

\end{proposition}

\noindent {\it Proof}. The finite-dimensional distributions coincide by Theorem \ref{t:RK} so we only have to check the regularity conditions (i), (ii) and (iii) in Definition \ref{d:BESQflow}. We have $S_{r,r}(a)=a$ indeed. Statement (ii) is a consequence of the continuity of the local times, and the observation that $a\to \tau_{a}^r$ is c\`adl\`ag by construction. We prove now (iii). We have $S_{r',r}(a')=L(\tau_{a'}^{r'},r)$. By definition of $\tau_{a}^r$, $S_{r',r}(a')>a$ is equivalent to  $\tau_{a'}^{r'}>\tau_a^r$. Therefore,  $S_{r',x}(a')=L(\tau_{a'}^{r'},x) \ge L(\tau_{a}^{r},x)=S_{r,x}(a)$. 
$\Box$

\bigskip

Conversely, we can recover $X$ from its local time flow. It is the content of the following proposition.
\begin{proposition}\label{p:Xflow}
The process $X$ is a measurable function of its local time flow ${\mathcal L}_X=(S_{r,x}(a),\, -\infty < r , x < \infty,\, a\ge 0)$.
\end{proposition}

\noindent {\it Proof}.  By the occupation times formula, for any $a\ge 0$ and $b$,  $\tau_a^b=\int_{\r} (S_{b,r}(a) -S_{0,r}(0))\d r$ hence $(\tau_a^b,\,a\ge 0,\, b\in \mathbb{R})$ is measurable with respect to ${\mathcal L}_X$.  Therefore $L(t,x), t\in\r, x\in \r,$ is measurable with respect to ${\mathcal L}_X$, which again by the occupation times formula yields that $X$ is measurable as well (for any $s<t$, $\int_s^t X_u \d u = \int_\r x (L(t, x)-L(s,x)) \d x$).  
$\Box$

\bigskip

We show now that any ${\rm BESQ}^\delta$  flow can be constructed from a countable number of flow lines. It allows us to identify  ${\rm BESQ}^\delta$ flows with flows embedded in a PRBM.

\begin{proposition}
Let $\mathcal S$ be a general ${\rm BESQ}^\delta$ flow. Let $(r_n,a_n)_n$ be a dense countable set in $\r\times \r_+$. Almost surely, for any $x\ge r$ and $a\ge 0$,
$$
{\mathcal S}_{r,x}(a)=\inf_{\{n\,:\, r_n\le r,\,{\mathcal S}_{r_n,r}(a_n)> a\} }\mathcal S_{r_n,x}(a_n).
$$
\end{proposition}

\noindent {\it Proof}. We first check that almost surely, for any $r\in \r$, and any $0\le a<a'$, one can find some $(r_n,a_n)$ such that $r_n\le r$ and ${\mathcal S}_{r_n,r}(a_n)\in (a,a')$. From Theorem \ref{t:RK}, one can reason on the local times of the PRBM   (notice that one can make this identification because we only look at a countable number of flow lines). The previous claim property follows since  it suffices to take $a<b_1<b_2<a'$, and $(r_n,a_n)$ such that $L(\tau^{b_1}_r,r_n)<a_n<L(\tau^{b_2}_r,r_n)$. We would then have $L(\tau_{r_n}^{a_n},r)\ge L(\tau^{b_1}_r,r)>a$ and $L(\tau_{r_n}^{a_n},r)\le L(\tau^{b_2}_r,r)<a'$.

Let us go back to the proof of the proposition. Fix $x\ge r$ and $a\ge 0$. We prove that ${\mathcal S}_{r,x}(a)=\inf_{\{n\,:\, r_n\le r,\,{\mathcal S}_{r_n,r}(a_n)> a\} }\mathcal S_{r_n,x}(a_n)$. 
Notice that the set over which the infimum is taken is not empty.  
For  any $(a_n,r_n)$ such that ${\mathcal S}_{r_n,r}(a_n)> a$, we have ${\mathcal S}_{r_n,x}(a_n) \ge {\mathcal S}_{r,x}(a) $ by (iii) of Definition \ref{d:BESQflow}, which proves one inequality. If $a'>a$, we take some $(a_n,r_n)$ such that  ${\mathcal S}_{r_n,r}(a_n)\in (a,a')$. Then   ${\mathcal S}_{r,x}(a) \le {\mathcal S}_{r_n,x}(a_n) \le {\mathcal S}_{r_n,x}(a') $. We then conclude by (ii) of Definition \ref{d:BESQflow}. $\Box$

\bigskip

\begin{proposition}\label{c:BESQdual}
Let $\delta \in \r$. Let $\mathcal S$ be a general ${\rm  BESQ}^\delta$  flow. Define its dual ${\mathcal S}^*$ by, for $r\le x$, 
$$
{\mathcal S}^*_{r,x}(a) := \inf\{b\ge 0\,:\, {\mathcal S}_{-x,-r}(b)>a\}.
$$
Then ${\mathcal S}^*$ is a ${\rm BESQ}^{2-\delta}$ flow (in the case $\delta\in (0,2)$, it is killed if ${\mathcal S}$ is not killed, and it is not killed if ${\mathcal S}$ is killed).  Moreover, $({\mathcal S}^*)^* = \mathcal S$.
\end{proposition}

This proposition gives the dual of a Feller CSBP with immigration. We refer to Foucart, Ma and Mallein \cite{fmm} for dual processes of CSBPs.

\bigskip
\noindent {\it Proof}. It is a direct consequence of Proposition \ref{p:Sembedded} and the following claim: almost surely, for all $r,x \in \r$ and $a\ge 0$,
\begin{equation}\label{eq:Sdual}
S_{r,x}(a) = \inf\{b\ge 0\,:\, S_{x,r}(b)>a\}.
\end{equation}

\noindent Let us prove this claim. Let $b\ge 0$  be
such that $S_{x,r}(b)>a$. By definition, it means that $L(\tau_b^x,r)>a$, and since $L(\tau_a^r,r)=a$,  we get $\tau_b^x>\tau_a^r$ so that $b=L(\tau_b^x,x)\ge L(\tau_a^r,x)$ which is $S_{r,x}(a)$ by definition. On the other hand, let $b\ge 0$ such that  $S_{x,r}(b)\le a$, i.e., $L(\tau_b^x,r)\le a$. Take $s\ge 0$ such that $L(s,r)>a$, hence $s>\tau_b^x$. It implies that $L(s,x) \ge L(\tau_b^x,x)=b$ and $L(\tau_a^r,x)\ge b$ by making $s\downarrow \tau_a^r$.
$\Box$

\section{Space-time transformations}
\label{s:transform}

In this section, we consider a fixed positive continuous function $f:I\to (0,\infty)$ where $I$ is an interval (not necessarily open nor bounded).
 
\subsection{Transformations of flows}
\label{s:flow}

The results in this section hold in a deterministic setting.  They provide the transformations that will be used later in our disintegration results.

\bigskip

We introduce some general notation.  Let $c\in I$ and define 
$$
\eta_{f, c}(x) := \int_c^x {\d r \over f(r)},\qquad x\in I.
$$

\begin{definition}
We denote by $\Psi(\cdot,f,c)$ the transformation with domain $D_f:=\{(a,x)\,:\, a\in [0,f(x)],\, x\in I\}$ defined as
$$
\Psi(\cdot,f,c):
\left\{
\begin{array}{ccc}
D_f &\mapsto& [0,1]\times \eta_{f,c}(I) \\
(a,x) &\mapsto& \left( {a\over f(x)},\eta_{f,c}(x)  \right).
\end{array}
\right.
$$
If $g:I\to\r_+$ is a continuous function such that $g\le f$, we denote by $\Psi(g,f,c)$ the  function whose graph is the image of the graph of $g$ by $\Psi(\cdot,f,c)$. That is
$$
\Psi(g,f,c): 
\left\{
\begin{array}{ccc}
\eta_{f,c}(I) &\mapsto& [0,1] \\
x &\mapsto&  {g\circ \eta_{f,c}^{-1}(x)\over f\circ \eta^{-1}_{f,c}(x)} .
\end{array}
\right.
$$
When $c=0$, we will simply write $\eta_{f}$ for $\eta_{f,0}$ and $\Psi(\cdot,f)$ for $\Psi(\cdot,f,0)$.
\end{definition}

 \bigskip
 
In this section, we will not need a precise definition of flows.  Let $g:I\to\r_+$ be a continuous function. We will simply call $g$-flow a collection ${\mathcal L}:=\{{\mathcal S}_{r,\cdot}(a),\, r\in I,\, a\in [0,g(r)] \}$ of continuous functions ${\mathcal S}_{r,\cdot}(a):I \mapsto \r_+$ such that  for every $ r , x \in I$ and $a\in [0,g(r)]$, ${\mathcal S}_{r, r}(a)= a $ and ${\mathcal S}_{r, x}(a) \le g(x)$. We call it flow only because we will deal exclusively with the BESQ and Jacobi flows. More generally, we give the following definition.

\begin{definition}\label{def:flowdomain}
 We call flow in the domain $D_f$ any flow which is a $g$-flow  for some continuous nonnegative function $g\le f$.
\end{definition}

{\noindent\bf Remark}. We also allow the case $f\equiv\infty$, meaning that ${\mathcal L}$ is a flow in the domain $\r_+\times \r$.   

\bigskip

If ${\mathcal L}$ is a $g$-flow, the {\it forward flow of}  ${\mathcal L}$ is the collection of functions ${\mathcal S}:=({\mathcal S}_{r,x}(a),\, x\in I\cap [r,\infty))_{r\in I,a\in [0,g(r)]}$  and the {\it backward flow} of ${\mathcal L}$ is the collection  of functions ${\mathcal S}^*:=  ({\mathcal S}_{-r,-x}(a),\, x\in (-I)\cap[r,+\infty)_{r\in -I,a\in [0,g(-r)]}$ where $-I:=\{-x,\, x\in I\}$. We already mentioned the Jacobi flow and the BESQ flow, which are in our terminology forward flows. The first one is a forward flow  in the domain $[0, 1] \times \r$ (here $g(r)\equiv 1$). The second one is a forward flow in the domain $\r_+\times \r$ (formally taking $g(r)\equiv \infty$).

\begin{definition}\label{def:imageflow}
Let $g\le f$ be a continuous nonnegative function and $c\in I$. Consider a $g$-flow ${\mathcal L}$. We denote by $\Psi({\mathcal L},f,c)$ the flow  such that for any $(v,r)$ in the image of $D_{g}:=\{(a,x)\,:\, a\in [0,g(x)],\, x\in I\}$ by $\Psi(\cdot,f,c)$, its flow line  passing through $(v, r)$ is the image  of the flow line in ${\mathcal L}$ passing through the  preimage of $(v,r)$. It is defined through the following equation: 
$$
\left(\Psi({\mathcal L},f,c)\right)_{r,x}(v) =  {{\mathcal S}_{\eta_{f,c}^{-1}(r),\eta_{f,c}^{-1}(x)}(vf \circ \eta_{f,c}^{-1}(r)) \over f\circ\eta_{f,c}^{-1}(x)}$$
for all $r , x \in\eta_{f,c}(I),\, v\in [0,\Psi(g,f,c)(r)]$.
When $c=0$, we will only write $\Psi({\mathcal L},f)$ for $\Psi({\mathcal L},f,0)$.
\end{definition}

\noindent {\bf Remark}. The image flow $\Psi({\mathcal L},f,c)$ is a $\Psi(g,f,c)$-flow.

\bigskip

If we only look at the forward flow ${\mathcal S}$ (for example in  \eqref{def:Y}), we will still write $\Psi({\mathcal S},f,c)$ for the forward flow obtained as in Definition \ref{def:imageflow}, restricted to $ r\le x$.

\subsection{Transformations of processes}
\label{s:process}

Let $J$ be an interval of $\r$ containing $0$ and ${\mathcal X}=({\mathcal X}_t,\, t\in J)$ be some real-valued continuous process such that ${\mathcal X}_t \in I$ for all $t\in J$. We suppose that ${\mathcal X}$ admits a version of bicontinuous local times $L_{\mathcal X}(t,x)$, $t\in J$, $x\in I$, defined as the densities of the occupation times: for any Borel nonnegative function $h$ and any $t\in J$, $$ \int_{J \cap (-\infty, t]} h({\mathcal X}_s) \d s= \int_{I} h(x) L_{\mathcal X}(t,x) \d x .$$

\noindent  We let $g(x):= L_{\mathcal X}(\infty, x)$ for $x\in I$  denote  the total local time of ${\mathcal X}$ at position $x$ and we suppose  that $g(x) \le f(x)$ on $I$.  Define ${\mathcal L}_{\mathcal X}=(\mathcal S_{r,x}(a))_{r,x,a}$ as, for any $x$ and $r$ in $I$ and $0 \le a \le g(r)$,  $$
{\mathcal S}_{r, x}(a):= L_{\mathcal X}(\tau^r_a({\mathcal X}), x), \  
$$

\noindent where $$\tau^r_a({\mathcal X}):= \inf\{t\in J\,:\, L_{\mathcal X}(t,r)>a\}.$$

\noindent We used the convention that $\inf\emptyset=\infty$, so that when $a=g(r)$, $\tau_a^r({\mathcal X})=\infty$ and ${\mathcal S}_{r, x}(a)=g(x)$ for all $x$. Then ${\mathcal L}_{\mathcal X} $ is a $g$-flow, which (in agreement with the previous section) we will refer to as the {\it local time flow of} ${\mathcal X}$. Its forward, resp. backward flow is called {\it forward}, resp.  {\it backward local time flow of ${\mathcal X}$}. By Definition \ref{def:flowdomain}, ${\mathcal L}_{\mathcal X}$ is a flow in the domain $D_f$.  

\begin{proposition}\label{p:locZ}
We set  $$C_f(t):=\int_0^t {{\rm d} s \over f({\mathcal X_s})^2}, \qquad t\in  J.$$

\noindent   For $c\in I$, we let $\Upsilon({\mathcal X},f,c)$ be the process  defined as
\begin{equation}\label{def:Z}
\Upsilon({\mathcal X},f,c)_t 
:=
\eta_{f,c}\left({\mathcal X}_{C_f^{-1}(t)}\right) = \int_c^{{\mathcal X}_{C_f^{-1}(t)}} {\d r \over f(r)},\qquad t\in C_f(J).
\end{equation}

Then, the process $\Upsilon({\mathcal X},f,c)$  possesses bicontinuous local times (given by \eqref{eq:locZ}) and its local time flow is $\Psi({\mathcal L}_{\mathcal X},f,c)$. 
\end{proposition}

\noindent For sake of brevity, we will write in the rest of this section ${\mathcal Z}:=\Upsilon({\mathcal X},f,c)$. 

\bigskip

{\bf Remarks}. 
(i) The definition of $\Upsilon$ does not depend on the choice of the interval $I$ on which is defined $f$. \\
(ii) Observe that $|C_f(t)|<\infty$ for all $t\in J$ and $t\to C_f(t)$ is strictly increasing so that ${\mathcal Z}$ is well defined and is a continuous process. Moreover, for all $t\in C_f(J)$, ${\mathcal Z}_t \in \eta_{f,c}(I)$. \\
(iii)
In the course of the proof, we will prove that the local time $L_{\mathcal Z}(t,x)$ of ${\mathcal Z}$ at time $t\in C_f(J)$ and position $x\in \eta_{f,c}(I)$ is given by 
\begin{equation}\label{eq:locZ}
L_{\mathcal Z}(t,x) = { L_{{\mathcal X}}(C_f^{-1}(t),\eta_{f,c}^{-1}(x)) \over f\circ \eta_{f,c}^{-1}(x) }.\end{equation}
In particular, the total local time of ${\mathcal Z}$ at position $x$ is $\Psi(g,f,c)(x)\le 1$ for all $x\in \eta_{f,c}(I)$. \\
(iv) Explicitly, the proposition means that 
$$
L_{{\mathcal Z}}(\tau^{\hat r}_v({\mathcal Z}), {\hat x})=\Psi({\mathcal L}_{\mathcal X},f,c)_{{\hat r},{\hat x}}(v)
$$

\noindent for all ${\hat r},{\hat x}\in {\eta_{f,c}(I)}$ and $v\in [0, \Psi(g,f,c)(\hat r)]$ where $\tau^{\hat r}_v({\mathcal Z}):= \inf\{t \in C_f(J) \,: \,L_{\mathcal Z}(t,{\hat r})>v\}$. \\

\bigskip

\noindent {\it Proof of the proposition}. For simplicity, we suppose that $c=0$. First we prove \eqref{eq:locZ} by following the proof of Lemma 2 of \cite{warren-yor}. Let $h: \r \to \r_+$ be a Borel function. Let $t_1:=\inf J$ and $t_2:=\sup J$. Using change of variables, we have for $t \in (t_1,t_2)$,
\begin{eqnarray*}
\int_{C_f(t_1)}^{C_f(t)} h({\mathcal Z}_u) \d u 
&=&
\int_{t_1}^t h({\mathcal Z}_{C_f(s)} ) \d C_f(s)
\\
%&=&
%\int_{t_1}^t h({\mathcal Z}_{C_f(s)}) {\d s \over f({\mathcal X}_s)^2} \\
&=&
\int_{t_1}^t h \circ \eta_f({\mathcal X}_s ) {\d s \over f({\mathcal X}_s )^2} \\
&=&
\int_{I_f}   h\circ \eta_f(x) {L_{\mathcal X}(t,x) \over f(x )^2} \d x
\end{eqnarray*}

\noindent by the occupation times formula. Let $\widehat L(C_f(t),\eta_f(x)) := { L_{{\mathcal X}}(t,x) \over f(x )}$.  We get by the change of variables $z=\eta_f(x)$,
$$
\int_{C_f(t_1)}^{C_f(t)} h({\mathcal Z}_u) \d u 
=
\int_{-\infty}^\infty  h(z)\widehat L(C_f(t),z) \d z.
$$

\noindent Therefore ${\mathcal Z}$ has local time given by $\widehat L$. It proves \eqref{eq:locZ}. Let $r,x\in I $ and  $0\le v\le \frac{g(r)}{f(r)}$. Set  $\hat r:= \eta_f(r)$, $\hat x := \eta_f(x)$, $\hat \tau_v^{\hat r}  := C_f(\hat t)$ where $\hat t:= \tau_{v f(r)}^r ({\mathcal X})$. Then $\widehat L(\hat \tau_v^{\hat r}  , \hat r) = v$, and  since $v\to \hat \tau_v^{\hat r}$ is right-continuous, we deduce that $\hat \tau$ is the inverse local time of ${\mathcal Z}$. Moreover, by definition,
$$
\widehat L(\hat \tau_v^{\hat r}, \hat x)
=
{ L_{\mathcal X}(\tau_{v f(r)}^r({\mathcal X}),x) \over f(x)} = {  {\mathcal S}_{r,x} (vf(r)) \over f(x)}.
$$

\noindent  By Definition \ref{def:imageflow}, it is $\Psi({\mathcal L}_{\mathcal X},f)_{{\hat r},{\hat x}}(v)$ indeed. $\Box$

\bigskip

The transformation $\Upsilon$ behaves well under composition. It is the content of the following lemma. For a process $X$, we let $X_{s+\cdot}$ denote the process $ X_{s+t}$, $t\ge 0$ (as long as $X_{s+t}$ is well-defined).

\begin{lemma}\label{l:comp}
Let $s\in J $ and $c'\in I$. We suppose that $f-L_{\mathcal X}(s,\cdot)$ is positive on an interval containing $c'$ and $\{{\mathcal X}_{s+\cdot}\}$.
%\begin{enumerate}[(i)]
%\item $x\to f(x)-L_{\mathcal X}(s,x)$ is a good function,
%\item 
%$f(x)-L_{\mathcal X}(s,x)>0$ for  $x=c'$ or $x\in \{{\mathcal X}_{s+\cdot}\}$.
%\item $\int_s^{s+t} { {\rm d}u \over (f({\mathcal X}_u)-L_{\mathcal X}(s,{\mathcal X}_u))^2}<\infty$ for all $t\ge 0$ such that $s+t\in J$.
%\end{enumerate} 
Then, $\Upsilon\left({\mathcal Z}_{C_f(s)+\cdot}, 1-L_{\mathcal Z}(C_f(s),\cdot),\eta_{f,c}(c')\right)$ is well-defined and equals $ \Upsilon\left({\mathcal X}_{s+\cdot}, f-L_{\mathcal X}(s,\cdot),c'\right)$. In particular, the values at time $0$ of the two processes are identical, i.e.
\begin{equation}\label{eq:XZ0}
\int_{\eta_{f,c}(c')}^{\mathcal Z_{C_f(s)}} {\d r\over 1-L_{\mathcal Z}(C_f(s),r) } = \int_{c'}^{{\mathcal X}_s} {\d r \over f(r)-L_{\mathcal X}(s,r)}.
\end{equation}
\end{lemma}

\bigskip

\noindent {\bf Remark}. By saying that $\Upsilon\left({\mathcal Z}_{C_f(s)+\cdot}, 1-L_{\mathcal Z}(C_f(s),\cdot),\eta_{f,c}(c')\right)$ is well-defined, we mean that
\begin{enumerate}[(i)]
\item the total local time of ${\mathcal Z}_{C_f(s)+\cdot}$ at position $u$ is smaller than $1-L_{\mathcal Z}(C_f(s),u)$ (which is clear since the total local time of ${\mathcal Z}$ is smaller than $1$)
%\item $u\to 1-L_{\mathcal Z}(C_f(s),u)$ is a good function,
\item $1- L_{\mathcal Z}(C_f(s),\cdot) >0$ on an interval containing   $\eta_{f,c}(c')$ and $\{{\mathcal Z}_{C_f(s)+\cdot}\}$.
%$\eta_{f,c}(c')$ and the range of $\mathcal Z_{C_f(s)+\cdot}$ are contained in $({\mathfrak g}_{1-L_{\mathcal Z}(C_f(s),\cdot)},{\mathfrak d}_{1-L_{\mathcal Z}(C_f(s),\cdot)})$,
%\item $\int_{C_f(s)}^{C_f(s+t)} { {\rm d}u \over (1-L_{\mathcal Z}(C_f(s),{\mathcal Z}_u))^2}<\infty$ for all $t\ge 0$ such that $s+t\in J$.
\end{enumerate}

\bigskip

\noindent {\it Proof of  Lemma \ref{l:comp}}.  
We first check that $\Upsilon\left({\mathcal Z}_{C_f(s)+\cdot}, 1-L_{\mathcal Z}(C_f(s),\cdot),\eta_{f,c}(c')\right)$ is well-defined, hence we check statement (ii) above since (i) is clear. 
By  \eqref{eq:locZ}, for any $u\in \eta_{f,c}(I)$,
$$
L_{\mathcal Z}(C_f(s),u) = {L_{\mathcal X}(s,\eta_{f,c}^{-1}(u)) \over f(\eta_{f,c}^{-1}(u))}.
$$

\noindent  We deduce that (ii) is satisfied since $\{{\mathcal Z}_{C_f(s)+\cdot}\}=\eta_{f,c}(\{{\mathcal X}_{s+\cdot}\})$. Let us prove the second statement of the lemma. Let $t\ge 0$ with $s+t\in J$. Let 
$$
t'':= \int_s^{s+t} { {\rm d}u \over (f({\mathcal X}_u)-L_{\mathcal X}(s,{\mathcal X}_u))^2}, \qquad c'':=\eta_{f,c}(c').
$$
 
 \noindent Substituting $v$ for $C_f(u)$ in the above integral and using that ${\mathcal Z}_{C_f(u)} = \eta_{f,c}({\mathcal X}_u)$ by \eqref{def:Z}, we get 
\begin{equation}\label{eq:Cfst}
t''=
\int_{C_f(s)}^{C_f(s+t)} { \d v \over (1-L_{\mathcal Z}(C_f(s),{\mathcal Z}_v))^2}.
\end{equation}

\noindent Observe that by the change of variables $x=\eta_{f,c}^{-1}(r)$ in the integral,
\begin{equation}\label{eq:etaZ}
\eta_{1-L_{\mathcal Z}(C_f(s),\cdot),c''}(u) = \int_{c''}^{u} { {\rm d} r \over 1 - L_{\mathcal Z}(C_f(s),r)}
=
\int_{c'}^{\eta_{f,c}^{-1}(u)} { {\rm d} x \over f(x) - L_{\mathcal X}(s,x)}.
\end{equation}

%\noindent Together with \eqref{eq:gdZ}, and since $f-L_{\mathcal X}(s,\cdot)$ is a good function by assumption, we deduce that $1-L_{\mathcal Z}(C_f(s),\cdot)$ is also a good function. Moreover, the range of ${\mathcal Z}_{C_f(s)+\cdot}$, which is the image of the range of ${\mathcal X}_{s+\cdot}$ by $\eta_{f,c}$, is contained in $({\mathfrak g}_{1-L_{\mathcal Z}(C_f(s),\cdot)},{\mathfrak d}_{1-L_{\mathcal Z}(C_f(s),\cdot)})$ by \eqref{eq:gdZ}, and $c''$ is also included by the assumption on $c'$. Finally, 

%\noindent which shows in particular that 
%$$
%\int_{C_f(s)}^{C_f(s+t)} { \d v \over (1-L_{\mathcal Z}(C_f(s),{\mathcal Z}_v))^2}<\infty.
%$$

%\noindent It proves that $\Upsilon({\mathcal Z}_{C_f(s)+\cdot}, 1-L_{\mathcal Z}(C_f(s),\cdot),c'')$ is well-defined. 

 Let us  consider $\Upsilon({\mathcal Z}_{C_f(s)+\cdot}, 1-L_{\mathcal Z}(C_f(s),\cdot),c'')$ at time $t''$.  By definition and using equation \eqref{eq:Cfst},  
$$
\Upsilon({\mathcal Z}_{C_f(s)+\cdot}, 1-L_{\mathcal Z}(C_f(s),\cdot),c'')(t'')= \eta_{1-L_{\mathcal Z}(C_f(s),\cdot),c''}({\mathcal Z}_{C_f(s+t)}).
$$

 \noindent Equation \eqref{eq:etaZ} with $u={\mathcal Z}_{C_f(s+t)} = \eta_{f,c}({\mathcal X}_{s+t})$ yields
$$
\eta_{1-L_{\mathcal Z}(C_f(s),\cdot),c''}({\mathcal Z}_{C_f(s+t)}) = \int_{c'}^{{\mathcal X}_{s+t}} { {\rm d} x \over f(x) - L_{\mathcal X}(s,x)}= \eta_{f-L_{\mathcal X}(s,\cdot),c'}({\mathcal X}_{s+t}).
$$

\noindent By definition, it is the value of the process $\Upsilon({\mathcal X}_{s+\cdot}, f-L_{\mathcal X}(s,\cdot),c')$ at time $t''$. The proof is complete. $\Box$

\bigskip

\begin{notation}\label{n:upsilon}
When $c={\mathcal X}_0$, we will write $\Upsilon(\mathcal X,f)$ for $\Upsilon(\mathcal X,f,{\mathcal X}_0)$. 

(i) Direct computations also show that  $\Upsilon(\mathcal X,f)$ stays unchanged when replacing $\mathcal X$ by $\mathcal X+a$ and $f$ by $f(\cdot-a)$ for $a\in\r$ arbitrary.

(ii) Since $\eta_{f,c}({\mathcal X}_s)={\mathcal Z}_{C_f(s)}$ by definition, Lemma \ref{l:comp} in the case that $c'={\mathcal X}_s$ reads 
$$  
\Upsilon\left({\mathcal Z}_{C_f(s)+\cdot}, 1-L_{\mathcal Z}(C_f(s),\cdot)\right) 
=
\Upsilon\left({\mathcal X}_{s+\cdot}, f-L_{\mathcal X}(s,\cdot)\right)
$$

\end{notation}

\bigskip

  We finish this section by a continuity lemma on the transformation $\Upsilon$ that will be used later on. 

\begin{lemma}\label{l:continuity}
We fix a couple $(\mathcal X,f)$ where $\mathcal X: J=[0,\infty)\to \r$, $f:I\to (0,\infty)$ are continuous and $\mathcal X_t \in I$ for all $t\ge 0$.

We consider a sequence $(\mathcal X^n,f^n)_{n\ge 1}$ such that $\mathcal X^n:[0,j_n)\to \r$, $j_n\in \r_+\cup\{+\infty\}$, $f^n:I_n\to (0,\infty)$ and $\mathcal X_t \in I_n$ for all $t\in [0,j_n)$. The sets $I$, $I_n$ are intervals.  

We make the following assumptions:
\begin{enumerate}[(i)]
\item for all $n$, $\mathcal X^n_0=\mathcal X_0=0$  and $\int_0^{j_n} {\d s \over f^n(\mathcal X^n_s)^2}=\int_0^\infty {\d s \over f(\mathcal X_s)^2}=\infty$;
\item  $\lim_{n\to\infty} j_n=\infty$;
\item $\mathcal X^n$ converges uniformly to $\mathcal X$ on any compact set of $\r_+$;
\item for any compact set $K$ of $\r_+$,  there exists a compact set $L$ in $\r$ such that
\begin{enumerate}
\item for $n$ large enough, $\mathcal X(K)$ and $\mathcal X^n(K)$ are contained in the set $L$ which is itself contained in $I$ and in $I_n$;
\item $f^n$ converges uniformly to $f$ on $L$.
\end{enumerate}
\end{enumerate}
Then $\Upsilon({\mathcal X}^n,f^n)$ converges to $\Upsilon({\mathcal X},f)$ (for the topology of uniform convergence on compact sets of $\r_+$).
\end{lemma}

See Appendix \ref{s:continuitylemmas} for the proof.

\section{Jacobi flows}
\label{s:perkins}

In this section, we define Jacobi flows, and show how they arise by disintegration of the BESQ flows. 

Let $\lambda$ denote the Lebesgue measure.  \begin{definition}\label{d:M}
We introduce the covariance functional
$$
Q_h(A,\widehat A) =  h\left(\lambda(A\cap \widehat A) -\lambda(A\cap[0,1])\lambda(\widehat A\cap[0,1])\right)
$$

\noindent for $A$, $\widehat A$ Borel sets of $\r_+$, and $h\ge 0$. Let $({\mathcal E}_s,\,-\infty<s<\infty)$ be a right-continuous filtration. We consider  a collection of random variables $({\mathcal M}(A\times[s,t]),\, s\le t,\, A \hbox{ Borel set of }\r_+)$ such that for any $s\in \r$, $A$, $\widehat  A$ Borel sets of $\r_+$, $({\mathcal M}(A\times[s, s+h]))_{h\ge 0}$ and $({\mathcal M}(\widehat  A\times[s, s+h]))_{h\ge 0}$ are two continuous martingales with respect to the filtration $({\mathcal E}_{s+h},\, h\ge 0)$
%$\sigma\{{\mathcal M}(\cdot\times [s, s+u]), 0\le u \le h\}, h\ge0, $ 
and 
\begin{equation}\label{eq:bracket}
\langle {\mathcal M}(A\times[s, s+\cdot]), {\mathcal M}(\widehat  A\times[s, s+\cdot])\rangle_h= Q_h(A, \widehat  A).
\end{equation} 
\end{definition}

From \eqref{eq:bracket}, we deduce that for fixed $s$, $h$, and $(A_n)_{n}$  disjoint Borel sets with $\lambda(\bigcup_n A_n)<\infty$, we have  ${\mathcal M}(\bigcup_n A_n,[s,s+h])=\sum_n {\mathcal M}(A_n,[s,s+h])$ in $L^2(\Omega)$. In other words, ${\mathcal M}$ is an $({\mathcal E}_s)$-martingale measure with covariance functional $Q$ on $\r_+\times \r$ in the sense of \cite{walsh84}.

\bigskip

 The martingale measure ${\mathcal M}$ can be defined through a white noise ${\mathcal W}$ on $\r_+\times \r$:  \begin{equation} \label{def-mathcalM}
 {\mathcal M}(A\times[s, t]):= {\mathcal W}(A \times [s, t])- \lambda(A\cap [0,1]) {\mathcal W}([0, 1] \times [s, t]).
 \end{equation}

 We now introduce the Jacobi flows.

\begin{definition}\label{d:jacobi}
Let $\delta>0$ and $\delta' \in \r$. Let ${\mathcal Y}$ be a collection of continuous processes $({\mathcal Y}_{s,t}(v),\,t\ge s)_{s\in\r,\, v\in [0,1]}$.  We say that ${\mathcal Y}$ is a (non-killed) {\rm Jacobi($\delta,\delta'$)} flow if it is solution of
\begin{equation}\label{def:Yflow}
 {\mathcal Y}_{s, t}(v)
=
v +  2 \int_s^t  {\mathcal M}\Big( [0, {\mathcal Y}_{s,r}(v)], \d r\Big) 
 + \int_s^t  \delta (1- {\mathcal Y}_{s, r}(v)) - \delta' {\mathcal Y}_{s,r}(v)  \; \d r ,\, t\in [0, T_1({\mathcal Y}_{s, \cdot}(v))]
\end{equation}

\noindent where $\mathcal M$ is a martingale measure on $\r_+\times \r$ with covariance functional $Q$ and $T_1({\mathcal Y}_{s, \cdot}(v)):=\inf\{t\ge s\,:\,{\mathcal Y}_{s, t}(v)=1 \}$ (it may be infinite). We require  the regularity conditions: almost surely,
\begin{enumerate}[(i)]
\item for all $s\in \r$ and $v\in [0,1]$, $ {\mathcal Y}_{s,s}(v)=v$,% and $\lim_{s'\uparrow s} {\mathcal Y}_{s',s}(v)=v$, 
\item for all $s\le t$, $v \in [0,1) \mapsto  {\mathcal Y}_{s,t}(v)$ is c\`adl\`ag,
\item  for all $s'\le s$ and all $v',v\in [0,1]$, if $ {\mathcal Y}_{s',s}(v')> v$, then $ {\mathcal Y}_{s',t}(v')\ge {\mathcal Y}_{s,t}(v)$ for all $t\ge s$,
\item \label{item:killatone} for all $s\le t\le t'$, and $v\in [0,1]$, if ${\mathcal Y}_{s, t}(v)=1$ then ${\mathcal Y}_{s,t'}(1)=1$.
\end{enumerate}
\end{definition}

\begin{definition}\label{d:jacobi-bis}
Let $\delta<2$ and $\delta'\in \r$. We say that ${\mathcal Y}$ is a killed  Jacobi($\delta,\delta'$) flow  if it is solution of \eqref{def:Yflow} but is absorbed at $0$ and satisfies the same regularity conditions. For $\delta \in \r$, we call general {\rm Jacobi($\delta,\delta'$)} flow a Jacobi flow which is either killed or non-killed.
\end{definition}

{\bf Remark}. Being non-killed means that the flow line is reflected at zero. With condition \eqref{item:killatone}, we arbitrarily decided to kill flow lines when they hit $1$. We could also let them reflect at $1$. Therefore, strictly speaking, our Jacobi flows are actually Jacobi flows killed at $1$.

\bigskip

 The existence of Jacobi flows will   be obtained as a consequence of Theorem \ref{t:perkins}. As for the case of BESQ flows, the Jacobi flow can be constructed from a countable number of flow lines.

\begin{definition}\label{d:BESQdd'}
Let $\delta,\delta'\in \r$ and $b \ge 0$. 
Let ${\mathcal W}$ be a white noise on $\r_+^2$. Consider the ${\rm BESQ}^\delta$ flow $({\mathcal S}_{r,x}(a),x\ge r)_{r\in \r,a\ge 0}$  driven by ${\mathcal W}$ as in Definition \ref{d:BESQflow} or Definition \ref{def:BESQflowgeneral}. Let $f=(f(x),\, x \in [0,{\mathfrak d}_f))$ be the pathwise unique solution of  
\begin{equation}\label{def:f1}
f(x) = 
b + 2 \int_0^x  {\mathcal W}([0,f(s)], \d s) + 
 ( \delta + \delta') x
\end{equation}

\noindent where the process is absorbed at $0$ if it hits $0$ at a positive time (we denote by ${\mathfrak d}_f$ this absorption time, possibly infinite).  We call ${\rm BESQ}^{\delta,\delta'}_b$ flow driven by ${\mathcal W}$ the collection of continuous processes $\Theta=(\Theta_{r,x}(a),\, 0\le r \le x <{\mathfrak d}_f),\, a \in [0,f(r)])$  where  $\Theta_{r,\cdot}(a)={\mathcal S}_{r,\cdot}(a)$  until it meets $f$, and is equal to $f$ afterwards.

The ${\rm BESQ}^{\delta,\delta'}_b$ flow is called killed when ${\mathcal S}$ is killed.

When $b=0$, by convention the flow will be defined for $x\ge r>0$.
\end{definition}

\bigskip

{\bf Remark}.
(i) To avoid trivial situations, we suppose either that $b>0$ or $b=0$ and $\delta+\delta'\ge 2$ (otherwise ${\mathfrak d}_f=0$). 

(ii) The flow lines may hit $f$ only in the case $\delta'< 2$.

 \bigskip

We suppose now that $b>0$. We have ${\mathfrak d}_f<\infty$ if and only if $\delta+\delta'<2$. We will use the transformation $\Psi$ of Section \ref{s:flow} with $f$.  We recall that $\eta_f(x):= \int_0^x \frac{\d r}{f(r)}$.  In all cases, $\lim_{x\uparrow {\mathfrak d}_f} \eta_f(x)=+\infty$ (see Lemma \ref{l:eta}).  Define the martingale measure $\widetilde {\mathcal W}$ by, for any Borel set $A\subset \r_+$ with finite Lebesgue measure $\lambda(A)$,
\begin{equation}
\widetilde {\mathcal W}(A \times [0,t]) := \int_{r=0}^{\eta_f^{-1}(t)} {1\over f(r)} {\mathcal W}\left(A f(r),\, \d r \right) , \qquad t\ge 0, \label{def:tildeW}
\end{equation}

\noindent where $Af(r) :=\{af(r),\,a \in A\}$.  If $({\mathscr E}_x,\, x\ge 0)$ denotes the natural filtration of ${\mathcal W}$, with ${\mathscr E}_x:= \sigma({\mathcal W}(\cdot \times [0, r]), 0\le r \le x)$, then
 the process $(\widetilde {\mathcal W}(A \times [0,t]))_{t\ge0} $ is a continuous martingale with respect to the filtration $(\widetilde {\mathscr E}_t,\, t\ge 0):=({\mathscr E}_{\eta_f^{-1}(t)},\, t\ge 0)$ whose quadratic variation process is equal to $\lambda(A) t$.  We deduce that $\widetilde {\mathcal W}$ is a white noise on $\r_+^2$ with respect to the filtration $(\widetilde {\mathscr E}_t,\, t\ge 0)$.

We define the martingale measure ${\mathcal M}^+$ on $\r_+^2$ by the identity for any Borel set $A\subset \r_+$:
\begin{equation}\label{def:M}
 {\mathcal M}^+ ( A \times [0,t])
:=
  \widetilde {\mathcal W}\left(A \times  [0,t] \right) - \lambda(A\cap[0,1]) \,\widetilde {\mathcal W}\left([0,1]\times [0,t]\right), \qquad t\ge 0.
\end{equation}

\noindent Notice that $({\mathcal M}^+(A\times [0,t]))_{t,A}$ is a Gaussian process. It  is actually an $(\widetilde {\mathscr E}_t)$-martingale measure with covariance functional $Q$ as defined in Definition \ref{d:M}, from $s=0$.  From the definition of ${\mathcal M}^+$, we see that for any deterministic Borel set $A\subset \r_+$, and $t\ge 0$,  
\begin{equation}\label{eq:Meta}
{\mathcal M}^+ (A\times [0,t])
=
\int_{r=0}^{\eta_f^{-1}(t)} {1\over f(r)} \Big({\mathcal W}\left(A f(r),\, \d r \right)- \lambda(A\cap[0,1]) \, {\mathcal W}\left([0, f(r)],\, \d r \right)\Big). 
\end{equation}

\noindent We can rewrite it as
\begin{eqnarray}\label{eq:Mg}
&& \int_{r=0}^{\infty}\int_{u=0}^\infty g(u,r) {\mathcal M}^+(\d u, \d r) \\
&=&
\int_{r=0}^{\infty} \int_{u=0}^{\infty} \frac{1}{f(r)}g\Big(\frac{u}{f(r)},\eta_f(r)\Big) \left( {\mathcal W}(\d u,\d r) - {\bf 1}_{[0,f(r)]}(u)\frac{\d u}{f(r)} {\mathcal W}\left([0, f(r)],\, \d r \right)\right)\nonumber
\end{eqnarray}

\noindent for $g={\bf 1}_{A\times [0,t]}$. By definition of the stochastic integral, it is also true when $g= Z {\bf 1}_{A\times (s,t]}$ where $0\le s\le t$ and $Z$ is a bounded $ \widetilde {\mathscr E}_s$-measurable function. By linearity and density, \eqref{eq:Mg} is true for any $(\widetilde {\mathscr E}_t)$-predictable square integrable function  $g$ with respect to the covariance of ${\mathcal M}^+$, i.e. 
$$
\e\left[\int_{t\ge 0}\int_{u\ge 0} g(u,t)^2{\rm d} u {\rm d} t - \int_{t\ge 0}\left(\int_{u\in[0,1]} g(u,t) {\rm d} u\right)^2{\rm d} t \right]<\infty . 
$$

Define 
\begin{equation}\label{def:gamma}
\gamma_x:= \int_0^x{1\over \sqrt{f(r)}}{\mathcal W}([0,f(r)],\d r), \qquad  0\le x< {\mathfrak d}_f.
\end{equation} 

\noindent The process $\gamma$ is an $({\mathscr E}_x)$-Brownian motion stopped at ${\mathfrak d}_f$ when ${\mathfrak d}_f < \infty$.

\begin{proposition}\label{p:M}
The martingale measure ${\mathcal M}^+$ is independent of the process $f$. 
\end{proposition}

\noindent {\it Proof}. The proof is an extension of that of Proposition 8 of \cite{warren-yor} to the setting of martingale measures. Let $\gamma$ be defined by \eqref{def:gamma}. From \eqref{eq:Meta}, we see that $({\mathcal M}^+(A\times[0,\eta_f(x)]),\, 0\le x < {\mathfrak d}_f)$ is a martingale   which is orthogonal to $\gamma$. On the other hand, its increasing process is $x\to (\lambda(A)-\lambda(A\cap[0,1])^2)\eta_f(x)$. We deduce by Knight's theorem that $({\mathcal M}^+(A\times[0,t]),\, t\ge 0)$ is a Brownian motion (with multiplicative constant) which is independent of $\gamma$. More generally, the multidimensional Knight's theorem implies that $({\mathcal M}^+(A_1\times[0,t]),\ldots, {\mathcal M}^+(A_n\times[0,t]) \, t\ge 0)$ are independent of $\gamma$ for any $n\ge 1$ and  disjoint Borel sets $A_1,\ldots,A_n$ of $\r_+$. We deduce that $({\mathcal M}^+(A\times[0,t]))_{t,A}$ is independent of $\gamma$. We observe that the process $f$ is measurable with respect to the filtration of $\gamma$ (it is the strong solution of $\d f(x)= 2 \sqrt{f(x)}\d \gamma_x + (\delta+\delta') \d x$), therefore $({\mathcal M}^+(A\times[0,t]))_{t,A} $ is independent of $f$.  $\Box$

\bigskip

In the notation of Section \ref{s:flow} with $I=[0,{\mathfrak d}_f)$ and $c=0$, we define the flow ${\mathcal Y}^+= ({\mathcal Y}_{s, t}(v), 0 \le s \le t, v \in [0,1])$ by  \begin{equation} \label{def:Y}
{\mathcal Y}^+:=\Psi(\Theta,f)
\end{equation}

\begin{figure}
\begin{minipage}{0.4\textwidth}
\begin{tikzpicture}[line cap=round,line join=round,>=triangle 45,x=.3cm,y=.5cm]
\clip(-10,-2.) rectangle (11.4,6.88);
\node at (-6,-1.5) {$0$};
\node at (-7.,0.6) {$v$};
\node at (-6.7, 3.) {$1$};
\draw[->] (-8,-1.) -- (11.,-1.);
\draw[line width=2pt]  (-6.,3.)-- (-4.9,3.3)-- (-3.94,2.6)-- (-3.,3.)-- (-2.,2.)-- (-0.96,2.76)-- (0.02,2.2)-- (0.96,2.74)-- (1.28,2.34)-- (1.86,2.98)-- (2.22,2.3)-- (2.72,2.66)-- (3.68,1.98)-- (4.78,2.16)-- (5.5,1.48)-- (6.18,2.16)-- (7.06,2.16)-- (7.22,3.02)-- (7.9,2.84)-- (8.38,3.4)-- (9.5,2.64)-- (10.06,3.24)-- (10.76,2.62);
\draw[<-]  (-6.,4.5)-- (-6,-1.);
\draw  (-6.,0.6)-- (-5.54,0.78)-- (-4.92,0.66)-- (-4.,1.28)-- (-3.,1.)-- (-2.4,0.76)-- (-1.74,1.18)-- (-0.54,0.7)-- (-0.06,1.26)-- (1.,1.)-- (1.6,1.5)-- (2.36,1.58)-- (2.62,1.84)-- (3.2,1.58)-- (3.68,1.98);
\end{tikzpicture}
\end{minipage}
\begin{minipage}{0.4\textwidth}
\begin{tikzpicture}[line cap=round,line join=round,>=triangle 45,x=.3cm,y=.5cm]
\clip(-10,-2.) rectangle (11.4,6.88);
\node at (-6,-1.5) {$0$};
\node at (-7.,0.6) {$v$};
\node at (-6.7, 3.) {$1$};
\draw[->] (-8,-1.) -- (11.,-1.);
\draw[line width=2pt]  (-6,3.) -- (11.,3.);
\draw[<-]  (-6.,4.5)-- (-6,-1.);
\draw  (-6.,0.6)-- (-5.5,0.5)-- (-4.9,1.)-- (-4.,0.8)-- (-3.,1.)-- (-2.4,0.6)-- (-2,0.4)-- (-1.5,1.)-- (-1.,1.7)-- (-0.5,1.5)-- (0.,2.3)-- (0.5,1.8)-- (1.,1.5)-- (1.5,1.)-- (2.,1.7)--(2.5,1.8)--(3.,2.2)--(3.5,1.9)--(4.,2.3)--(4.5,2.7)--(5.,2.3)--(5.5,2.8)--(6.,3.);
\end{tikzpicture}
\end{minipage}
\caption{%\leftskip=1.8truecm \rightskip=1.8truecm 
\small  On the left: ${\rm BESQ}^{\delta,\delta'}$ flow. The process $x\to f(x)$ is traced in bold. A flow line $x\to S_{0,x}(v)$ is also represented. On the right: Jacobi($\delta,\delta'$) flow after a space-time transformation. The process $x\to f(x)$ is mapped to the constant function $1$, and the image of the flow line at $v$ is a flow line of the Jacobi flow.}
\label{f:perkins0}
\end{figure}
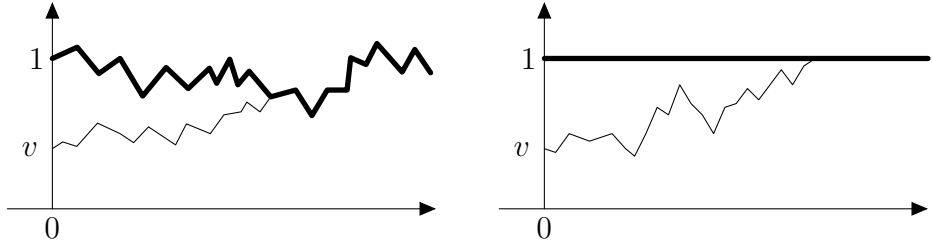

\noindent  (we add the superscript in ${\mathcal Y}^+$ to stress that the flow is only defined for $0\le s\le t$ whereas the flow ${\mathcal Y}$ is defined for all $-\infty< s\le t$), see Figure \ref{f:perkins0}. By definition,
\begin{equation}\label{def:Y2}
{\mathcal Y}_{s,t}(v) = { {\Theta }_{\eta_f^{-1}(s),\eta_{f}^{-1}(t)}\Big(v f\circ \eta_f^{-1}(s)\Big)\over f\circ \eta_{f}^{-1}(t)}, \qquad t\ge s\ge 0,\,  v \in [0,1].
\end{equation}

\begin{theorem}\label{t:perkins}
Let $\delta,\delta'\in \r$. Let $b>0$ and $\Theta $ be a ${\rm BESQ}^{\delta,\delta'}_b$ flow. Then, the flow ${\mathcal Y}^+$ is a {\rm Jacobi($\delta,\delta'$)} flow, restricted to the positive time-axis. It is the $[0,1]$-valued pathwise unique solution, which is strong, of the following SDE.  Let $v\in [0,1]$ and $s\ge 0$.  Almost surely, for $t\ge s$,  
\begin{equation}\label{eq:Y1}
 {\mathcal Y}_{s, t}(v)
=
v +  2 \int_s^t  {\mathcal M}^+( [0, {\mathcal Y}_{s,r}(v)], \d r) 
 +  \int_s^t  \left(\delta(1- {\mathcal Y}_{s, r}(v)) - \delta'{\mathcal Y}_{s, r}(v) \right) \d r  
\end{equation}

\noindent with the required absorption conditions\footnote{The process is absorbed if it hits $1$. If $\Theta$ is killed, then the process is also absorbed if it hits $0$.}.   In particular, ${\mathcal Y}^+$ is independent of $f$. 
\end{theorem}

{\bf Remark}. The theorem will be used only in the case $\delta'=0$ in the rest of the paper. The general case will be used in a follow-up paper.

\bigskip

\noindent {\it Proof}.  The regularity conditions (i), (iii) and (iv) of ${\mathcal Y}^+$ come from the properties of $\Theta$. Let us check that (ii) is also satisfied. Let $u\in [0,1)$ and $s\le t$. If ${\mathcal Y}_{s,t}(u)=1$, then ${\mathcal Y}_{s,t}(v)=1$ for all $v\ge u$ hence $v\to {\mathcal Y}_{s,t}(u)$ is c\`adl\`ag at $u$. Suppose then that ${\mathcal Y}_{s,t}(u)<1$. There exists $v>u$ such that ${\mathcal Y}_{s,t}(v)<1$. Indeed, the  flow $({\mathcal S}_{r,x}(a),\, x\ge r)_{r\in \r,\, a\ge 0}$ verifies: almost surely, for all $r<x$ and $a\ge 0$, $\sup_{r'\in [r,x]}|{\mathcal S}_{r,r'}(b)-{\mathcal S}_{r,r'}(a)|\to 0$ as $b\downarrow a$. This statement follows from the identification ${\mathcal S}_{r,x}(a)=L(\tau_a^r,x)$ (see Section \ref{s:BESQflow}), the fact that $b\to \tau_b^{r}$ is c\`adl\`ag, and the bicontinuity of the local times.  Since for $v>u$ close enough to $u$, one has ${\mathcal Y}_{s,t}(v)<1$, we can use the regularity of the flow $({\mathcal S}_{r,x}(a),\, x\ge r)_{r\in \r,\, a\ge 0}$ to see that $v\to {\mathcal Y}_{s,t}(v)$ is c\`adl\`ag at $u$.

We show now the SDE. To avoid too much notation, we prove it for $s=0$. For brevity, we suppose that $\Theta$ is non-killed. The other case is similar, one just has to look at the times before the hitting time of $0$. By  \eqref{def:Y2}, we have
$$
{\mathcal Y}_{0,t}(v) = {{\Theta }_{0,\eta_f^{-1}(t)} (bv) \over f\circ \eta_f^{-1}(t)}.
$$

\noindent Equation \eqref{eq:Y1} is equivalent to 
\begin{equation}\label{eq:proofY1a}
{{\Theta }_{0,x} (b v) \over f(x)} = 
v +  2 \int_0^{\eta_f(x)}  {\mathcal M}^+( [0, {\mathcal Y}_{0,r}(v)], \d r) 
 +  \int_0^{\eta_f(x)} \left( \delta (1- {\mathcal Y}_{0, r}(v))   -\delta'{\mathcal Y}_{0, r}(v) \right)\d r.
\end{equation}

\noindent Let $\xi_{x}(v):={{\Theta }_{0,x}(b v)\over f(x)}$. By \eqref{eq:Mg} with $g(u,r):={\bf 1}_{[0,{\mathcal Y}_{0,r}(v)]\times [0,\eta_f(x)]}(u,r)$ (hence $g(\frac{u}{f(r)},\eta_f(r))={\bf 1}_{[0,{\Theta }_{0,r}(bv)]\times [0,x]}$), the first integral of the right-hand side of \eqref{eq:proofY1a} is equal to 
\begin{equation}\label{eq:proofY1b}
\int_0^{x}  \frac{1}{f(r)}\left( {\mathcal W}( [0, {\Theta }_{0,r}(bv)] , \d r)  
 -  \xi_{r}(v)    {\mathcal W}([0,f(r)], \d r) \right).
\end{equation}

\noindent By a change of variables, the last term of the right-hand side in \eqref{eq:proofY1a} is
\begin{equation}\label{eq:proofY1c}
\int_0^{\eta_f(x)}  \delta (1- {\mathcal Y}_{0, r}(v)) -\delta' {\mathcal Y}_{0,r}(v) \d r = \int_0^{x}  \frac{\delta(1- \xi_r(v))-\delta' \xi_r(v)}{f(r)}  \d r.
\end{equation}

\noindent Let us verify \eqref{eq:proofY1a}. By \eqref{def:f1}, we have $ \langle f, f\rangle_x = 4  \int_0^x   f(s) \d s$  and 
$$
\langle {\Theta }_{0, \cdot}(b v), f \rangle_x = 4  \int_0^x   {\Theta }_{0, s}(b v) \d s.
$$

\noindent  We deduce from It\^{o}'s formula that, as long as ${\Theta }_{0,x}(b v)< f(x)$,  ($\d \equiv \d_x$ below) 
\begin{eqnarray*}
&& \d \xi_x(v)\\
&=&
\frac{\d {\Theta }_{0,x}(b v)}{f(x)} -  \frac{{\Theta }_{0,x}(b v)}{f(x)^2} \d f(x) + {\Theta }_{0,x}(b v) \frac{4f(x)}{f^3(x)}\d x - \frac{1}{f^2(x)} 4 {\Theta }_{0,x}(bv)\d x \\
&=&
\frac{\d {\Theta }_{0,x}(b v)}{f(x)} -  \frac{{\Theta }_{0,x}(b v)}{f(x)^2} \d f(x)  \\
&=&
\frac{1}{f(x)}  \big(2 {\mathcal W}([0,{\Theta }_{0,x}(b v)], \d x) + \delta \d x\big) 
 - \frac{  \xi_{x}(v)}{f(x)}  \big( 2 {\mathcal W}( [0, f(x)] ,\d x) +( \delta +\delta') \d x\big) 
\\
&=&
\frac{2}{f(x)}\left( {\mathcal W}( [0, {\Theta }_{0,x}(b v)] , \d x)  
 -  \xi_{x}(v)    {\mathcal W}([0,f(x)], \d x) \right)
+ \frac{ \delta (1 -  \xi_{x}(v)) -\delta '\xi_x(v)}{f(x)}  \d x.
\end{eqnarray*}

\noindent Comparing with \eqref{eq:proofY1b} and \eqref{eq:proofY1c}, we get \eqref{eq:proofY1a} indeed.

Now we show that ${\mathcal Y}^+$ is the pathwise unique solution (by Yamada--Watanabe's theorem it will imply that it is a strong solution). This can be achieved by imitating the usual proof of the pathwise uniqueness of a one-dimensional SDE with non-Lipschitz coefficient. We give the details here for completeness. As before we take $s=0$ for notational brevity. %and show that ${\mathcal Y}_{0, t}(v)$ is measurable with respect to the martingale measure ${\mathcal M}^+$.  
  Define $a_0:=1$ and $a_k:= a_{k-1} e^{-k}$ for $k\ge1$. Let $\psi_k$ be a continuous function on $\r$ with support in $(a_k, a_{k-1})$ such that $\int_{a_k}^{a_{k-1}} \psi_k(x) \d x=1$ and $0\le \psi_k(x)  \le \frac{2}{k x}$ for $x\in (a_k,a_{k-1})$. Define $$\phi_k(y):= \int_0^{|y|} \d x \int_0^x \psi_k(r) \d r, \qquad y \in \r.$$

\noindent Observe that $\phi_k$ is twice continuously differentiable and $|\phi'_k|\le 1$. Moreover $\phi_k(y) \to |y|$ as $k\to \infty$.  Suppose that ${\mathcal Y}^{(1)}_{0, t}(v)$ and ${\mathcal Y}^{(2)}_{0, t}(v)$, $t\ge0$, are two solutions of 
\eqref{eq:Y1} [with $s=0$ there] with respect to the same martingale measure ${\mathcal M}^+$. Let $\Delta {\mathcal Y}_t:= {\mathcal Y}^{(1)}_{0, t}(v)- {\mathcal Y}^{(2)}_{0, t}(v)$ for $t\ge0$. By \eqref{def:M}, we have    $\int_0^t  {\mathcal M}^+( [0, {\mathcal Y}_{0,r}^{(i)}(v)], \d r)= \int_0^t \big(\widetilde {\mathcal W}([0,{\mathcal Y}_{0,r}^{(i)}(v)], \d r)- {\mathcal Y}_{0,r}^{(i)}(v) \widetilde {\mathcal W}([0,1], \d r)\big)$ for $i=1,2$. Then $$ \Delta {\mathcal Y}_t
=
2 \int_{r=0}^t \int_{s=0}^1 \big[ ( 1_{\{s\le {\mathcal Y}^{(1)}_{0, r}(v)\}}- {\mathcal Y}^{(1)}_{0, r}(v))- ( 1_{\{s\le {\mathcal Y}^{(2)}_{0, r}(v)\}}- {\mathcal Y}^{(2)}_{0, r}(v))\big] \widetilde {\mathcal W}(\d s, \d r) - (\delta+\delta') \int_0^t \Delta {\mathcal Y}_r \d r.$$

\noindent Since $0\le {\mathcal Y}^{(1)}_{0, r}(v), {\mathcal Y}^{(2)}_{0, r}(v)\le 1$,  $\d \langle \Delta {\mathcal Y} \rangle_r= 4 \d r \int_{s=0}^1 \big[ ( 1_{\{s\le {\mathcal Y}^{(1)}_{0, r}(v)\}}- {\mathcal Y}^{(1)}_{0, r}(v))- ( 1_{\{s\le {\mathcal Y}^{(2)}_{0, r}(v)\}}- {\mathcal Y}^{(2)}_{0, r}(v))\big]^2 \d s $ and we observe that $ \int_{s=0}^1 \big| ( 1_{\{s\le {\mathcal Y}^{(1)}_{0, r}(v)\}}- {\mathcal Y}^{(1)}_{0, r}(v))- ( 1_{\{s\le {\mathcal Y}^{(2)}_{0, r}(v)\}}- {\mathcal Y}^{(2)}_{0, r}(v))\big|^2 \d s  \le 2  |\Delta {\mathcal Y}_r| $.  Applying It\^{o}'s formula to $\phi_k(\Delta {\mathcal Y}_t)$ gives that \begin{eqnarray*} \e(\phi_k(\Delta {\mathcal Y}_t))
&=&
-(\delta+\delta')\, \e \int_0^t \phi_k'(\Delta {\mathcal Y}_r)\,  \Delta {\mathcal Y}_r \d r + \frac12 \e \int_0^t \psi_k(|\Delta {\mathcal Y}_r|) \d \langle \Delta {\mathcal Y} \rangle_r\\
&\le&
|\delta+\delta'| \,  \e \int_0^t |\Delta {\mathcal Y}_r| \d r + \frac{8t}{k},
\end{eqnarray*}

\noindent where in the last inequality we have used the facts that $|\phi_k'|\le 1$ and $\psi_k(x) x \le \frac{2}{k}$ for any $x\ge0$.  Letting $k\to\infty$ we deduce from Fatou's lemma  that $\e(|\Delta {\mathcal Y}_t|) \le |\delta+\delta'| \, \e \int_0^t |\Delta {\mathcal Y}_r| \d r$, yielding that $\Delta {\mathcal Y}_t=0$ by Gronwall's inequality. This shows the pathwise uniqueness and completes the proof of the theorem.  $\Box$

 \bigskip
 
Let $\delta >0$. Recall that $X$ is a PRBM defined in \eqref{def:muprocess} and $S$ is its local time flow introduced in \eqref{def:S}. Recall that $S$ satisfies the SDEs \eqref{eq:S1} and \eqref{eq:S2} driven by $W$.  The following corollary of Theorem \ref{t:perkins} gives a pathwise construction of a Jacobi flow starting from a PRBM.

\begin{corollary} \label{t:Sperkins} Let $b>0$ and $\mu>0$. 
\begin{enumerate}[1)]
\item  Let  $f=(f(x),\, x\ge 0)$   defined by $f(x)=S_{0,x}(b)=L(\tau_{b}^0,x)$. Notice that ${\mathfrak d}_f=\inf\{x\ge 0\,:\, L(\tau_b^0,x)=0\}$. We have ${\mathfrak d}_f<\infty$ if and only if $\mu > 1$. Let $\Theta := (S_{r,x}(a),\, 0\le r\le x < {\mathfrak d}_f,\, a \in [0,f(r)])$. The flow $\Psi(\Theta,f)$ is a {\rm Jacobi($2/\mu,0$)} flow in the positive time-axis, independent of $f$.
\item  Let  $f=(f(x),\, x\ge 0)$   defined by $f(x)=S_{0,-x}(b)=L(\tau_b^0,-x)$. Notice that ${\mathfrak d}_f= - \inf_{[0,\tau_b^0]} X$. Let $\Theta := (S_{-r,-x}(a),\, 0\le r\le x < {\mathfrak d}_f,\, a \in [0,f(r)])$. The flow $\Psi(\Theta,f)$ is a killed {\rm Jacobi($2-2/\mu,0$)} flow in the positive time-axis, independent of $f$. 
%\item Let $f^0(x):=S_{0,x}(0)$ be the local time process of $(X_t,\,t\le 0)$ and ${\mathcal S}^0$ its forward local time flow. Let $c>0$. The flow  $\Psi({\mathcal S}^0,f^0,c)$ is a {\rm Jacobi(${\delta,0}$)} flow.
\end{enumerate}
\end{corollary}

\begin{figure}[h!]
\begin{tikzpicture}[line cap=round,line join=round,>=triangle 45,x=0.7cm,y=0.7cm]
%\draw [xstep=1.0cm,ystep=1.0cm] (-4.3,-4.72) grid (15.04,6.3);
\clip(-6,-1.5) rectangle (15.04,2.7);
\draw[->] (-4,-1)--(15,-1);
\draw (-4,2)--(15,2);
\node at (-4,-0.7) {$0$};
\node at (-4,2.3) {$1$};
\node at (3.7,0.2) {$v$};
\node at (4.,-1.25) {$0$};
%\draw [domain=-4.3:15.04] plot(\x,{(-11.748--0.04*\x)/11.4});
%\draw [domain=-4.3:15.04] plot(\x,{(--22.56-0.*\x)/11.28});
\draw  (4.,2.)-- (4.,-1.);
\draw[->] (4.,0.4) -- (4.48,0.835);
\draw[->] (4.,0.4) -- (3.445,0.715);
\draw (4.48,0.835)-- (4.64,0.98)-- (5.34,0.28)-- (5.72,0.56)-- (6.32,-0.42)-- (6.68,-0.06)-- (7.28,-0.7)-- (7.66,-0.36)-- (8.,-1.)-- (8.4,-0.5)-- (9.,-0.64)-- (9.3,-0.18)-- (10.,-1.)-- (10.4,-0.36)-- (10.72,-0.52)-- (11.16,-0.06)-- (11.36,0.68)-- (11.86,0.42)-- (12.3,0.9)-- (12.68,0.78)-- (13.1,1.36)-- (13.56,1.16)-- (14.,2.);
\draw (3.445,0.715)-- (3.26,0.82)-- (2.62,0.38)-- (2.24,0.74)-- (1.78,0.14)-- (1.14,0.32)-- (0.8,-0.36)-- (0.22,-0.22)-- (-0.16,0.08)-- (-0.56,-0.46)-- (-1.24,-0.26)-- (-1.58,-0.76)-- (-2.,-0.56)-- (-2.34,-1.);
\end{tikzpicture}
\caption{%\leftskip=1.8truecm \rightskip=1.8truecm 
\small The flow $Y$. A forward flow line in the positive time-axis is a Jacobi($2/\mu,0$) process. A backward flow line in the negative time-axis is a Jacobi($2-2/\mu,0$) process absorbed at $0$.}
\label{f:perkins}
\end{figure}
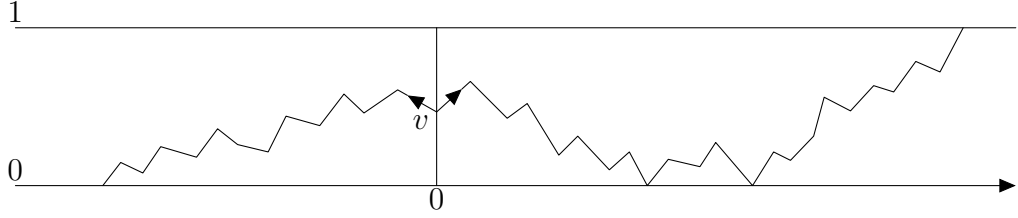

\begin{remark}\label{r:Y} We can rephrase the corollary in the following way.  Let $b>0$. Define the flow $Y= (Y_{s, t}(v), -\infty < s, t<\infty, v \in [0,1])$ by 
$Y:=\Psi({\mathcal L}_X^b,f)$ where $f(x)=S_{0,x}(b)$ and ${\mathcal L}_X^b$ is the local time flow of $(X_t,\,t\le \tau_b^0)$. The forward flow of $Y$ in the positive time-axis (i.e.
$(Y_{s,t}(v), 0\le s\le t,\, v\in [0,1])$)  is a Jacobi($2/\mu,0$) flow. Similarly, the backward flow of $Y$ in the negative time-axis (i.e.
$(Y_{-s,-t}(v), 0\le s\le t,\, v\in [0,1])$) is a killed Jacobi($2-2/\mu,0$) flow. See Figure \ref{f:perkins}. 
\end{remark}

\bigskip

For later use, let for all $s\le t\in \r$, and $v\in [0,1]$,
$$
Y_{s,t}^*(v) := \inf\{ u \in[0,1]\, : \, Y_{-t,-s}(u) > v\}\land 1.
$$

\begin{lemma}\label{l:Jdual}
For all $s,t\in \r$, and $v\in [0,1]$, $Y_{s,t}(v) = \inf\{ u \in[0,1]\, : \, Y_{t,s}(u) > v\}\land 1$. In particular, for $s\le t$ and $v\in [0,1]$, $Y^*_{s,t}(v)=Y_{-s,-t}(v)$. 
\end{lemma}

\noindent {\it Proof}.  Let $s,t\in \r$. The case $v=1$ is immediate so we can suppose $v<1$. By definition of $Y$, we have
$$
Y_{s,t}(v) = {S_{\eta_f^{-1}(s),\eta_f^{-1}(t)} (v S_{0,\eta_f^{-1}(s)}(b)) \over S_{0,\eta_f^{-1}(t)}(b)}.
$$
 
\noindent By the duality \eqref{eq:Sdual} for $S$, 
$$
S_{\eta_f^{-1}(s),\eta_f^{-1}(t)} (c)
=
\inf\{ a\ge 0\,:\,S_{\eta_f^{-1}(t),\eta_f^{-1}(s)} (a) >c \}.
$$

\noindent Taking $c=v S_{0,\eta_f^{-1}(s)}(b)$, we get
\begin{eqnarray*}
Y_{s,t}(v) &=& {1\over  S_{0,\eta_f^{-1}(t)}(b)} \inf\left\{ a\ge 0\,:\,{S_{\eta_f^{-1}(t),\eta_f^{-1}(s)} (a) \over S_{0,\eta_f^{-1}(s)}(b)} > v \right\}\\
&=&
 \inf\left\{ u\ge 0\,:\,{S_{\eta_f^{-1}(t),\eta_f^{-1}(s)} (u S_{0,\eta_f^{-1}(t)}(b)) \over S_{0,\eta_f^{-1}(s)}(b)} > v \right\} 
\end{eqnarray*}
 
\noindent which is $\inf\{ u \in[0,1]\, : \, Y_{t,s}(u) > v\}$ indeed. $\Box$

\section{Disintegration of the PRBM with respect to its occupation field}

\subsection{The burglars}
\label{s:burglar}

Let $b>0$. Consider the PRBM $X$ defined in \eqref{def:muprocess} up to $\tau_b^0$. Recall from \eqref{def:S} that 
$ S_{r,x}(a) = L(\tau_a^r, x )$. Set
$$
t^*_0:=\sup\left\{t\le 0\,:\, L(\tau_b^0,X_t)=0\right\}, \qquad 
t^*:=\inf\left\{t\ge 0\,:\, L(\tau_b^0,X_t)=0 \right\}.
$$ 

\noindent Notice that $t^*$ is almost surely the hitting time of $\inf_{[0,\tau_b^0]} X$ and by the Ray--Knight theorems, $|t^*_0|< \infty$ if and only if $\mu >1$ (if $\mu\le 1$, $\{t\le 0\,:\, L(\tau_1^0,X_t)=0\}=\emptyset$ hence $t_0^*=-\infty$).  % The time $t^*$ is the time when $(X_t,\, t\in [0,\tau_1^0])$ attains its minimum. 

 We let $f(x):=S_{0,x}(b)=L(\tau_b^0, x)$ restricted to the interval $I$ with boundaries 
\begin{equation}\label{eq:defI}
\inf I:= \inf_{[0, \tau^0_b]}X =X_{t^*} \qquad \sup I:=\inf\{x\ge0: L(\tau_b^0, x)=0\}
\end{equation}

\noindent ($\sup I =X_{t_0^*}$ when $|t_0^*|<\infty$ and $\infty$ otherwise).  We consider the processes
$$
X^{(1)}_s:= X_s,\, s\in (t_0^*,t^*), \qquad X^{(2)}_s:= X_{\tau_b^0-s},\, s\in [0,\tau_b^0-t^*).
$$

\begin{figure}[h!]
\begin{minipage}{\textwidth}
\center
\begin{tikzpicture}[line cap=round,line join=round,>=triangle 45,scale=0.7, x=0.5cm,y=0.5cm]
\clip(-6,-3.2) rectangle (16.48,6.3);
\draw[->] (3.,-4.) -- (3.,5.); % axe y
\draw[->] (-2.5,0)--(14.5,0.); % axe x
\draw[dashed] (8.64,-3.06) -- (8.64,0.); % minimum
\node at (8.64,0.8) {\scriptsize $t^*$};
\node at (14.,-0.9) {\scriptsize $\tau_1^0$};
\node at (2.5, -0.5) {\scriptsize $0$};
\draw (3.,0.)-- (2.58,0.76)-- (2.22,0.24)-- (1.82,1.26)-- (1.4,1.44)-- (1.,1.)-- (0.6,2.56)-- (0.22,1.74)-- (-0.22,2.34)-- (-0.62,1.74)-- (-1.,3.)-- (-1.4,2.44)-- (-1.78,3.)-- (-2.18,4.28)-- (-2.58,3.52);
\draw  (3.,0.)-- (3.4,0.36)-- (3.78,-0.86)-- (4.2,-0.2)-- (4.6,-0.62)-- (5.,0.38)-- (5.42,-0.24)-- (5.78,0.96)-- (6.18,0.56)-- (6.6,-0.46)-- (7.,-0.82)-- (7.44,-2.1)-- (7.8,-1.62)-- (8.24,-2.06)-- (8.64,-3.06)-- (8.98,-2.4)-- (9.36,-2.64)-- (9.74,-2.)-- (9.98,-1.18)-- (10.4,-1.7)-- (10.8,-1.4)-- (11.18,-0.62)-- (11.58,-1.06)-- (11.98,-0.32)-- (12.38,0.6)-- (12.66,-0.38)-- (13.,0.24)-- (13.38,-0.4)-- (13.64,0.);
\end{tikzpicture}
\end{minipage}

\begin{minipage}{0.4\textwidth}
\begin{tikzpicture}[line cap=round,line join=round,>=triangle 45, scale=0.7, x=0.5cm,y=0.5cm]
\draw [ xstep=1.0cm,ystep=1.0cm] ;
\clip(-6,-5.) rectangle (10.7,6.3);
\draw[->](-3.,0.) -- (10.5,0.); %axe x
\draw[->] (3.,-5) -- (3.,4.5); % axe y
\node at (2.5, -0.5) {\scriptsize $0$};
\draw  (3.,0.)-- (2.62,0.5)-- (2.2,0.14)-- (1.8,1.18)-- (1.42,0.74)-- (1.,1.)-- (0.6,2.28)-- (0.24,1.1)-- (-0.18,2.)-- (-0.62,1.82)-- (-0.78,2.86)-- (-1.2,3.72)-- (-1.62,2.9)-- (-2.,2.72)-- (-2.38,3.08)-- (-2.8,3.96);
\draw  (3.,0.)-- (3.38,0.54)-- (3.6,-1.12)-- (4.,-0.32)-- (4.38,-0.62)-- (4.78,-1.12)-- (5.18,-0.46)-- (5.54,0.26)-- (5.84,-0.48)-- (6.22,0.76)-- (6.58,0.28)-- (7.,-0.78)-- (7.44,-3.02)-- (7.76,-2.58)-- (8.24,-3.2)-- (8.6,-2.24)-- (8.98,-4.16)-- (9.38,-3.58)-- (9.84,-4.38);
\draw[dashed] (9.84,-4.38) -- (10.7,-4.7);
\end{tikzpicture}
\end{minipage}
\begin{minipage}{0.4\textwidth}
\begin{tikzpicture}[line cap=round,line join=round,>=triangle 45,scale=0.7, x=0.5cm,y=0.5cm]
\draw [ xstep=1.0cm,ystep=1.0cm] ;
\clip(-9,-5.) rectangle (4.,6.3);
\draw[<-](-5,0.) -- (4.5,0.); %axe x
\draw[->] (3.,-5) -- (3.,4.5); % axe y
\node at (3.5, -0.5) {\scriptsize $0$};
\draw  (-4.5,-4)--(-4.,-3.2)--(-3.75,-3.5)--  (-3.5, -2.96)--(-3.25,-3.5) -- (-3.,-3.12) -- (-2.75,-2.4)--(-2.5,-2.90)-- (-2.25,-2.75)-- (-2.,-2.25)--(-1.75,-2)  --(-1.5,-2.2)--(-1.25,-1.3)--(-1.,-1.44) -- (-0.75,-1)--(-0.5,-1.20)--(-0.25,-0.8) -- (0.,-0.50)-- (0.25,-1)-- (0.5, -1.61)--(0.75,0.8)-- (1.,0.08)--(1.25,-0.5)-- (1.5,0.2) --(1.75,-0.7)--(2.,0.7)-- (2.25,0.5)--(2.5,0.10 ) --(2.75, -0.5)--(3.,0.00);
\draw[dashed] (-4.6,-4.1)--(-5.,-4.5);
\end{tikzpicture}
\end{minipage}
\caption{%\leftskip=1.8truecm \rightskip=1.8truecm 
\small Illustration of the two burglars $Z^{(1)}$ and $Z^{(2)}$. Top: a PRBM with $\mu\in (0,1]$ up to $\tau_1^0$, with $t_0^*=-\infty$  and $t^*$ being the time associated with the minimum. Bottom: the process before time $t^*$ is mapped through a space-time transformation to the process $Z^{(1)}$ (left). That after time $t^*$ is mapped to the process $Z^{(2)}$ (right).}
\label{f:burglar0}
\end{figure}
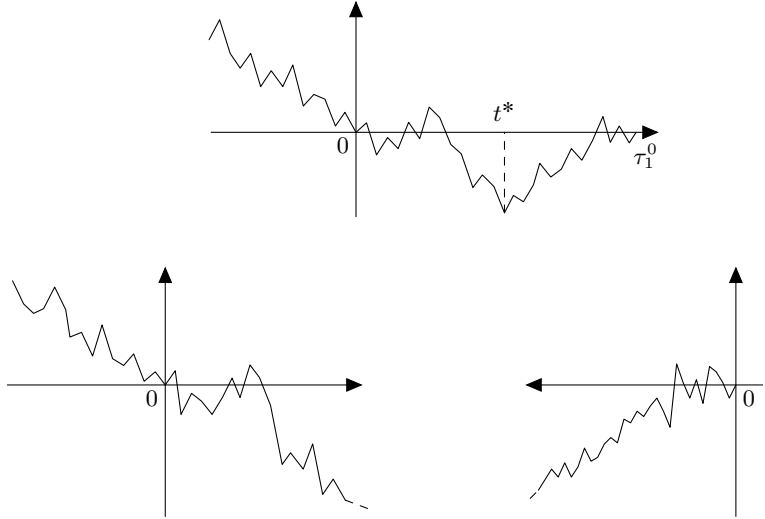

\noindent Recall the notation $\Upsilon$ in Section \ref{s:process}. We introduce the processes 
\begin{equation}\label{def:Z12}
Z^{(1)}:=\Upsilon(X^{(1)},f), \qquad  Z^{(2)}:= \Upsilon(X^{(2)},f),
\end{equation}

\noindent see Figure \ref{f:burglar0}. In the case when $\mu=1$, $Z^{(2)}$ defines a variant of the Brownian burglar introduced by Warren and Yor in \cite{warren-yor} (its time-change $\widehat Z^{(2)}$ in Section \ref{s:bassburdzy} is distributed as the one appearing in \cite{warren}, compare Theorem \ref{p:bassburdzy} and Proposition \ref{p:parameters} with Theorem 1 in \cite{warren}).%, using for example the forthcoming Lemma \ref{l:JA}).  

\bigskip

\begin{lemma}\label{l:eta}
Almost surely, the process $Z^{(1)}$ is defined on $\r$ with limits given by $\lim_{t\to -\infty} Z^{(1)}_t=+\infty$ and $\lim_{t\to \infty} Z_t^{(1)}=-\infty$; the process $Z^{(2)}$ is defined on $\r_+$ and $\lim_{t\to\infty} Z^{(2)}_t=-\infty$.
\end{lemma}

\noindent {\it Proof}.  We only prove it for $Z^{(1)}$. First we check that the transformation $\eta_f$ maps the interval $I$ defined in \eqref{eq:defI} onto $\r$. Consider $\sup \eta_f(I)$. When $\mu\le 1$, $\sup I=\infty$ and $\lim_{y\to \infty} \eta_f(y)=\int_0^{\infty} {\d r \over f(r)}$ which is infinite (Exercise X.3.20 in \cite{revuz-yor}). In the case $\mu>1$, we need to show that $\int_0^{\sup I} \frac{\d r}{f(r)}=\infty$. By time reversal, $r\to f(\sup I-r)$ is a ${\rm BESQ}^{4-\delta}_0$ process, which, by scaling arguments, implies that the integral is infinite indeed.  Consider $\inf \eta_f(I)$. It is $\int_0^{\inf I}\frac{\d r}{f(r)}$ which is $-\infty$ by the same argument. Finally, we need to show that $\int_{0}^{t_0^*} \frac{\d s}{f(X_s^{(1)})^2}=-\infty$ and $\int_{0}^{t^*} \frac{\d s}{f(X_s^{(1)})^2}=\infty$. By the occupation time formula, $\int_0^{t_0^*} \frac{\d s}{f(X_s^{(1)})^2} = \int_{\inf I}^{\sup I} \frac{L(t_0^*,r)-L(0,r)}{f(r)^2}\d r = -\int_{0}^{\sup I} \frac{L(0,r)}{f(r)^2}{\d r}$. When $r<\sup I$ is close enough to $\sup I$, $L(0,r)=f(r)$, so that the integral is $-\infty$ indeed.
Similarly, 
\begin{align*}
\int_{0}^{t^*} \frac{\d s}{f(X_s^{(1)})^2}&= \int_{\inf I}^{\sup I} \frac{L(t^*,r)-L(0,r)}{f(r)^2}\d r \\
&\ge \int_{\inf I}^{0} \frac{L(t^*,r)-L(0,r)}{f(r)^2}\d r \\
&=
 \int_{\inf I}^{0} \frac{L(t^*,r)}{f(r)^2}\d r.
\end{align*}

\noindent On an event of probability arbitrarily close to $1$ as $\varepsilon\to 0$, $L(t^*,\inf I+y)$ and $f(\inf I +y)-L(t^*,\inf I+y)$ on $[0, \varepsilon]$ are independent BESQ processes of respective dimensions $\delta$ and $2$, starting at $0$ (see Theorem 4.3 in \cite{eyzmupro}). Scaling arguments ensure again that the integral is infinite indeed. $\Box$

\bigskip

Using the notation $\Psi$ of Section \ref{s:flow}, define the process $\Xi$ by 
\begin{equation}\label{def:Xi}
\Xi := \Psi(L(t^*,\cdot),f).
\end{equation}

\noindent Recall from Remark \ref{r:Y} that $Y$ denotes the flow $\Psi({\mathcal L}_X^b,f)$. Its forward flow is a Jacobi($2/\mu,0$) flow and its backward flow is a killed Jacobi($2-2/\mu$) flow. The following lemma shows that the process $\Xi$ is measurable with respect to the flow $Y$.  
\begin{lemma}
Almost surely, for all $r\in \r$,   
\begin{equation}\label{def-Xi}
\Xi_r = \inf\{ v\in[0,1]\,:\, Y_{r,x}(v)=1  \hbox{ for some } x < r \}.
\end{equation}
\end{lemma}

\noindent {\it Proof}. By definition, $\Xi_r = \frac{L(t^*, \eta^{-1}_f(r))}{f(\eta^{-1}_f(r))}, $ with $\eta_f(s):= \int_0^s {\d r \over f(r)}$ for $s\in I$.  By Definition \ref{def:imageflow}, $Y_{r, x}(v)= \frac1{f(\eta^{-1}_f(x))} L(t_v, \eta_f^{-1}(x))  $, with $t_v:= \tau_{v f(\eta^{-1}_f(r))}^{\eta_f^{-1}(r)}$. We have $t_v > t^*$ for any $1\ge v>  \Xi_r$ and  $\eta_f^{-1}(x) \to X_{t^*}$ as $x \to -\infty$. Observe that $Y_{r, x}(v)=1$ for all $x$ such that $\eta_f^{-1}(x) \le \inf_{[t_v,\tau_b^0]} X$. On the other hand, if $v < \Xi_r$, then $t_v < t^*$ and $Y_{r, x}(v)<1$ for all $x < r$. $\Box$

\bigskip
The law of the process $(\Xi_r, r \in \r)$, called the ``primitive Eve'' process in Bertoin and Le Gall \cite{bertoin-legall-I}, will be given in Proposition \ref{p:eve}. 

Let $L^{(1)}(t,r)$, resp. $L^{(2)}(t, r)$, denote the local time of $Z^{(1)}$, resp. $Z^{(2)}$, at time $t$ and position $r$. Observe from \eqref{eq:locZ} that $L^{(1)}(\infty,r)+L^{(2)}(\infty,r)=1$ for all $r\in \r$. On the other hand,  \eqref{def:Xi} implies that  $\Xi_r= L^{(1)}(\infty, r)$.

\begin{theorem}\label{p:loctime}

(i) The process $Z^{(1)}$ possesses local time flow $Y$ ``seen from the left'', meaning that $L^{(1)}(\tau^r_v(Z^{(1)}), x)=Y_{r,x}(v)$ for $v \in [0, \Xi_r]$ and $r,x\in\r$. 

(ii)  The process $Z^{(2)}$ possesses local time flow $Y$ ``seen from the right'', meaning that $L^{(2)}(\tau^r_{u-}(Z^{(2)}), x)=1-Y_{r,x}(1-u)$ for $u \in [0,1-\Xi_r]$ and $r,x\in\r$. 

(iii) The processes $Z^{(1)}$ and $Z^{(2)}$ are independent of $(L(\tau_b^0,x), \, x\in \r)$ (and their distributions do not depend on the value of $b>0$).
\end{theorem}

{\noindent\it Proof}.  (i)   Note that for any $ r \in I$, $L_{X^{(1)}}(t, r)= L(t, r)$ for all $t \le t^*$, hence the local time flow of $X^{(1)}$ is the flow ${\mathcal S}^1:=\{S_{r,x}(a),\, r,x \in I, \, 0\le a\le L(t^*,r)\}$. The image flow by $\Psi(\cdot,f)$ is the flow $Y$ ``seen from the left'' by definition.  

(ii)  Similarly  for any $ r \in I$, $L_{X^{(2)}}(t, r)= f(r) - L(\tau_b^0-t, r)$ for all $0\le t \le \tau_b^0-t^*$. We deduce that  the local time flow of $X^{(2)}$ is the flow $\{f(x)-S_{r,x}((f(r)-a)-),\, r,x \in I, \, 0\le a\le f(r)-L(t^*,r)\}$. The image flow by $\Psi(\cdot,f)$ is the flow $Y$ ``seen from the right''.  
 
 (iii)   Again $Z^{(1)}$ and $Z^{(2)}$ are measurable with respect to their local time flow $Y$ (by a  proof similar to Proposition \ref{p:Xflow}). By Lemma \ref{l:Jdual}, the flow $Y$ is measurable with respect to its forward flow. The forward flow of $Y$ is independent of $(L(\tau_b^0,x), \, x\in \r)$ by Theorem \ref{t:perkins}.  It completes the proof. $\Box$

\subsection{A Markov property for the process $Z^{(2)}$}

The following proposition gives a way of constructing $Z^{(2)}$ by means of a Brownian motion stopped at a hitting time and an independent squared Bessel process.

\begin{proposition}\label{p:Z2A}
Let $g<0$ be a constant. Take a Brownian motion up to the hitting time of $g$ (call it $B=(B_t, t\in [0,T_{g}^B])$). Take $h$ such that $h(x)=0$ if $x\le g$ and $h(g+r)$, $r\ge 0$ is  an independent ${\rm BESQ}^\delta_0$. Let $f(x):= L_B(T_{g}^B,x)+h(x)$, $x\in \r$ where $L_B$ is the local time of the Brownian motion $B$. Then $\Upsilon(B,f)$ has the law of $Z^{(2)}$ and is independent of $f$. 
\end{proposition}

\noindent {\it Proof}. We construct a probability measure $\widetilde \p$ under which $X^{(1)}$ and $X^{(2)}$ have simple descriptions. The process $B$ will stand for the process $X^{(2)}$ while $h$ will stand for the local time process of $X^{(1)}$ under $\widetilde \p$ conditioned on some event (which is $\{J(A)=g\}$ in the notation below).

Let $m>|g|$ be an arbitrary constant (the exact value of $m$ plays no role). We consider a probability $\widetilde \p$ on the product space $\r_+\times C(\r,\r)$ such that, under $\widetilde \p$, $(A,X)$ has the following distribution: $A$ is gamma(${\delta\over 2},2m$) distributed (meaning it has density ${(2m)^{-{\delta \over 2}}\over \Gamma\left({\delta\over 2}\right)} x^{{\delta\over 2}-1} e^{- x/2m}$ on $\r_+$) and conditionally on $A=b$, $X$ has distribution 
$$
c_b |J(b)|^{{\delta\over 2}-1} {\bf 1}_{\{J(b)>-m\}}\cdot \p
$$

\noindent where  $J(b):=\inf_{t\in [0,\tau_b^0]} X_t$ and $c_b=\Gamma(\frac\delta2) \left({b\over 2}\right)^{1-{\delta\over 2}}e^{{b\over 2m}}$ is the renormalizing constant. Notice that $J(b)$ is measurable with respect to $L(\tau_b^0,x)$, $x\in \r$.  Let $X^{(1)}$, $X^{(2)}$, $Z^{(1)}$, $Z^{(2)}$ be under $\widetilde \p$ the processes of Section \ref{s:burglar} with $b=A$. From Theorem \ref{p:loctime} (iii), and using that $J(A)$ is measurable with respect to $(L(\tau_A^0,x), x \in \r)$, we deduce that under $\widetilde \p$, $Z^{(2)}$ has the same law as under $\p$,  and  is independent of  $(L(\tau_A^0,x), x \in \r)$. 

  Denote by
\begin{equation}
    \label{def-T} 
    T_r :=\inf\{ t\in \r : \, X_t = r \},
\end{equation} 

\noindent the hitting time of $r \in \r$. %Note that $T_r \to -\infty$ as $r \to \infty$. 

\begin{lemma}\label{l:JA}
%(i) Under $\widetilde \p$, ${1\over m} J(A)$ is  beta(${\delta\over 2},1$) distributed.

(i) Under $\widetilde \p(\cdot | J(A)=g)$, the law of the process $X^{(1)}$ is the one of $(X_s,\,s\le T_g)$ under $\p$.

(ii) Under $\widetilde \p(\cdot | J(A)=g)$, the law of the process $X^{(2)}$ is that of a Brownian motion stopped when hitting $g$.

(iii) Under $\widetilde \p(\cdot | J(A)=g)$, the two processes $X^{(1)}$ and $X^{(2)}$ are independent.

\end{lemma}
{\noindent\it Proof of Lemma \ref{l:JA}}. These results can be deduced from \cite{eyzmupro}. First, the law of $(X_s, \,s\le 0)$ is the same under $\p$ and $\widetilde \p$, since $(X_s, \, s \le 0)$ and $(X_s, \, s\ge 0)$ are independent. We restrict our attention to $(X_s, \, s\in [0,\tau_A^0])$. By Proposition 9.1 p.123, Section 9.2 of \cite{yor92} and Corollary 3.4 in \cite{eyzmupro}, we observe that $(A,X_s, \, s\in [0,\tau_A^0])$ under $\widetilde \p$ has the distribution of $(L(T_{-m},0),X_s, \, s\in [0,{\tt g}_m])$ under $\p$,  where ${\tt g}_m:= \sup\{t\le T_{-m}: X_t=0\}$. %The law of $(X_s, \, s\in [0,\tau_b^0])$ under $\widetilde \p(\cdot | A=b)$ is the one of $(X_s, \, s\in [0, {\tt g}_m])$ under $\p(\cdot | L(T_{-m},0)=b)$,. We deduce that the law of $(X_s, \, s\in [0,\tau_A^0])$ under $\widetilde \p$ is the one of  $(X_s, \, s\in [0,{\tt g}_m])$ under $\p$. 
The lemma is then Theorem 3.2 of \cite{eyzmupro}. $\Box$

\bigskip

We can now complete the proof of Proposition \ref{p:Z2A}. Let $g\in (-m,0)$. We recall that $Z^{(2)}$ is independent of  $(L(\tau_A^0,x), x \in \r)$ under $\widetilde \p$ so that the law of $Z^{(2)}$ is identical under $\widetilde  \p(\cdot | J(A) = g)$ for any $g\in (-m,0)$ (and is equal to the law of $Z^{(2)}$ under $\p$). By definition, $Z^{(2)}=\Upsilon(X^{(2)},f)$ where $f(x):=L(\tau_A^0,x)$. The statement of the proposition comes from the description of $X$ under $\widetilde  \p(\cdot | J(A) = g)$ in Lemma \ref{l:JA}. This holds for any $g<0$ since $m$ is arbitrary. $\Box$

\bigskip

 We state now the Markov property for $Z^{(2)}$. It will be the key ingredient in the forthcoming connection with the Bass--Burdzy flow. Let as before $L^{(2)}(t,x)$ denote the local time of $Z^{(2)}$ at time $t$ and position $x$. Recall Notation \ref{n:upsilon}.

\begin{theorem}\label{t:markov}
 For any $t\ge 0$, $\Upsilon(Z^{(2)}_{t+\cdot}, 1-L^{(2)}(t,\cdot))$ is independent of $(Z^{(2)}_r, \, r\in [0,s])$ and distributed as $Z^{(2)}$. In particular, $(Z^{(2)}_{t}, (L^{(2)}(t,x))_{x\in \r})$ is a Markov process.  
 \end{theorem}

\bigskip

Before proving the theorem, we show that for fixed $t\ge 0$, the process $\Upsilon(Z^{(2)}_{t+\cdot}, 1-L^{(2)}(t,\cdot))$ is indeed well-defined.

\begin{lemma}\label{l:continuity2}
 For any fixed $t\ge 0$, the map $s \mapsto \Upsilon(Z^{(2)}_{s+\cdot}, 1-L^{(2)}(s,\cdot))$ is a.s. well-defined on a neighborhood (in $\r_+$) of $t$ and is continuous at $t$ (in the space $C(\r_+,\r)$ of continuous processes on $\r_+$ endowed with the usual topology of uniform convergence over all compacts).
\end{lemma}

See Appendix \ref{s:continuitylemmas} for the proof of Lemma \ref{l:continuity2}. 

\bigskip

\noindent {\it Proof of the theorem}. Take $B$ and $f$ as in Proposition \ref{p:Z2A}, with, say, $g=-1$. Recall the notation $C_f(t):=\int_0^t {{\rm d} s \over f({B_s})^2}$. Set $Z^{(2)}:=\Upsilon(B,f)$. Applying Lemma \ref{l:comp} to ${\mathcal X}=B$, $s\in J=[0, T^B_{-1}), I= (-1, \max_{0\le t \le T^B_{-1}}B_t]$ and $c'=B_s$, we have
$$
\Upsilon(B_{s+\cdot},f - L_B(s,\cdot)) = \Upsilon(Z^{(2)}_{C_f(s)+\cdot},1-L^{(2)}(C_f(s),\cdot)).
$$

\noindent On the other hand, by Proposition \ref{p:Z2A}, the Markov property of the Brownian motion at time $s$ and Notation \ref{n:upsilon} (i), {conditionally on $\{s < T^B_{-1}\}$,  the left-hand side is distributed as $Z^{(2)}$ and is independent of $f-L_B(s,\cdot)$ and of $(B_r,\,r\in [0,s])$. 

 It follows that, conditionally on $\{s < T^B_{-1}\}$,  $\Upsilon(Z^{(2)}_{C_f(s)+\cdot},1-L^{(2)}(C_f(s),\cdot))$ is independent of $\sigma\{f, B_r,\,r\in [0,s]\}$ and has the law of $Z^{(2)}$. Using the continuity of $s\to \Upsilon(Z^{(2)}_{s +\cdot},1-L^{(2)}(s,\cdot))$ stated in Lemma \ref{l:continuity}, and 
since $C_f(s)$ and $(Z^{(2)}_t, \,t\le C_f(s))$ are $\sigma\{f, B_r,\,r\in [0,s]\}$-measurable, we deduce  that  for all $t\ge 0$, $\Upsilon(Z^{(2)}_{t+\cdot}, 1-L^{(2)}(t,\cdot))$ is independent of $(Z^{(2)}_r, \, r\in [0,t])$ and is distributed as $Z^{(2)}$.  $\Box$

 \begin{remark}\label{R:markovfort} As a consequence of Theorem \ref{t:markov}, if $T$ is a stopping time with respect to the natural filtration of $Z^{(2)}$ such that almost surely $ s \to \Upsilon(Z^{(2)}_{s+\cdot}, 1-L^{(2)}(s,\cdot))$ is continuous at $T$ in $C(\r_+,\r)$, then $\Upsilon(Z^{(2)}_{T+\cdot}, 1-L^{(2)}(T,\cdot))$ is independent of $(Z^{(2)}_r, \, r\in [0,T])$ and distributed as $Z^{(2)}$.
 \end{remark}

\bigskip

The next result gives an invariance principle for $Z^{(2)}$ at small times.   

\begin{theorem}\label{t:scaling}
The process $\left( {1\over \sqrt{a}} Z^{(2)}_{at},\, t\ge 0\right)$ converges in distribution as $a \downarrow 0$ to a standard Brownian motion.
\end{theorem}
  
\noindent {\it Proof}.  We use again the representation  $Z^{(2)}=\Upsilon(B,f)$ from Proposition \ref{p:Z2A} with, say, $g=-1$. By definition of $\Upsilon$ in Section \ref{s:process}, $Z_t^{(2)} = \eta_f(B_{C_f^{-1}(t)})$ with $\eta_f(x)=\int_0^x {\d r \over f(r)}$, $C_f(t)= \int_0^t {\d s\over f(B_s)^2}$. Let $a>0$. Define $\widetilde f(r):=f(r\sqrt{a})$, $\widetilde B_t:= {1\over \sqrt{a}} B_{at}$, $\widetilde Z_t^{(2)}:= {1\over \sqrt{a}} Z_{at}^{(2)}$. We check that $\eta_{\widetilde f}(\widetilde B_t) = {1\over \sqrt{a}} \eta_f(B_{at})$ which is ${1\over \sqrt{a}} Z^{(2)}_{C_f(at)}$. Moreover $C_{\widetilde f}(t) = \int_0^t \frac{\d s}{\widetilde f(\widetilde B_s)^2}={1\over a} C_f(at)$ so that $\eta_{\widetilde f}(\widetilde B_t) = {1\over \sqrt{a}} Z^{(2)}_{a C_{\widetilde f}(t)}$. We proved that $\Upsilon(\widetilde B, \widetilde f) = \widetilde Z^{(2)}$. By scaling $\widetilde B$ is a standard Brownian motion stopped when hitting $-{1\over \sqrt{a}}$. Moreover, $\widetilde f$ converges when $a\downarrow 0$ to the constant function $f(0)$,  and $\widetilde B$ and $f(0)$ are asymptotically independent.  Use Lemma \ref{l:continuity} (together with Skorokhod's representation theorem to suppose that the convergence of $(\widetilde B, \widetilde f)$ is almost sure instead of in distribution) to complete the proof. $\Box$

\bigskip
We end this subsection by describing the law of the process $(\Xi_r)_{r\in \r}$ defined in \eqref{def-Xi}.  Let $a, b>0$. Denote by beta($a,b$) the distribution with density  ${\Gamma(a+b)\over \Gamma(a)\Gamma(b)}x^{a-1}(1-x)^{b-1}{\bf 1}_{(0,1)}(x)$.

\begin{proposition}\label{p:eve} The random variable $\Xi_0$ is distributed as \mbox{\rm beta}$(\delta/2, 1)$.  Conditionally on $\Xi_0$,  $(\Xi_{-t})_{t\ge 0}$ and $(\Xi_t)_{t\ge0}$ are independent  {\rm Jacobi($\delta,2$)} and  {\rm Jacobi($\delta,0$)} processes. 
\end{proposition}  

\noindent {\it Proof}. Using the representation $Z^{(2)}=\Upsilon(B,f)$ in Proposition \ref{p:Z2A} with $g=-1$, we deduce from  \eqref{eq:locZ} that $L^{(2)}(t,x)= \frac{L^B(C_f^{-1}(t), \eta_f^{-1}(x))}{f(\eta_f^{-1}(x))}$, with $C_f(t):=\int_0^t \frac{\d s}{f(B_s)^2}$, $\eta_f(x):=\int_0^x \frac{\d r}{f(r)}$. Note that $f$ is a ${\rm BESQ}^{2+\delta}_0$ process on $[-1, 0]$, and a ${\rm BESQ}^\delta$ on $\r_+$. Since ${\mathfrak d}_f:=\inf\{x>0: f(x)=0\}= \inf\{x>\sup_{0\le s \le T^B_{-1}} B_s: h(x)=0\}$, where we recall that $h$ is a ${\rm BESQ}^\delta$ on $[-1, \infty)$ with $h(-1)=0$, independent of $B$,    $C_f^{-1}(t)\to T^B_{-1}$ as $t\to\infty$ and  $L^{(2)}(\infty,x)= \frac{L^B(T^B_{-1}, \eta_f^{-1}(x))}{f(\eta_f^{-1}(x))}$ for $x\in \r$. By the classical Ray--Knight theorem, $(L^B(T^B_{-1}, y), y\ge -1)$ is a ${\rm BESQ}_0^2$ on $[-1, 0]$, and ${\rm BESQ}^0$ on $[0, \infty)$.  Since $f(y)= L^B(T^B_{-1}, y) + h(y)$, we deduce that  $L^{(2)}(\infty,0)= \frac{L^B(T^B_{-1}, 0)}{f(0)}$ is distributed as beta$(1, \delta/2)$. By  Warren and Yor \cite{warren-yor}, $L^{(2)}(\infty, \cdot)$ is a Jacobi process of parameters $(0, \delta)$ on  $\r_+$ starting from $L^{(2)}(\infty, 0)$.   Let $\varepsilon \in (0, 1)$. Again by Warren and Yor \cite{warren-yor}, $(L^{(2)}(\infty, \eta_f(-1+\varepsilon)+x),   0\le x\le -\eta_f(-1+\varepsilon))$ is a
Jacobi process of parameters $(2, \delta)$. Recall that $\eta_f$ is independent of $L^{(2)}(\infty, \cdot)$.  Observe that $L^{(2)}(\infty, \eta_f(-1+\varepsilon))$ is still distributed as beta$(1, \delta/2)$ which is in fact the stationary distribution of a Jacobi process of parameters $(2, \delta)$. By time-reversal, $(L^{(2)}(\infty, -x),   0\le x \le -\eta_f(-1+\varepsilon))$ is a Jacobi process of parameters $(2, \delta)$.    Let $\varepsilon\to 0$ we get the proposition by using the fact that $\Xi_r= 1- L^{(2)}(\infty, r)$ for all $r \in \r$.   $\Box$

\subsection{Link with the Bass--Burdzy flow and proof of Theorem III}

\label{s:bassburdzy}

 Let $\beta_1,\, \beta_2\in \r$, $\sigma>0$ and $\gamma$ be a standard Brownian motion. Following \cite{bassburdzy}, one defines the Bass--Burdzy flow of parameters $(\beta_1,\beta_2)$ and diffusivity $\sigma$ as 
 the collection of homeomorphisms of the real line $({\mathcal R}_t,\,t\ge 0)$ such that  for any $x\in \r$, the process $({\mathcal R}_t(x),\, t\ge 0)$ is the strong solution of the SDE
$$
{\mathcal R}_t(x) := x+ \sigma \gamma_t + \beta_1\int_0^t {\bf 1}_{\{{\mathcal R}_s(x)< 0\}} \d s +\beta_2\int_0^t {\bf 1}_{\{{\mathcal R}_s(x)> 0\}} \d s.
$$

\noindent  When $\sigma=1$, we just call it Bass--Burdzy flow with parameters $(\beta_1,\beta_2)$ as defined in the introduction.

We continue to look at the process $Z^{(2)}$ and show that it is linked to the Bass--Burdzy flow via a time-change, extending the result in \cite{warren} to all parameters of the PRBM. Using the same time-change as in  \cite{warren}, we set
\begin{equation}\label{eq:zeta}
\zeta(t):=\int_0^t {{\rm d} s \over (1-L^{(2)}(s,Z^{(2)}_s))^2}, \qquad \widehat Z^{(2)}_t :=  Z^{(2)}_{\zeta^{-1}(t)}.
\end{equation}

\begin{lemma} \label{l:tmax}
We have $\lim_{t\to t_{\max}} \zeta(t)=\infty$ where $t_{\max}:=\inf\{s\ge 0\,:\,L^{(2)}(s,Z^{(2)}_s)=1\}\in(0,\infty]$ (it is finite almost surely if $\mu>1$ and infinite otherwise). Consequently, if $t_{\max}<\infty$, $\lim_{t\to \infty} \widehat Z^{(2)}_t = Z^{(2)}_{t_{\max}}$, and if $t_{\max}=\infty$, $\lim_{t\to \infty} \widehat Z^{(2)}_t = -\infty$.
\end{lemma}

 \noindent {\it Proof}. It suffices to  treat the case $\mu>1$. Recall the construction of $Z^{(2)}$ in Proposition \ref{p:Z2A} and recall from  \eqref{eq:locZ} that $L^{(2)}(t,x)= \frac{L_{B}(C_f^{-1}(t), \eta_f^{-1}(x))}{f(\eta_f^{-1}(x))}$, with $C_f(t):=\int_0^t \frac{\d s}{f(B_s)^2}$, $\eta_f(x):=\int_0^x \frac{\d r}{f(r)}$ and $f(x)=L_B(T_g^B,x)+h(x)$ where $h(g+\cdot)$ is a ${\rm BESQ}^\delta_0$ independent of $B$. Let $\widetilde t_{\max}:=C_f^{-1}(t_{\max})$ which is $\inf\{s\in [0,T_g^B]\,:\, L_B(s,B_s)=f(B_s)\}$. Using the occupation time formula, we can express $\zeta(t_{\max})$ as
 $$
 \zeta(t_{\max})= \int_\r \frac{L_B(\widetilde t_{\max},y)}{f(y)(f(y)-L_B(\widetilde t_{\max},y))} \d y.
 $$
 
 \noindent Let $y_0:= B_{\widetilde t_{\max}}$. We observe that  $y_0$ is the last zero of $h$ before $M:=\max_{[0,T_g^B]} B$ and $\widetilde t_{\max}$ is  the last passage time of $B$  at $y_0$ before $T_g^B$. Hence, for any $y\ge y_0$, $L_B(\widetilde t_{\max},y)=L_B(T_g^B,y)$ so that 
 $$
 \zeta(t_{\max})\ge \int_{y_0}^{\infty} \frac{L_B(T_g^B,y)}{f(y)(f(y)-L_B(T_g^B,y))} \d y
 =
 \int_{y_0}^{\infty} \frac{L_B(T_g^B,y)}{f(y)h(y)} \d y.
 $$
 
 \noindent Notice that in a neighborhood of $y_0$,  the numerator is positive and $f$ is bounded from above while $\int_{y_0}^{y_0+\varepsilon} \frac{\d y}{h(y)}=\infty$ for all $\varepsilon$ positive.\footnote{
 For example, we can express $h$ in terms of the square of a Bessel meander whose law is absolutely continuous with respect to BESQ, see Section 3.6 in \cite{yor92}.} $\Box$

 \bigskip

 The local time of $\widehat Z^{(2)}$ can be computed since
\begin{align*}
\int_0^t {\bf 1}_{\left\{\widehat Z^{(2)}_{s} < x\right\}} \d s
&=
\int_0^{\zeta^{-1}(t)} {\bf 1}_{\left\{ Z^{(2)}_{r} < x\right\}} \frac{\d r}{(1-L^{(2)}(r,Z^{(2)}_r))^2}\\
&=
\int_{-\infty}^x \d y \int_0^{\zeta^{-1}(t)}  \frac{\d_r L^{(2)}(r,y)}{(1-L^{(2)}(r,y))^2}
\\
&= \int_{-\infty}^x \d y \frac{L^{(2)}(\zeta^{-1}(t),y)}{1-L^{(2)}(\zeta^{-1}(t),y)}
\end{align*}

\noindent which proves, with natural notation, that $\widehat L^{(2)}(t,y)=\frac{L^{(2)}(\zeta^{-1}(t),y)}{1-L^{(2)}(\zeta^{-1}(t),y)}$. Consequently,
$$
\zeta^{-1}(t)= \int_0^t  (1-L^{(2)}(\zeta^{-1}(r),\widehat Z^{(2)}_r))^2 {\rm d} r = \int_0^t  {{\rm d} r \over (1+\widehat L^{(2)}(r,\widehat Z^{(2)}_r))^2} .
$$

\noindent It implies $\zeta^{-1}$ is adapted with respect to the natural filtration $\widehat {\mathcal F}^{(2)}$ of $\widehat Z^{(2)}$. We also introduce for $x\in \r$,
\begin{equation}\label{eq:R}
R_t(x):=\int_{Z^{(2)}_{\zeta^{-1}(t)}}^x { {\rm d} y \over  1-L^{(2)}(\zeta^{-1}(t),y)} = \int_{\widehat Z^{(2)}_t}^x \left(1+ \widehat L^{(2)}(t,y) \right){\rm d} y.
\end{equation}

\noindent Observe that $x=\widehat Z^{(2)}_{t}$ if and only if $R_t(x)=0$.  Note that in the case $\mu>1$, the flow $(R_t(x),\, x\in \r)$ only gives the reconstruction of the burglar $Z_\cdot^{(2)}$ until the finite time $t_{\max}$.

\bigskip

\begin{theorem}\label{p:bassburdzy}
The process $(R_t(x),\, x\in \r)$ is a Bass--Burdzy flow for some parameters $\beta_1$ and $\beta_2$ and diffusivity $\sigma$. Moreover   $\beta_2=\beta_1+1$. 
\end{theorem}

\noindent {\it Proof}.  First we note that  a.s. for all $t\ge 0 $, 
\begin{equation}\label{eq:Ht}
H_t:= R_t(x)- x - \int_0^t {\bf 1}_{\{R_s(x) > 0\}} \d s, \qquad x\in \r, 
\end{equation}

\noindent is well defined. In other words, the right-hand side of \eqref{eq:Ht} does not depend on $x$.  This can be seen as follows: 
$$
\int_0^t {\bf 1}_{\left\{R_s(x) > 0\right\}} \d s
= \int_0^t {\bf 1}_{\left\{\widehat Z^{(2)}_{s} < x\right\}} \d s
=
\int_{-\infty}^x  \widehat L^{(2)}(t,y) \d y
$$

\noindent showing that the derivative with respect to $x$ of the right-hand side of \eqref{eq:Ht} vanishes.

 Clearly $(H_t)$ is  $(\widehat {\mathcal F}^{(2)}_t)$-adapted. We are going to show that $(H_t)$ is  an   $(\widehat {\mathcal F}^{(2)}_t)$-L\'evy process. As we shall see below,  this  boils down to showing the strong Markov property for $R_t$. We use equation \eqref{eq:XZ0} with ${\mathcal X}:=Z^{(2)}_{\zeta^{-1}(s)+\cdot}$, $f(r) = 1-L^{(2)}(\zeta^{-1}(s),r)$, $c=\widehat Z^{(2)}_s$, $c'=x$ and $ \zeta^{-1}(t+s)-\zeta^{-1}(s)$ in lieu of $s$  to get that
\begin{equation}\label{eq:integral}
\int_{R_s(x)}^{{\widetilde Z}^{(2)}_{\tilde t}} {\d r \over 1- {\widetilde L^{(2)}}(\tilde t,r)}
=
\int_{x}^{Z^{(2)}_{\zeta^{-1}(s+t)}} \frac{\d r}{1-L^{(2)}(\zeta^{-1}(s+t),r)}
\end{equation}

\noindent where we used the observation that $\eta_{f,c}(x)=R_s(x)$,    the process ${\widetilde Z}^{(2)}$ is defined by 
$$
{\widetilde Z}^{(2)} := \Upsilon(Z^{(2)}_{\zeta^{-1}(s)+\cdot},1-L^{(2)}(\zeta^{-1}(s),\cdot), Z^{(2)}_{\zeta^{-1}(s)})\equiv \Upsilon(Z^{(2)}_{\zeta^{-1}(s)+\cdot},1-L^{(2)}(\zeta^{-1}(s),\cdot)),
$$

\noindent and ${\widetilde L}^{(2)}$ is the local time of the process ${\widetilde Z}^{(2)}$,   $\tilde t = C_f(\zeta^{-1}(t+s)-\zeta^{-1}(s))$ with
$$
 C_f(v):=\int_{0}^v {\d u \over (1-L^{(2)}(\zeta^{-1}(s), Z^{(2)}_{\zeta^{-1}(s)+u}))^2}.
$$

\noindent We get
$$
\tilde t =\int_{\zeta^{-1}(s)}^{\zeta^{-1}(s+t)} {\d u \over (1-L^{(2)}(\zeta^{-1}(s), Z^{(2)}_u))^2}.
$$

\noindent We can rewrite \eqref{eq:integral} as
$$
\widetilde R_t(R_s(x))= R_{t+s}(x)
$$

\noindent with
$$
\widetilde R_t(y) = \int_{\widetilde Z^{(2)}_{\widetilde t}}^y {\d r\over 1 - \widetilde L^{(2)}(\widetilde t,r)}.
$$

\noindent Let $\widetilde \zeta$ denote the $\zeta$ associated with $\widetilde Z^{(2)}$, i.e.,
$$
{\widetilde \zeta}(v):=\int_0^v {\d u \over (1-\widetilde L^{(2)}(u,\widetilde Z^{(2)}_u))^2}. 
$$

\noindent Suppose for the time being that we know that 
\begin{equation}\label{eq:ttilde}
\widetilde t=\widetilde \zeta^{-1}(t).
\end{equation}

\noindent   By Remark \ref{R:markovfort} applied to the stopping time $\zeta^{-1}(s)$, $\widetilde Z^{(2)}$ is distributed as $Z^{(2)}$ and is independent of $\sigma\{Z^{(2)}_r, r\le \zeta^{-1}(s)\}$, thus independent of $\widehat {\mathcal F}_s^{(2)}$. It 
implies that conditionally on $\widehat {\mathcal F}_s^{(2)}$ and $R_s(x)=y$,  the process $(R_{t+s}(x))_{t\ge0}$ has the law of $(R_t(y))_{t\ge0}$.  It follows that conditionally on $\widehat {\mathcal F}_s^{(2)}$ and $R_s(x)=y$, for all $t\ge0$, $H_{t+s}-H_s= \widetilde  R_t(y)- y -\int_0^t 1_{\{\widetilde  R_s(y) > 0\}} \d s= \widetilde  H_t$, where and $\widetilde  H$ is defined from $\widetilde  R$ as $H$ is from $R$. In other words, $H_{t+s}-H_s, t\ge0$ is independent of $\widehat {\mathcal F}_s^{(2)}$ and has the same distribution as $H$. Then  $H$ is a (continuous) L\'evy process hence of the form $\sigma \gamma_t+ \beta_1 t$ with $\gamma$ an  $(\widehat {\mathcal F}_t^{(2)})$-Brownian motion, which yields the theorem.

It remains to prove \eqref{eq:ttilde}. We know from equation \eqref{eq:locZ} that
$$
\widetilde L^{(2)}(u,\widetilde Z^{(2)}_u)
=
{L^{(2)}( C_f^{-1}(u),Z^{(2)}_{\zeta^{-1}(s)+ C_f^{-1}(u)} ) - L^{(2)}(\zeta^{-1}(s),Z^{(2)}_{\zeta^{-1}(s)+ C_f^{-1}(u)} )
\over
1 - L^{(2)}(\zeta^{-1}(s),Z^{(2)}_{\zeta^{-1}(s)+ C_f^{-1}(u)})
}.
$$

\noindent Hence, by definition of $\widetilde \zeta$ and a change of variables,
$$
\widetilde \zeta(\tilde t)
=
\int_0^{\tilde t} {\d u \over (1-\widetilde L^{(2)}(u,\widetilde Z^{(2)}_u))^2}
=
\int_{\zeta^{-1}(s)}^{\zeta^{-1}(t+s)} { \d u \over (1- L^{(2)}(u,Z^{(2)}_{u}))^2} =t
$$

\noindent by definition of $\zeta$, which completes the proof. $\Box$

\bigskip

  We identify the parameters of the Bass--Burdzy flow in the following proposition. 

\begin{proposition}\label{p:parameters}
The parameters of the Bass--Burdzy flow are $\sigma=1$, $\beta_1={\delta\over 2}-1$ and $\beta_2={\delta\over 2}$.
\end{proposition}

\noindent {\it Proof}. We first show that $\sigma=1$. Recall that $\widehat Z^{(2)}_t=Z^{(2)}_{\zeta^{-1}(t)}$. Since ${1\over t}\zeta(t) \overset{t\to 0}{\longrightarrow} 1$, one can deduce from Theorem \ref{t:scaling} that $\left( {1\over \sqrt{a}} \widehat Z^{(2)}_{at},\, t\ge 0 \right)$ converges in distribution as $a\downarrow 0$ to a standard Brownian motion. Therefore, the hitting time of $x$ by $\widehat Z^{(2)}$, divided by $x^2$, converges in law when $x\downarrow 0$ to the hitting time of $1$ by a standard Brownian motion. But the hitting time of $x$ by $\widehat Z^{(2)}$ is the hitting time of $0$ by $R_t(x)$. We deduce that necessarily $\sigma=1$. 

We now determine $\beta_2$ (or equivalently $\beta_1$ since we know that $\beta_2=1+\beta_1$).  By definition of the Bass--Burdzy flow, if we start the flow $R_t(x)$ from $x^*= \sup_{t\ge0}  \widehat Z^{(2)}_t$, then $R_t (x^*)\ge 0$ for all $t\ge0$ and $R_t(x^*)= x^*+  \gamma_t + \beta_2 t $. Note that the unique zero is attained when $\widehat Z^{(2)}_t$ reaches its maximum $x^*$, that means $$x^*= - \inf_{t\ge0} (\gamma_t + \beta_2 t).$$

\noindent It implies that $\beta_2> 0$ and $x^*$ is exponentially distributed with parameter $2\beta_2$. On the other hand, observe that $x^*=\sup_{t\ge 0} Z^{(2)}_t$, so that,  recalling from Section \ref{s:burglar} that $\Xi_r=L^{(1)}(\infty,r)=1-L^{(2)}(\infty,r)$,
$$
x^*:=\inf\{r\ge 0\,:\, \Xi_r=1\}.
$$

\noindent In other words, by Proposition \ref{p:eve}, it is the hitting time of $0$ by a Jacobi process with parameters $(0,\delta)$ with initial distribution beta($1,{\delta\over 2}$).  The proof will be complete once we prove that this hitting time is exponentially distributed with parameter $\delta$. Let us prove it. Under some probability $P^x$, denote by $V$ a Jacobi($0,\delta$) process starting from $x$, and let $T_0$ be its hitting time of $0$.  We can check that $M_t:={V_t \over x} e^{\delta t}$ is a martingale. We can then define the probability measure $Q^x$ with Radon--Nikodym derivative $M_t$ with respect to $P^x$ on the $\sigma$-algebra $\sigma(V_s,\,s\in [0,t])$. Under $Q^x$, the process $V$
 is a Jacobi($4,\delta$) process. We have $ P^x (T_0> t)= e^{- \delta t} E_{Q^x}\left[\frac{x}{V_t}\right]$. Taking  $x$ with the  beta$(1, \frac\delta2)$ distribution, and setting $P:=\int_0^1  \frac\delta2\, (1-x)^{\frac\delta2-1} P^x \d x $, we get 
 \begin{eqnarray*}
  P (T_0> t)
  &=&
   e^{- \delta t} \int_0^1 \frac\delta2\, (1-x)^{\frac\delta2-1} \,   E_{Q^x}\left[\frac{x}{V_t}\right] \d x \\ 
   &=&
   e^{- \delta t} \int_0^1 \frac\delta2\, (1-x)^{\frac\delta2-1}x \,   E_{Q^x}\left[\frac{1}{V_t}\right] \d x.
\end{eqnarray*}

\noindent Since beta$(2, \frac\delta2)$ is the invariant distribution of  a Jacobi($4, \delta$) process, one has 
$$
\int_0^1 \frac\delta2\, (1-x)^{\frac\delta2-1}x \,   E_{Q^x}\left[\frac{1}{V_t}\right] \d x
=
\int_0^1 \frac\delta2\, (1-x)^{\frac\delta2-1}x \,   E_{Q^x}\left[\frac{1}{V_0}\right] \d x = 1.
$$

\noindent Hence $P(T_0>t)=e^{-\delta t}$ so $T_0$ is indeed exponentially distributed with parameter $\delta$. $\Box$

\bigskip

\noindent {\bf Remark}.   In light of Theorem III, Theorem \ref{t:perkins} gives a flow version of the Ray--Knight theorems appearing in \cite{hw00}.

\subsection{Contour function of a Fleming--Viot forest}
\label{s:FV}

Let $\delta>0$ and $b>0$. Let $X$ be the PRBM associated with $\mu= \frac2\delta$, which we consider up to time $\tau_b^0$ when the local time at position $0$ hits $b$.

\begin{notation}\label{n:<h}  
We define for $u\in (t_0^*,\tau_b^0)$ and $t\in (- \int_{t_0^*}^0 {\bf 1}_{\{ X_s >  0\} } {\rm d} s, \int_0^{\tau_b^0} {\bf 1}_{\{ X_s >  0\} } {\rm d} s ] $, 
$$
A_u^{+} := \int_{0}^u {\bf 1}_{\{ X_s >  0 \}} {\rm d} s,\, 
\qquad  
\alpha_t^{+}  := \inf\{u \,:\, A_u^{+} > t\}, 
\qquad
 X^{+}_t :=  X_{\alpha_t^{+}}.
$$

%\noindent We will call  ${\mathcal X}^{+}$ the process ${\mathcal X}$ seen above $0$. 
\end{notation}

 The process $X^{+}$ is the process $X$ looked above zero. The processes $(X^+_t,\,t\ge 0)$ and $(X^+_t,\,t\le 0)$ are independent, the former being a reflecting Brownian motion, see \cite{perman-werner}.
 We let as before
\begin{equation}\label{def:f}
f(r):=L(\tau_b^0,r), \qquad r\in [0,{\mathfrak d}_f)
\end{equation} 

\noindent be the local time at position $r$ at time $\tau_b^0$, where ${\mathfrak d}_f$ denotes the hitting time of $0$ by $L(\tau_b^0,r)$, $r\ge 0$ (which is $X_{t_0^*}$ if $|t_0^*|<\infty$ and $\infty$ otherwise).  We want to reconstruct $X^{+}$ conditionally on $f$ up to a certain random time. Recall the definition of the forward local time flow in Section \ref{s:process}. Let $S^+$ denote that of $X^{+}$, i.e. by definition
$$
S_{r,x}^+(a)  =L_{X^{+}}(\tau^r_a(X^{+}), x) = L(\tau_a^r,x),\qquad  0\le r\le x < {\mathfrak d}_f,\, a\in [0, f(r)]. 
$$

\noindent  By Theorem \ref{t:RK}  and Definition \ref{d:BESQdd'},  $S^+$ is a ${\rm BESQ}^{\delta,0}_b$ flow. In the notation of Section \ref{s:process}, let $Z^+:=\Upsilon((X^{+}_t,\, t\in [0,\tau_b^0(X^+)]),f)$.  Specifically, 
\begin{equation}\label{eq:Z+}
Z^{+}_{C_f(t)} = \int_0^{X^{+}_t} \frac{\d r}{f(r)}, \qquad
C_f(t):=\int_0^t \frac{\d s}{f(X^{+}_s)^2},
\qquad
 0\le t \le \tau_b^0(X^+).
\end{equation}

\noindent  The flow $Y^+:=\Psi(S^+,f)$ is a ${\rm Jacobi}(\delta,0)$ flow independent of $f$ by Theorem \ref{t:perkins}.  In particular, $Z^+$ can be thought of as the contour function of the Fleming--Viot forest embedded in $Y^+$, rooted at level $0$. The following theorem is the analog of Theorem \ref{p:loctime} (iii).

\begin{theorem}\label{t:RK1}
The process $Z^+$ is independent of $f$.
\end{theorem}

\noindent {\it Proof}. The forward local time flow of $Z^+$ is composed of the flow lines of $Y^+$ located at the right of $Y^+_{0,\cdot}(0)$. The flow  $Y^+$ is  independent of $f$, and $Z^+$ is measurable with respect to its forward local time flow. $\Box$

\bigskip

Let us give now a description in time of the process $Z^+$. We introduce a variant of the Bass--Burdzy flow. Let $\gamma^{+,\delta}$ be a reflecting Brownian motion with drift  $1-\frac{\delta}2$, see \cite{GS00}.  It is the absolute value of the unique strong solution of the SDE
$$
\d U_t = \Big(1-\frac{\delta}2\Big) {\rm sgn} U_t\,  \d t + \d B_t,\qquad U_0=0
$$

\noindent where $B$ is a standard Brownian motion. When $\delta = 2$, it is the usual reflecting Brownian motion. It is recurrent when $\delta \ge 2$, and transient when $\delta\in (0,2)$. It is proved in \cite{GS00} that $\gamma^{+,\delta}:=|U|$ is distributed as the process $(\sup_{s\in [0,t]} (B_s+(\frac{\delta}2-1) s) - (B_t+(\frac{\delta}2-1) t) ,\, t\ge 0)$.

One defines for any $x\ge 0$, the process $({\mathcal R}_t^{+,\delta}(x),\, t\ge 0)$ which is adapted to the filtration of $\gamma^{+,\delta}$ and is  the strong  solution of the SDE \footnote{To obtain the existence and uniqueness of the (strong) solution, we may use Zvonkin's method: Let   $h(x):= \int_0^x e^{2 \max(y,0)} \d y, x \in \r,$ be an increasing (and convex) function on $\r$. Applying It\^{o}--Tanaka's formula to $h({\mathcal R}_t^{+,\delta}(x))$, we see that $\eta_t:=h({\mathcal R}_t^{+,\delta}(x))$  satisfies the SDE: $\eta_t= h(x) + \int_0^t \sigma(\eta_s) \d \gamma_s^{+,\delta}$ with $\sigma(x):= - h'(h^{-1}(x))= - 1 - 2 \max(x,0)$.  In particular $\sigma$ is a Lipschitz function, we may apply Theorem V.6 in Protter \cite{protter} to get the existence and  uniqueness of the (strong) solution  $\eta$ and then that of   ${\mathcal R}^{+,\delta}(x)$.}
 $$
{\mathcal R}_t^{+,\delta}(x) = x -  \gamma_t^{+,\delta}  + \int_0^t {\bf 1}_{\{{\mathcal R}_s^{+,\delta}(x)> 0\}} \d s.
$$

\noindent Theorem \ref{t:FV} will construct such a solution. Let $(R_t^+(x),\, x\ge 0,\, t\ge 0)$ be the flow associated with $Z^{+}$ via \eqref{eq:zeta} and \eqref{eq:R}, i.e. set
\begin{equation}\label{eq:zeta+}
\zeta^{+}(t):=\int_0^t {{\rm d} s \over (1-L^+(s,Z^+_s))^2}, \qquad \widehat Z^+_t :=  Z^+_{(\zeta^+)^{-1}(t)}
\end{equation}
 
\noindent where $L^+$ denotes the local times of $Z^+$, and for $x\ge 0$,
\begin{equation}\label{eq:R+}
R_t^+(x):=\int_{\widehat Z^{+}_t}^x { {\rm d} r \over  1-L^+((\zeta^+)^{-1}(t),r)}=\int_{\widehat Z^{+}_t}^x \left(1+ \widehat L^{+}(t,y) \right){\rm d} y.
\end{equation} 

\noindent In the (second) equality,  $\widehat L^{+}$ is the local time of $\widehat {Z}^+$. This equality is proved along the lines of \eqref{eq:R}. Notice that $R^+$ is measurable with respect to $X^+$.

With the notation  $t_{\max}^+ := \inf\{s\ge 0\,:\,L^{+}(s,Z^{+}_s)=1\}\in(0,\infty]$,  we have the analog of Lemma \ref{l:tmax}, meaning that  
$\lim_{t\to t_{\max}^+} \zeta^+(t)=\infty$. The time $t_{\max}^+$  is finite almost surely if $\delta\in(0,2)$ and infinite otherwise. The proof follows the lines of the proof of Lemma \ref{l:tmax}. We need to replace there $B$ by $(X^+_t, \, t\in[0,\tau_b^0(X^+)])$, $T_g^B$ by $\tau_b^0(X^+)$, and $(h(x),\, x\ge 0)$ is the local time at $x$ of $(X^+_t,\, t\le 0)$, i.e. a ${\rm BESQ}^\delta_0$ process.      
%In the case $\delta\ge 2$, $\zeta^{+}(t)<\infty$ for all $t\ge 0$. In the case $\delta \in (0,2)$, $\zeta^+(t)$ blows up in finite time. In this case, the process $\widehat Z^{+}$ explores only a part of the Fleming--Viot trees rooted at $0$.  
\begin{theorem}\label{t:FV}
The process $(R_t^+(x),\, x\ge 0,\, t\ge 0)$ has the law of $({\mathcal R}_t^{+,\delta}(x),\, x\ge 0,\, t\ge 0)$. In particular, $(\widehat Z^+_t,\,t\ge 0)$ is distributed as $({\mathcal R}^{+,\delta}_t)^{-1}(0)$.
\end{theorem}

\noindent {\it Proof}. We observe that $(X_t^+,\, t\in [0,\tau_b^0(X^+)],f)$ and $(X_{\tau_b^0(X^+)-t}^+,\, t\in [0,\tau_b^0(X^+)],f)$ have the same distribution. Therefore, we will take without loss of generality $Z^+=\Upsilon((X_{\tau_b^0(X^+)-t}^+,\, t\in [0,\tau_b^0(X^+)]),f)$ (and take $\widehat Z^+$, $\widehat L^+$, $R^+$ as in \eqref{eq:zeta+}, \eqref{eq:R+} associated with this $Z^+$).

As in the proof of Theorem \ref{p:bassburdzy}, the process
$ R_t^+(x)- x - \int_0^t {\bf 1}_{\{R_s^+(x) > 0\}} \d s$ does not depend on $x$, since $
\int_0^t {\bf 1}_{\left\{R_s^+(x) > 0\right\}} \d s
= \int_0^t {\bf 1}_{\left\{\widehat Z^{+}_{s} < x\right\}} \d s
=
\int_{-\infty}^x  \widehat L^{+}(t,y) \d y$. Therefore, it suffices to show that $(R_t^+(0),\, t\ge 0)$ is distributed as $(-\gamma_t^{+,\delta},\, t\ge 0)$ (observe that by definition, $R_t^+(0)\le 0$ for all $t$).

In Section \ref{s:burglar}, taking $b=1$ in  equation \eqref{def:Z12}, we constructed the process $Z^{(2)}$ associated with
$$
X^{(2)}_s:= X_{\tau_1^0-s},\, s\in [0,\tau_1^0-t^*),
$$

\noindent where $t^*$ is the time such that $X_{t^*}=\inf_{[0,\tau_1^0]} X$. We proved in Proposition \ref{p:parameters} that the time-changed process $\widehat Z^{(2)}$, defined in equation \eqref{eq:zeta}, is associated with the Bass--Burdzy flow $(R_t(x))_{t,x}$ with parameters $\beta_1=\frac{\delta}2-1$ and $\beta_2=\frac{\delta}2$.

We define the process $X^{(2),+}$ as the process $X^{(2)}$ looked above $0$ (which should be denoted by $(X^{(2)})^{+}$ in the notation \ref{n:<h}). Specifically, $X^{(2),+}_t=X^+_{\tau_1^0(X^+)-t}$, $t\in [0,\tau_1^0(X^+) - \tau_\ell^0(X^+)]$ where $\ell:=L(t^*,0)$. We now introduce the process $Z^{(2),+} := \Upsilon(X^{(2),+},f)$  where $f$ is given by \eqref{def:f} with $b=1$. We mention that $Z^{(2),+}$ is also the process $Z^{(2)}$ looked above $0$. [We omit the details. Intuitively, each excursion of $Z^{(2),+}$ is obtained from the corresponding excursion of $X^{(2)}$ in $\r_+$ by the transformation $\Upsilon(\cdot,f)$.] Analogously to \eqref{eq:zeta+}, we let
$$
\widehat Z^{(2),+}_t :=  Z^{(2),+}_{(\zeta^{(2),+})^{-1}(t)},\, t < \int_0^\infty {\bf 1}_{\{\widehat Z^{(2)}_s>0\}} \d s=:T^+
$$
 
\noindent where $\zeta^{(2),+}(t):=\int_0^t {{\rm d} s \over (1-L^{(2),+}(s,Z^{(2),+}_s))^2}$, $L^{(2),+}$ being the local times of $Z^{(2),+}$. Again, the process $\widehat Z^{(2),+}$ is the process $\widehat Z^{(2)}$ looked above $0$.

We claim that for any fixed $M>0$, \begin{equation}\label{T+M}
\p\big( T^+ \le M \, |\, L(t^*, 0) =\ell \big) \to 0, \qquad \ell \to 0.
\end{equation}

%\noindent In fact by time-change, $T^+$ can be represented by the local times of  the $\mu$-process $X$:  \begin{equation}\label{T+} T^+= \int_{t^*}^{\tau_1^0} \frac{1_{\{ X_s > 0\}}}{ (L(s, X_s))^2}{\rm d} s .\end{equation}

\noindent Let us prove \eqref{T+M}. First, By Lemma \ref{l:tmax}, if $t_{\max}<\infty$ and $Z^{(2)}_{t_{\max}}>0$, then $T^+=\infty$. Suppose now that $t_{\max}=\infty$ or $t_{\max}<\infty$ and $Z^{(2)}_{t_{\max}}\le 0$. By \eqref{eq:zeta}, we have $$T^+=\int_0^{t_{\max}}    \frac{1_{\{Z^{(2)}_s >0\}}}{(1- L^{(2)}(s, Z^{(2)}_s))^2} \d s .
$$

\noindent Recall from  \eqref{eq:locZ} that $L^{(2)}(t,x)= \frac{L_{X^{(2)}}(C_f^{-1}(t), \eta_f^{-1}(x))}{f(\eta_f^{-1}(x))}$, with $C_f(t):=\int_0^t \frac{\d s}{f(X^{(2)}_s)^2}$, $\eta_f(x):=\int_0^x \frac{\d r}{f(r)}$.  Therefore $$T^+= \int_0^{C_f^{-1}(t_{\max})} \frac{1_{\{ X^{(2)}_r > 0\}}}{ (f(X^{(2)}_r)- L_{X^{(2)}}(r, X^{(2)}_r))^2}{\rm d}r
=
 \int_{\tau_1^0-C_f^{-1}(t_{\max})}^{\tau_1^0} \frac{1_{\{ X_s > 0\}}}{ (L(s, X_s))^2}{\rm d} s
.$$

\noindent If $t_{\max}=\infty$ (i.e. $\delta\ge 2$), $\tau_1^0-C_f^{-1}(t_{\max})=t^*$. If $t_{\max}<\infty$ and $Z^{(2)}_{t_{\max}}\le 0$, the process $Z^{(2)}$ stays negative after $t_{\max}$, which is equivalent with saying that $X$ is always negative between $t^*$ and $\tau_1^0-C_f^{-1}(t_{\max})$. Hence, in both cases,  $T^+=\int_{\tau_1^0-C_f^{-1}(t_{\max})}^{\tau_1^0} \frac{1_{\{ X_s > 0\}}}{ (L(s, X_s))^2}{\rm d} s= \int_{t^*}^{\tau_1^0} \frac{1_{\{ X_s > 0\}}}{ (L(s, X_s))^2}{\rm d} s$. From the occupation times formula, we get that
$$
T^+ = \int_0^\infty \bigg(\frac{1}{L(t^*,x)} - \frac{1}{f(x)}\bigg)\d x \ge \int_0^1 \bigg(\frac{1}{L(t^*,x)} - \frac{1}{f(x)}\bigg)\d x
$$

\noindent (the upper boundary $1$ is arbitrary). Using \cite{perman} Lemma 2.3, conditionally on $L(t^*,0)=\ell$, the processes $x\in \r_+\to L(t^*,x)$  and $x\in\r_+\to f(x)-L(t^*,x)$ are independent, with distribution respectively ${\rm BESQ}^{\delta}_{\ell}$ and ${\rm BESQ}^{0}_{1-\ell}$. The integral $\int_0^1 \frac{1}{f(x)} \d x$ is tight as $\ell\to 0$ while the integral $\int_0^1 \frac{1}{L(t^*,x)} \d x$ tends in law to $\infty$ by scaling. We  deduce \eqref{T+M}.

\bigskip

As in \eqref{eq:R+}, we define for $x\ge 0$,
$$
R_t^{(2),+}(x):=\int_{\widehat Z^{(2),+}_t}^x { {\rm d} r \over  1-L^{(2),+}((\zeta^{(2),+})^{-1}(t),r)}=\int_{\widehat Z^{(2),+}_t}^x \left(1+ \widehat L^{(2),+}(t,y) \right){\rm d} y.
$$ 

\noindent By construction, $\widehat Z^{(2),+}$ is the process $(\widehat Z^+_t,\, t\in [0,T^+))$, hence for all $x\ge 0$, the process $R^{(2),+}(x)$ coincides with $R^+(x)$ up to time $T^+$.

\bigskip

On the other hand, the process $R^{(2),+}(x)$ is deduced from  $R(x)$ by a time-change: define   $$
\tau^{(2),+}_t := \inf\Big\{s\ge 0\,:\, \int_0^s {\bf 1}_{\{ R_u(0)<0 \}}\d u > t\Big\}=\inf\Big\{s\ge 0\,:\, \int_0^s {\bf 1}_{\{ \widehat Z_u^{(2)} >0 \}}\d u > t\Big\}.
$$

\noindent Then we have for any $x\ge 0$ and $ 0\le t  < T^+$, \footnote{In fact,   $\widehat Z^{(2),+}_{t} = Z^{(2),+}_{(\zeta^{(2),+})^{-1}(t)} = Z^{(2)}_{\alpha^{(2),+} \circ (\zeta^{(2),+})^{-1}(t)}$ where $\alpha^{(2),+}$ is the inverse of $\int_0^u {\bf 1}_{\{Z^{(2)}_s>0\}}\d s$. We notice that $L^{(2),+}(s,x)=L^{(2)}(\alpha^{(2),+}_s,x)$, $x\ge 0$.  Hence $R_t^{(2),+}(x)= R_{\zeta^{(2)}\circ \alpha^{(2),+} \circ (\zeta^{(2),+})^{-1}(t)}(x)$. We check $\zeta^{(2)}\circ \alpha^{(2),+} \circ (\zeta^{(2),+})^{-1}(t)= \tau^{(2),+}_t$ by noticing that $
\int_0^{\zeta^{(2)}\circ \alpha^{(2),+} \circ (\zeta^{(2),+})^{-1}(t) } {\bf 1}_{\{\widehat Z_u^{(2)} >0\}} \d u 
=
\int_0^{\alpha^{(2),+} \circ (\zeta^{(2),+})^{-1}(t) } {\bf 1}_{\{ Z_s^{(2)} >0\}} \frac{\d s}{(1- L^{(2)}(s,Z_s^{(2)})  )^2} 
=
\int_0^{(\zeta^{(2),+})^{-1}(t) } \frac{\d s}{(1- L^{(2),+}(s,Z_s^{(2),+})  )^2} 
=
t. $ }
\begin{equation}\label{eq:R2+}
R^{(2),+}_t(x) = R_{\tau_t^{(2),+}}(x).
\end{equation}

\noindent Recall that we want to show that $(R_t^+(0),\, t\ge 0)$ is distributed as $(-\gamma_t^{+,\delta},\, t\ge 0)$. If we let $\gamma$ be the Brownian motion driving the Bass--Burdzy flow $(R_t(x))_{t,x}$, we have
$$
R_t(0) =  \gamma_t +\frac{\delta}2 t - \int_0^t {\bf 1}_{\{R_s(0)<0\}} \d s.
$$

\noindent By Tanaka's formula,
$$
\min(R_t(0),0) = \int_0^t {\bf 1}_{\{R_s(0)\le 0\}} \d R_s(0) - \frac12 L_{R(0)}(t,0)
$$

\noindent where $L_{R(0)}(t,0)$ is the local time of $R_{\cdot}(0)$ at time $t$ at position $0$. By \eqref{eq:R2+}, observing that $R^{(2),+}_t(0)\le 0$, 
$$
R^{(2),+}_t(0) = \int_0^{\tau_t^{(2),+}} {\bf 1}_{\{R_s(0)\le 0\}} \d R_s(0) - \frac12 L_{R(0)}(\tau_t^{(2),+},0).
$$

\noindent We have $L_{R(0)}(\tau_t^{(2),+},0)= L_{R^{(2),+}(0)}(t,0)$ which is the local time of $R_{\cdot}^{(2),+}(0)$ at time $t$ at position $0$. Using that $ \int_0^{\tau_t^{(2),+}} {\bf 1}_{\{R_s(0)\le 0\}} \d s = \int_0^{\tau_t^{(2),+}} {\bf 1}_{\{R_s(0)< 0\}} \d s =t$, we get
$$
R^{(2),+}_t(0) = \int_0^{\tau_t^{(2),+}} {\bf 1}_{\{R_s(0)\le 0\}} \d \gamma_s +(\frac{\delta}2 -1)t  - \frac12 L_{R^{(2),+}(0)}(t,0),\, 0\le t <T^+.
$$

\noindent By Dambis--Dubins--Schwarz theorem, there exists a  standard Brownian motion $B^{(+)}=(B^{(+)}_t,\,t\ge 0)$  such that 
$$
\int_0^{\tau_t^{(2),+}} {\bf 1}_{\{R_s(0)< 0\}} \d \gamma_s = B^{(+)}_t,\, 0\le t <T^+.
$$

\noindent By Skorokhod's lemma, 
$$
\frac12 L_{R^{(2),+}(0)}(t,0) = \sup_{s\le t} \bigg(B_s^{(+)} +\bigg(\frac{\delta}2 -1\bigg)s\bigg),\, 0\le t<T^+. 
$$

\noindent Define 
$$
\gamma^{+,\delta}_t :=  \sup_{s\le t} \bigg(B_s^{(+)} +\bigg(\frac{\delta}2 -1\bigg)s\bigg) -B^{(+)}_t -\bigg(\frac{\delta}2-1\bigg) t  ,\, t\ge 0.
$$

\noindent The process $-\gamma^{+,\delta}$ coincides with $R^{(2),+}(0)$ up to time $T^+$.  We now prove that $R^+(0)$ is distributed as $-\gamma^{+,\delta}$. Fix $M>0$ and a functional $F$ measurable and bounded. Notice that $L(t^*,0)$ is measurable with respect to the process $(X_s,\, s\in [0,\tau_1^0])$ looked below $0$, which is independent of $X^+$, see \cite{perman-werner}. Therefore  $R^+$ and $L(t^*,0)$ are independent, which implies that for all $\ell>0$,
$$
E[F(R^+_t(0),\, t\in [0,M])]
=
E[F(R^+_t(0),\, t\in [0,M])\, |\, L(t^*,0) \le \ell]
.
$$

\noindent We consider two cases. The first case is when $T^+>M$. In this case, $R^+_t(0)=R^{(2),+}_t(0)= -\gamma_t^{+,\delta}$, for $t\in [0,M]$. The second case is when $T^+\le M$, which has probability going to $0$ under $\p(\cdot \,|\, L(t^*,0) \le \ell)$ when $\ell \to 0$ by \eqref{T+M}. We deduce that 
$$
E[F(R^+_t(0),\, t\in [0,M])]
=
E[F(-\gamma_t^{+,\delta},\, t\in [0,M])]
$$
 
\noindent which is what we wanted to prove. $\Box$

\appendix
\section{Perfect flow property}
\label{s:perfectflow}

Recall the definition of a ${\rm BESQ}^\delta$ in Definition \ref{d:BESQflow}. We establish some flow properties for the ${\rm BESQ}^\delta$ flow.

\begin{proposition}
Let  ${\mathcal S}$ be either a non-killed ${\rm BESQ}^\delta$ flow with $\delta>0$ or  a killed ${\rm BESQ}^\delta$ flow with  $\delta\le 0$. The flow ${\mathcal S}$ satisfies the {\it perfect flow property}: 
almost surely, for every $ r \le x \le  y$, ${\mathcal S}_{r, y}= {\mathcal S}_{x,y} \circ {\mathcal S}_{r,x}$.
\end{proposition}

\noindent {\it Proof}. The property is true when one of the inequalities in $r\le x\le y$ is an equality since ${\mathcal S}_{z,z}(a)=a$. Let $a\ge 0$ and $r<x<y$. We first treat the case $\delta>0$. By Proposition \ref{p:Sembedded}, we can set ${\mathcal S} = (S_{r,x}(a),\, x\ge r)_{a,r\in\r}$.   We want to show that $S_{r,y}(a)=S_{x,y}\circ S_{r,x}(a)$.  Let $b:=S_{r,x}(a)$. By definition of $S$, $L(\tau_a^r,x)=b$ hence $\tau_b^x \ge \tau_a^r$.  If $t>\inf\{s>\tau_a^r\,:\, X_s>x\}=:\theta$, then $L(t,x)>b$ since a Brownian motion accumulates local time at any level that it visits. It implies that $\tau_b^x\le \theta$. If $t\in (\tau_a^r,\theta)$, one has $L(t,x)=L(\tau_a^r,x)=b$. We conclude that  $L(\tau_a^r,y)=L(\tau_b^x,y)$ indeed (and $\tau_b^x=\theta$). We now deal with the  case $\delta\le 0$. To be consistent with our setting, we will actually  consider a ${\rm BESQ}^{\delta'}$ flow with $\delta'\le 0$, and rather take $\delta:=2-\delta'$ which is now greater than $2$. Then we can set ${\mathcal S}={(S_{-r,-x}(a))_{r\le x, a\ge0}}$. 
We write  $S_{-r,-x}(a)=L(\tau_a^{-r},-x)=:b$  and see that $\tau_b^{-x} \ge \tau_a^{-r}$.  Again, if $t>\inf\{s>\tau_a^{-r}\,:\, X_s<-x\}$, then $L(t,-x)>b$ because the PRBM does not have point of monotonicity when $\delta\ge 2$, see \cite{perman-werner}. We deduce that $L(\tau_a^{-r},-y)=L(\tau_b^{-x},y)$. $\Box$

\bigskip

% The Jacobi($\delta,\delta'$) flow also possesses this property. We  will deduce it from the disintegration theorems of Section \ref{s:perkins}.  
When $\delta\in (0,2)$, the killed ${\rm BESQ}^\delta$ flow  loses this property.  There will be exceptional times when a flow line starting at $0$ will exit the boundary $0$. This comes from the fact that in the case $\delta\in (0,2)$, the PRBM has points of monotonicity, see \cite{perman-werner}. Still the flow satisfies a weaker form of the flow property stated in Proposition \ref{p:almost}. We need some preliminary results.

\bigskip

The following proposition shows that, for $\delta\in(0,2)$ the killed ${\rm BESQ}^\delta$ flow is naturally embedded in its non-killed version. 
\begin{proposition}\label{p:Stilde}
Let $\delta\in (0,2)$ and ${\mathcal S}=({\mathcal S}_{r,x}(a),\, x\ge r)_{r \in \r,a\ge 0}$ be a non-killed ${\rm BESQ}^\delta$ flow. 
Define a flow $\Pi(\mathcal S):=(\widetilde {\mathcal S}_{r,x}(a),\, x\ge r)_{r\in\r,a\ge 0} $ where for $r\in \r$ and $a\ge 0$,  $\widetilde {\mathcal S}_{r,x}(a)$ is equal to ${\mathcal S}_{r,x}(a)$ for $x\ge r$ up to the time $\inf\{r'> r\,:\, {\mathcal S}_{r,r'}(a)=0\}$, and is equal to $0$ afterwards. Then $\Pi({\mathcal S})$ is a killed ${\rm BESQ}^\delta$ flow. We call  $\Pi({\mathcal S})$ the killed version of ${\mathcal S}$. 
\end{proposition}

\noindent {\bf Remark}. Note that when $a=0$, at times $r$ when the flow line $({\mathcal S}_{r,x}(0),\, x\ge r)$ starts an excursion away from $0$, the flow line $(\widetilde {\mathcal S}_{r,x}(0),\, x\ge r)$ traces the same excursion then gets absorbed at $0$ when coming back.

\bigskip

\noindent{\it Proof}.  The flow lines of $\Pi ({\mathcal S})$ have the required distribution by construction so we need to check the regularity of $\Pi ({\mathcal S})$ imposed in Definition \ref{d:BESQflow}.  Assumption (i) is clear. We prove (ii).  The map $a\mapsto \widetilde {\mathcal S}_{r,x}(a)$ is nondecreasing. Indeed, if $\widetilde {\mathcal S}_{r,x}(a)>0$, then $\widetilde {\mathcal S}_{r,x}(a)= {\mathcal S}_{r,x}(a)\le {\mathcal S}_{r,x}(a')=\widetilde {\mathcal S}_{r,x}(a')$ for all $a'\ge a$. And obviously $\widetilde {\mathcal S}_{r,x}(a)\le \widetilde {\mathcal S}_{r,x}(a')$ if $\widetilde {\mathcal S}_{r,x}(a)=0$. We show the right-continuity. The case  $\widetilde {\mathcal S}_{r,x}(a)>0$ is immediate from the same arguments, using the right-continuity of  $a\mapsto{\mathcal S}_{r,x}(a)$. Suppose now that $\widetilde {\mathcal S}_{r,x}(a)=0$. We can suppose $x>r$. By construction, there exists $r'\in (r,x]$ such that ${\mathcal S}_{r,r'}(a)=0$. Since $b\mapsto \mathcal S_{r,r'}(b)$ is piecewise constant, we would have $\mathcal S_{r,r'}(b)=0$ for $b>a$ small enough, which implies ${\mathcal S}_{r',r}(b)=0$ for the same set of $b$, hence $\widetilde {\mathcal S}_{r,r'}(b)=0$. The right-continuity is therefore proved. Finally, we check condition (iii).  We suppose  $\widetilde {\mathcal S}_{r',r}(a')>a$.  
In particular, $\widetilde {\mathcal S}_{r',r}(a')>0$ so  $\widetilde {\mathcal S}_{r',r}(a')={\mathcal S}_{r',r}(a')$ hence by construction $\widetilde {\mathcal S}_{r',x}(a')= {\mathcal S}_{r',x}(a') \ge {\mathcal S}_{r,x}(a) = \widetilde {\mathcal S}_{r,x}(a)$ for all $x$ before $\widetilde {\mathcal S}_{r,\cdot}(a)$ hits $0$ when it will get absorbed. The case $\widetilde {\mathcal S}_{r',r}(a')<a$ is dealt with similarly. $\Box$

\bigskip

Conversely, for $\delta\in(0,2)$, a non-killed ${\rm BESQ}^\delta$ flow can also be constructed from its killed version. Recall the definition of the dual flow in Proposition \ref{c:BESQdual}.

\begin{proposition}
Let $\delta\in (0,2)$ and $\widetilde {\mathcal S}$ be a killed ${\rm BESQ}^\delta$ flow. The dual of $\widetilde {\mathcal S}$ is a non-killed ${\rm BESQ}^{2-\delta}$ flow, from which we can construct its killed version by Proposition \ref{p:Stilde}. Let $ {\mathcal S}$ denote the dual of the latter. Then $ {\mathcal S}$ is a non-killed ${\rm BESQ}^\delta$ flow, and  $\widetilde {\mathcal S}$ is the killed version of $\mathcal S$. 
\end{proposition}

\noindent {\it Proof}. We start with a lemma.
\begin{lemma}
Let $\delta \in (0,2)$. Let ${\mathcal S}$  be a non-killed ${\rm BESQ}^\delta$ flow. Then the dual of ${\mathcal S}$ is the killed version of the dual of $\Pi({\mathcal S})$, i.e. ${\mathcal S}^*=\Pi(\Pi({\mathcal S})^*)$.
\end{lemma}

\noindent {\it Proof}. Let $\mathcal S^*$, resp. $\Pi(\mathcal S)^*$, denote the dual of ${\mathcal S}$, resp. $\Pi(\mathcal S)$. Write $\widetilde {\mathcal S}$ for $\Pi(\mathcal S)$. We have by definition, for any $r\le x$,
$$
\mathcal S^*_{r,x}(a) := \inf\{b\ge 0\,:\, {\mathcal S}_{-x,-r}(b)>a\},
\qquad 
(\widetilde {\mathcal S})^*_{r,x}(a) := \inf\{b\ge 0\,:\, \widetilde {\mathcal S}_{-x,-r}(b)>a\}
.
$$

\noindent Let $r < x$ and $a\ge 0$. We first show that $(\widetilde {\mathcal S})^*_{r,x}(a)\ge {\mathcal S}_{r,x}^*(a)$. Let $b>(\widetilde {\mathcal S})^*_{r,x}(a)$, hence $\widetilde {\mathcal S}_{-x,-r}(b)>a$. In particular, $\widetilde {\mathcal S}_{-x,-r}(b)>0$ so that $\widetilde {\mathcal S}_{-x,-r}(b)= {\mathcal S}_{-x,-r}(b)$. Therefore  $ {\mathcal S}_{-x,-r}(b)>a$ hence the inequality.

We suppose now that $\mathcal S^*_{r,s}(a)\neq 0$ for all $s\in (r,x)$, and we want to show that  $\mathcal S^*_{r,x}(a)=(\widetilde {\mathcal S})^*_{r,x}(a)$. We only have to prove that $\mathcal S^*_{r,x}(a)\ge (\widetilde {\mathcal S})^*_{r,x}(a)$. Let $b>\mathcal S^*_{r,x}(a)$. Then ${\mathcal S}_{-x,-r}(b)>a$. If ${\mathcal S}_{-x,-s}(b)=0$ for some $s\in (r,x)$, we would have by the perfect flow property ${\mathcal S}_{-s,-r}(0)>a$ hence $\mathcal S^*_{r,s}(a)= 0$, which is a contradiction. Therefore ${\mathcal S}_{-x,-s}(b)\neq 0$ for all $s\in (r,x)$, hence ${\mathcal S}_{-x,-r}(b)=\widetilde {\mathcal S}_{-x,-r}(b)$ so $\widetilde {\mathcal S}_{-x,-r}(b) >a$. We deduce that $\mathcal S^*_{r,x}(a)\ge (\widetilde {\mathcal S})^*_{r,x}(a)$.  

Finally, we show that  flow lines in ${\mathcal S}^*$ are absorbed at $0$. So we suppose that $\mathcal S^*_{r,s}(a)= 0$ for some $s\in (r,x)$, and we want to show that  $\mathcal S^*_{r,x}(a)=0$. We have by definition, ${\mathcal S}_{-s,-r}(b)>a$ for all $b>0$. Since $b\to {\mathcal S}_{-s,-r}(b)$ is piecewise constant, we must have ${\mathcal S}_{-s,-r}(0)>a$, therefore ${\mathcal S}_{-x,-r}(0)>a$ by the perfect flow property. It implies that $\mathcal S^*_{r,x}(a)=0$.  $\Box$

\bigskip

We go back to the proof of the proposition. We use the notation of the proposition.  By definition, ${\mathcal S}=(\Pi( ( {\widetilde {\mathcal S}})^*))^*$. By the lemma (using $({\widetilde {\mathcal S}})^*$, which is a non-killed ${\rm BESQ}^{2-\delta}$ flow in place of $\mathcal S$ in the statement of the lemma), we have $\Pi( \mathcal S) = ((\widetilde {\mathcal S})^*)^* = \widetilde {\mathcal S}$. The fact that  ${\mathcal S}$ is a non-killed ${\rm BESQ}^\delta$ flow comes from Proposition \ref{c:BESQdual}. $\Box$

\bigskip

We can now state a flow property for a general ${\rm BESQ}^\delta$ flow.

\begin{proposition}\label{p:almost}
Let $\delta\in \r$. A general ${\rm BESQ}^\delta$ flow $\mathcal S$ satisfies the following properties almost surely:
\begin{itemize}
\item (Almost perfect flow property) If $r\le x\le y$, $a\ge 0$ and ${\mathcal S}_{r,x}(a)>0$, then ${\mathcal S}_{x,y} \circ {\mathcal S}_{r,x}(a) = {\mathcal S}_{r,y}(a)$.
\item (Coalescence) If $r',r<x$, $0\le a,a'$ and ${\mathcal S}_{r,x}(a)={\mathcal S}_{r',x}(a')$, then ${\mathcal S}_{r,y}(a)={\mathcal S}_{r',y}(a')$ for all $y\ge x$.
\end{itemize}
\end{proposition}

\noindent {\it Proof}. These properties hold for non-killed ${\rm BESQ}^\delta$ flows and killed ${\rm BESQ}^\delta$ flows with $\delta\le 0$  as a consequence of the perfect flow property, so we should only deal with $\delta\in (0,2)$ and killed ${\rm BESQ}^\delta$ flows. Write it $\widetilde {\mathcal S}$ and let rather ${\mathcal S}$ be the non-killed version of it (i.e. the flow such that $\Pi(\mathcal S)=\widetilde {\mathcal S}$). Let $r\le x \le y$, $a\ge 0$ and suppose $b:=\widetilde {\mathcal S}_{r,x}(a)>0$. From the definition of $\Pi(\mathcal S)$, we have $\widetilde {\mathcal S}_{r,s}(a)={\mathcal S}_{r,s}(a)$ for all $s$ from $r$ to the hitting time of $0$, hence we can apply the perfect flow property of $\mathcal S$. After hitting $0$, the flow line is absorbed. Let now $r,r'<x$, $0\le a,a'$ and suppose that $b:=\widetilde {\mathcal S}_{r,x}(a)={\mathcal S}_{r',x}(a')$. If $b>0$, we can apply the almost perfect flow property to conclude that flow lines coalesce. If $b=0$, both flow lines are absorbed at $0$ (we use here the fact that $r,r' <x$ so that none are flow lines which are departing from $0$ to trace an excursion). $\Box$ 

\bigskip

{\bf Remark}. The coalescence property does not hold in full generality if we let $r=x$ for killed ${\rm BESQ}^\delta$ flows when $\delta \in (0,2)$. Indeed, at times when a flow line starts at $r$ an excursion away from $0$, flow lines which are already absorbed at $0$ will stay stuck at $0$, while the flow line starting at $r$ will leave $0$. 

\begin{corollary}\label{c:almost}
Let $\delta,\delta'\in \r$. A general {\rm Jacobi($\delta,\delta'$)} flow  satisfies the following properties almost surely:
\begin{itemize}
\item (Almost perfect flow property) If $r\le s\le t$, $v\in [0,1]$ and ${\mathcal Y}_{r,s}(v)>0$, then ${\mathcal Y}_{s,t} \circ {\mathcal Y}_{r,s}(v) = {\mathcal Y}_{r,t}(v)$.
\item (Coalescence) If $r,r'<s$, $v,v'\in [0,1]$ and ${\mathcal Y}_{r,s}(v)={\mathcal Y}_{r',s}(v)$, then ${\mathcal Y}_{r,t}(v)={\mathcal Y}_{r',t}(v)$ for all $t\ge s$.
\end{itemize}
\end{corollary}

\noindent {\it Proof}. Proposition \ref{p:almost} stays true for ${\rm BESQ}^{\delta,\delta'}_b$ flows by construction (for $a\in [0,f(r)]$ and $r,r'\in [0,{\mathfrak d}_f)$ there)  then we use Theorem \ref{t:perkins} to prove the corollary when $r,r'\ge 0$. It is true for all $r,r'$ by stationarity. $\Box$

\section{A Girsanov theorem for Jacobi flows}
\label{s:girsanov}

Let $\delta \in \r$. Consider a general Jacobi$(\delta,0)$ flow ${\mathcal Y}$,  driven by the martingale measure ${\mathcal M}$ given by Definition \ref{d:M} under some measure ${\mathbb P}$. Fix $T>0$, and restrict to the time interval $[0,T]$. Call ${\mathcal Y}^T$ the collection $({\mathcal Y}_{s,t},\, 0\le s\le t \le T)$. Let $\delta'\in \r$. We want to relate  the Jacobi($\delta,\delta'$) flow to the Jacobi($\delta,0$) flow. We will use in this section the notation ${\mathcal M}_s(A):={\mathcal M}(A\times [0,s])$ for $s\ge 0$ and $A$ a Borel set of $[0,1]$.

\bigskip

For ${\mathfrak v} \in [0,1]$, let ${\mathcal F}_{{\mathfrak v}}^T$ be the $\sigma$-field generated by $({\mathcal M}_s(A),\, s\in [0,T],\, A \hbox{ Borel set of } [0,{\mathfrak v}])$. Observe that ${\mathcal F}^T:= ({\mathcal F}^T_{\mathfrak v}, {\mathfrak v} \in [0,1))$ forms a filtration.

We define 
$$
m_{\mathfrak v}^T := {\rm e}^{-\frac{\delta'}{2(1-{\mathfrak v})}{\mathcal M}_T([0,{\mathfrak v}]) - \frac{T}{8} \frac{\delta'^2}{1-{\mathfrak v}} {\mathfrak v} },\qquad {\mathfrak v} \in [0,1).
$$

\begin{theorem}\label{t:girsanov}
The process $(m_{\mathfrak v},\, {\mathfrak v} \in [0,1))$ is a ${\mathcal F}^T$-martingale. Introduce the measure ${\mathbb Q}$ on ${\mathcal F}_1^T$ such that 
$$
\frac {\d {\mathbb Q}}{\d {\mathbb P}}_{\big| {\mathcal F}_{\mathfrak v}^T} := m_{\mathfrak v}^T,\qquad {\mathfrak v} \in [0,1).
$$

\noindent Under ${\mathbb Q}$, ${\mathcal M}$ has drift $-\frac{\delta'}{2} {\bf 1}_{[0,1]\times [0,T]}(u,x) \d u \d x$, and ${\mathcal Y}^T$ is a general {\rm Jacobi($\delta,\delta'$)} flow  on $[0,T]$. 
\end{theorem}

\noindent {\it Proof}. We use the representation in \eqref{def:M}: ${\mathcal M}_T([0,u]) = \widetilde {\mathcal W}_T([0,u]) - u \widetilde {\mathcal W}_T([0,1])$ where $\widetilde {\mathcal W}$ is a white noise and we used the notation $\widetilde {\mathcal W}_s(A):=\widetilde {\mathcal W}(A\times [0,s])$ for $s\ge 0$ and Borel sets $A$.  The random variable
\begin{equation}\label{eq:changeW}
{\rm e}^{-\frac{\delta'}{2(1-{\mathfrak v})}\widetilde {\mathcal W}_T([0,{\mathfrak v}]) - \frac{T}{8} \frac{\delta'^2}{(1-{\mathfrak v})^2} {\mathfrak v} }
\end{equation}

\noindent has mean $1$. Take this random variable as the Radon--Nikodym derivative of a new probability measure  with respect to $\p$ on $\sigma(\widetilde {\mathcal W}_s(A), 0\le s \le T, A \hbox{ Borel set of } [0,1])$.  This change of measure  adds a drift $-\frac{\delta'}{2(1-\mathfrak v)}{\bf 1}_{[0,{\mathfrak v}]\times [0,T]}(u,x) \d u \d x$ to the white noise $\widetilde {\mathcal W}$. Therefore it adds a drift 
$$
-\frac{\delta'}{2}\Big(\frac{1}{(1-\mathfrak v)} {\bf 1}_{[0,{\mathfrak v}]\times [0,T]}(u,x) - \frac{\mathfrak v}{1-\mathfrak v} \Big)  = -\frac{\delta'}{2} \Big( {\bf 1}_{[0,{\mathfrak v}]\times [0,T]}(u,x) - \frac{\mathfrak v}{1-\mathfrak v} {\bf 1}_{({\mathfrak v},1]\times [0,T]}(u,x)\Big)
$$ 

\noindent  to the martingale measure ${\mathcal M}$. Use that ${\mathcal M}_T([0,{\mathfrak v}]) = \widetilde {\mathcal W}_T([0,{\mathfrak v}]) - {\mathfrak v} \widetilde {\mathcal W}_T([0,1])$, and average \eqref{eq:changeW} over $\widetilde {\mathcal W}_T([0,1])$ which is independent of ${\mathcal M}$ to obtain that it is also the distribution of ${\mathcal M}$ under $m_{\mathfrak v}^T \cdot {\mathbb P}$. Notice that on $[0,\mathfrak v]\times [0,T]$, the drift is just $-\frac{\delta'}{2}$. We conclude that $m^T_{\mathfrak v}$ is the Radon-Nykodim derivative  on ${\mathcal F}_{\mathfrak v}^T$ of the probability measure under  which ${\mathcal M}$ has drift $-\frac{\delta'}{2} {\bf 1}_{[0,1]\times [0,T]}(u,x)$. Hence $m^T$ is a martingale, and ${\mathcal M}$ under ${\mathbb Q}$ has the required distribution (we could directly check that $ m_{\mathfrak v}^T$ is a ${\mathcal F}_{\mathfrak v}^T$-martingale, by using the semimartingale decomposition of the bridge: ${\mathcal M}_T([0,u])= \sqrt{T} \widehat  B_u - \int_0^u \frac{\d s}{1-s} {\mathcal M}_T([0,s])$ where $\widehat  B$ is a standard ${\mathcal F}^T$-Brownian motion.).

Plug the drift into \eqref{def:Yflow} to conclude that a flow line  from $v\in [0,1]$ and $s\in [0,T]$ is  a (possibly absorbed at $0$) Jacobi($\delta,\delta'$) process up to  the minimum between  the hitting time of $1$ and $T$. We set, under ${\mathbb Q}$, ${\mathcal Y}_{s,t}^T(1):=1$ for all $0\le s\le t \le T$, and absorb at $1$ any flow line which hits $1$.

Let $\mathfrak v\in [0,1)$. From the absolute continuity of $\mathbb Q$ with respect to $\mathbb P$ on ${\mathcal F}_{\mathfrak v}^T$, one concludes that the almost perfect flow property of Corollary \ref{c:almost} holds as long as the flow lines do not touch $\mathfrak v$. Making $\mathfrak v$ tend to $1$, this property is also true as long as they do not touch $1$. With our convention, it is also true after the hitting time of $1$. Hence  ${\mathcal Y}^T$ possesses the almost perfect flow property also under $\mathbb Q$. We can similarly show that flow lines coalesce under ${\mathbb Q}$ in the sense of Corollary \ref{c:almost}. 

Let us check the regularity conditions of Definition \ref{def:Yflow}.  Statement (i) is clear. We prove (ii). We already know that  $v \in[0,1]\to {\mathcal Y}^T_{s,t}(v)$ is nondecreasing since flow lines coalesce. Again, for any $\mathfrak v\in (0,1)$, the right-continuity holds for any $0\le s\le t\le T$ at any $v\in [0,1)$ such that the flow line ${\mathcal Y_{s,\cdot}}(v)$ did not hit $\mathfrak v$ on $[s,t]$. Then, it holds if the flow line did not hit $1$. Since ${\mathcal Y_{s,t}}(v)=1$ if the flow line hit $1$, we conclude that it is right-continuous  at any $v\in [0,1)$. Condition (iii) is a consequence of the fact that flow lines coalesce and condition (iv) is satisfied by construction. $\Box$

%%%%%%%%%%%

\section{Proofs of Lemmas \ref{l:continuity} and \ref{l:continuity2}}
\label{s:continuitylemmas}

\noindent {\it Proof of Lemma  \ref{l:continuity}}.    Fix $M\ge 0$ and take a compact set $L$ as in the statement of the lemma, for $K=[0,M]$. By assumption, there exists $\varepsilon\in (0,1)$ such that $f^n$ and $f$ are contained in  $(\varepsilon,\varepsilon^{-1})$ on $L$ for $n$ large enough. Let $t\in [0,M]$. We have, since $\mathcal X([0,M])\cup \mathcal X^n([0,M])\subset L \subset I_n$ for $n$ large enough, 
$$
|\eta_{f}(\mathcal X_t)-\eta_{f^n}(\mathcal X_t^n)| \le \bigg|\int_{0}^{\mathcal X_t} \d r \left({ 1 \over f(r)}-{1\over f^n(r)}\right)\bigg| + \bigg|\int_{\mathcal X_t}^{\mathcal X_t^n}  {\d r\over f^n(r)}\bigg|.
$$

\noindent  The first term is bounded by $|\!|{1\over f}-{1\over f^n}|\!|_{L}|\!|\mathcal X |\!|_{[0,M]}$. The second term is bounded by ${1\over \varepsilon}|\!|\mathcal X-\mathcal X^n |\!|_{[0,M]}$. Both terms go to $0$ as $n\to \infty$. 
 It shows the uniform convergence of $\eta_{f^n}\circ \mathcal X^n$ to $\eta_{f}\circ \mathcal X$ on $[0,M]$. Since $|f(\mathcal X_t) -  f^n(\mathcal X^n_t)|\le |f(\mathcal X_t) -  f(\mathcal X^n_t)| + |f(\mathcal X^n_t) -  f^n(\mathcal X^n_t)|$, that $f$ is uniformly continuous on $L$, that $\mathcal X^n$ converges uniformly to $\mathcal X$ on $[0,M]$ and $f^n$ converges uniformly to $f$ on $L$, we have that $f^n\circ \mathcal X^n$ converges uniformly to $f\circ\mathcal X$ on $ [0,M]$. Let $C_{f_n}(t):=\int_0^t {\d s \over f^n(\mathcal X^n_s)^2}$ for $t\in [0,j_n)$ and $C_{f}(t):=\int_0^t {\d s \over f(\mathcal X_s)^2}$ for $t\ge 0$. Notice that the range of $C_{f_n}$ and of $C_f$ is $\r_+$ by our assumption (i).  Dominated convergence implies that $C_{f^n}(t)$ converges to $C_f(t)$ for all $t\in [0,M]$. Fix $t\in (0,M)$ and set $x:=C_f(t)$. Let $t_n:=C_{f_n}^{-1}(x)$. Since $C_{f^n}(v)-C_{f^n}(u)\ge \varepsilon^{2} (v-u)$ for $v\ge u$ in $[0,M]$, we have $|C_{f^n}(t)-x|\ge \varepsilon^{2} |t-t_n|$. We deduce that $t_n\to t$, i.e. $C_{f^n}^{-1}(x) \to C_f^{-1}(x)$. Since $M$ is arbitrary, it proves that $C_{f^n}^{-1}$ converges to $C_f^{-1}$ pointwise on $\r_+$. Since $C_{f_n}^{-1},C_f^{-1}$ are  increasing continuous functions on $\r_+$, it implies the uniform convergence on any compact set of $\r_+$. 
By composition of $\eta_{f^n}\circ \mathcal X^n$ and $C_{f_n}^{-1}$, we deduce that $\Upsilon({\mathcal X}^n,f^n)$ converges to $\Upsilon({\mathcal X},f)$ uniformly on any compact. 
 $\Box$}

\medskip
\noindent {\it Proof of Lemma  \ref{l:continuity2}}.  
First, for any $t\ge 0$, a.s.,  $1-L^{(2)}(t,x)>0$ for all $x\in \{Z^{(2)}_{t+u},\, u\ge 0\}$. Let us prove this statement. Because $Z^{(2)}$ accumulates some local time at each level that it crosses, if  $L^{(2)}(t,x)=1$, it implies that the process $Z^{(2)}$ will never visit $x$ again. Hence necessarily, the only possibility to have $L^{(2)}(t,x)=1$ for some $x\in \{Z^{(2)}_{t+u},\, u\ge 0\}$ is to have $x=Z^{(2)}_t$. So we need to show that for all $t\ge 0$, a.s.,  $L^{(2)}(t,Z^{(2)}_t)<1$. Let $\varepsilon\in (0,1)$. We have by the occupation times formula, for any $s<t$ and $M\ge 0$, $\int_s^t {\bf 1}_{\{L^{(2)}(u,Z^{(2)}_u) \ge 1-\varepsilon,\, |Z_u^{(2)}| \le M \}} \d u=\int_{x\in [-M,M]} \int_s^t {\bf 1}_{\{L^{(2)}(u,x) \ge 1-\varepsilon\}} \d_u L^{(2)}(u,x) \d x\le  2M \varepsilon (t-s)$. Taking expectation, we deduce that  $t$ is an accumulation of points $s$ such that  $\p(L^{(2)}(s,Z^{(2)}_s) \ge 1-\varepsilon,\, |Z^{(2)}_s|\le M) \le 2 M\varepsilon$, hence by continuity of $s\to L^{(2)}(s,Z^{(2)}_s)$ and of $Z^{(2)}$, $\p(L^{(2)}(t,Z^{(2)}_t)=1, |Z^{(2)}_t|< M)\le 2 M\varepsilon $ which proves the claim by making $\varepsilon\to 0$ then $M\to \infty$. It implies that    $\Upsilon(Z^{(2)}_{t+\cdot}, 1-L^{(2)}(t,\cdot))$ is well-defined. Notice that by continuity of $Z^{(2)}$, it implies that $\Upsilon(Z^{(2)}_{s+\cdot}, 1-L^{(2)}(s,\cdot))$ is well-defined on a neighborhood of $t$.

Let us go back to the proof of the continuity and fix $t\ge 0$. Let $t'<t<t''$ and $x'':=\inf\{x\,:\,  L^{(2)}(t'',x)=1\}$. We suppose that $t'$ and $t''$ are close enough to $t$ so that $x''> \sup_{\{u\ge t'\}}Z^{(2)}_u$ (such  $t'$ and $t''$ exist a.s. from what we just proved). Let $(t_n)_n$  be a sequence which converges to $t$. We can suppose that $t_n\in (t',t'')$ for all $n$ for simplicity.  We apply Lemma \ref{l:continuity} with $\mathcal X_u:=Z^{(2)}_{t+u}-Z^{(2)}_t$, $f(x):=1-L^{(2)}(t,x+Z^{(2)}_t)$, $I=(-\infty,x''-Z^{(2)}_t)$, and $\mathcal X^n$, $f^n$, $I_n$ defined similarly by replacing $t$ by $t_n$.  Let us check the assumptions of  Lemma \ref{l:continuity}. In the notation of that lemma, (i) and (ii) are verified since all processes start at $0$, $j_n = \infty$ and $f(x)=f^n(x)=1$ for $x\to -\infty$. Assumption (iii) comes from the uniform continuity of $Z^{(2)}$ on compact sets. Let $K$ be a compact set of $\r_+$. Take $\ell_+ \in \r$ such that $\inf_{(t',t'')} Z^{(2)} + \ell_+ > \sup_{u\ge t'} Z^{(2)}_u$ and $\sup_{(t',t'')} Z^{(2)} + \ell_+ < x''$  (we can suppose that $t$ and $t''$ are close enough so that $\ell_+$ exists).  Then for all $s\ge 0$,  $Z^{(2)}_{t+s}-Z^{(2)}_t \le \ell_+$. Indeed, $Z^{(2)}_{t+s}-Z^{(2)}_t \le \sup_{u\ge t'} Z^{(2)}_u - \inf_{(t',t'')}Z^{(2)} \le \ell_+$. The same lines replacing $t$ by $t_n$ imply that for all $s\ge 0$, $Z^{(2)}_{t_n+s}-Z^{(2)}_{t_n} \le \ell_+$ for all $n$. Similarly, take $\ell_-\in \r$ such that $\sup_{(t',t'')} Z^{(2)} + \ell_- < \inf_{u\le t'' + \sup K} Z^{(2)}_{t''+ u}$. We have that for all $s\in K$,  $Z^{(2)}_{t+s}-Z^{(2)}_t \ge \ell_-$ since $Z^{(2)}_{t+s}-Z^{(2)}_t \ge \inf_{u\le t''+\sup K} Z^{(2)}_{u} - \sup_{t',t''}Z^{(2)} \ge \ell_-$, similarly for $t_n$. The compact set $L=[\ell_-,\ell_+]$ satisfies assumption (iv). We deduce the convergence of $\Upsilon(Z^{(2)}_{t_n +\cdot},1-L^{(2)}(t_n,\cdot))$ towards $\Upsilon(Z^{(2)}_{t +\cdot},1-L^{(2)}(t,\cdot))$ (on an event of probability $1$ which does not depend on the choice of $(t_n)_n$). $\Box$

\bigskip
{\noindent\bf Acknowledgments.}
 We thank  Hui He, Zenghu Li, Titus Lupu, Emmanuel Schertzer, Anton Wakolbinger for helpful discussions.


\begin{thebibliography}{99}
\label{references}

%\baselineskip=14pt
\def\&{et}



\bibitem{elie}
	A\"id\'ekon, E.\ (2020).
Cluster explorations of the loop soup on a metric graph related to Gaussian free field. \href{https://arxiv.org/abs/2009.05120}{\tt arXiv:2009.05120}
	
	
\bibitem{eyzRK}
	A\"id\'ekon, E., Hu, Y.\ and Shi, Z.\ (2020).
	An infinite-dimensional representation of the Ray--Knight theorems.
	{\it Science China Mathematics} (to appear).


\bibitem{eyzmupro}
    A\"id\'ekon, E., Hu, Y.\ and Shi, Z.\ (2020).
    Path decompositions of perturbed reflected Brownian motions. {\it Chaumont, L. and Kyprianou, A.E. (ed.), A lifetime of excursions through random walks and Lévy processes. A volume in honour of Ron Doney's 80th birthday.}  Cham: Birkh\"{a}user. Prog. Probab. 78, 13--42. 



\bibitem{aldous}
	Aldous, D.J.\ (1998).
	Brownian excursion conditioned on its local time.
	{\it Elect.\ Comm.\ Probab.} {\bf 3}, 79--90.


\bibitem{arratia}
	Arratia, R.A. (1979).
	Coalescing Brownian motions on the line.
	PhD thesis. University of Wisconsin.

\bibitem{bassburdzy}
	Bass, R.F., and Burdzy, K. \ (1999).
	Stochastic Bifurcation Models.
	{\it Ann.\ Probab.}, {\bf 27}, 50--108. 
	
\bibitem{berestycki}		
	Berestycki, J., and Berestycki, N.\ (2009).
	Kingman's coalescent and Brownian motion.
	{\it ALEA} {\bf 6}, 253--259.



\bibitem{bertoin-legall-I}
	Bertoin, J.\ and Le Gall, J.-F.\ (2003).
	Stochastic flows associated to coalescent processes.
	{\it Probab. Theory Relat. Fields} {\bf 126}, 261--288.


\bibitem{bertoin-legall-II}
	Bertoin, J.\ and Le Gall, J.-F.\ (2005).
	Stochastic flows associated to coalescent processes II: stochastic differential equations.
	{\it Ann.\ Inst.\ H.\ Poincar\'e Probab.\ Statist.} {\bf 41}, 307--333.

 

\bibitem{carmona-petit-yor94a}
	Carmona, P., Petit, F.\ and Yor, M. (1994).
	Sur les fonctionnelles exponentielles de certains processus de L\'evy.
	{\it Stochastics and Stochastic Reports} {\bf 42}, 71--101.
	
 

\bibitem{dawson-li12}
	Dawson, D.A.\ and Li, Z.\ (2012).
	Stochastic equations, flows and measure-valued processes.
	{\it Ann.\ Probab.} {\bf 40}, 813--857.


\bibitem{DuLG02}
	Duquesne, T. and Le Gall, J.-F.\ (2002).
	{Random Trees, L\'evy Processes and Spatial Branching Processes}. 
	{\it Ast\'erisque} {\bf 281}.

\bibitem{foucart12}
	Foucart, C. (2012).
	Generalized Fleming-Viot processes with immigration via stochastic flows of partitions.
	{\it ALEA.} {\bf 9}, 451--472.

 
\bibitem{fmm}
	Foucart, C.,\, Ma, C.\ and Mallein, B.\ (2019).
	Coalescences in Continuous-State Branching Processes.
	{\it Electron. J. Probab.}, {\bf 24}, Paper no.~103, 1--52.
	
\bibitem{gy03} Going-Jaeschke, A.\  and Yor, M.\ (2003). 
    A survey and some generalizations of Bessel processes. 
    {\it Bernoulli} {\bf 9}, 313--349.

\bibitem{GS00}
	Graversen, S., E. and Shiryaev, A.N. (2000).
	An extension of P. L\'evy's distributional properties to the case of a Brownian motion with drift.
	{\it Bernoulli}, {\bf 6} (4), 615--620.

\bibitem{GKW}
	Gufler, S., Kersting, G.\ and Wakolbinger, A.\ (2020).
	A decomposition of the Brownian excursion.
	(in preparation)

\bibitem{hw00}
	Hu, Y.\ and Warren. J.\ (2000).
	Ray--Knight theorems related to a stochastic flow.
	{\it Stoch.\ Proc.\ Appl.} {\bf 86}, 287--305.

 


\bibitem{legall-yor}
	Le Gall, J.-F.\ and Yor, M.\ (1986).
	Excursions Browniennes et carr\'es de processus de Bessel.
	{\it C.\ R.\ Acad.\ Sci.\ Paris} {\bf 302}, 641--643.

 \bibitem{lupu16}
    Lupu, T.\ (2016).
    From loop clusters and random interlacements to the free field.
    {\it Ann.\ Probab.} {\bf 44}, 2117--2146. 

\bibitem{lupu18}
    Lupu, T.\ (2018).
    Poisson ensembles of loops of one-dimensional diffusions.
    {\it M\'em.\ Soc.\ Math.\ France}, {\bf 158}.

 \bibitem{lst19}
    Lupu, T., Sabot, C.\ and Tarr\`es, P.\ (2019).
    Inverting the coupling of the signed Gaussian free field with a loop-soup.
    {\it Electron.\ J.\ Probab.} {\bf 24}, Paper no.~70, 1--28.

 \bibitem{lst20}
    Lupu, T., Sabot, C.\ and Tarr\`es, P. \ (2021).
    Inverting the Ray--Knight identity on the line.
    %\href{https://arxiv.org/abs/1910.06836}{\tt arXiv:1910.06836}
{\it Electron.\ J.\ Probab.} {\bf 26}, Paper no.~96, 1--25.


\bibitem{pal}
	Pal, S.\ (2013).	
	Wright--Fisher diffusion with negative mutation rates.
	{\it Ann.\ Probab.} {\bf 41}, 503--526.
 
\bibitem{perkins}
    Perkins, E.\ (1992).	
    Conditional Dawson--Watanabe processes and Fleming--Viot processes.
    {\it Sem.\ Stoch.\ Proc.\ 1991} pp.~143--156,
    Progr.\ Probab.\ {\bf 29}, Birkh\"auser, Boston.



\bibitem{perman}	
	Perman, M.\  (1996).
	An excursion approach to Ray--Knight theorems for perturbed Brownian motion.
	{\it Stoch.\ Proc.\ Appl.} {\bf 63}, 67--74.
	

\bibitem{perman-werner}
    Perman, M.\ and Werner, W.\ (1997).
    Perturbed Brownian motions.
    {\it Probab.\ Theory Related Fields} {\bf 108}, 357--383.


	
\bibitem{protter}
	Protter, Ph. E.\ (2004). Stochastic integration and differential equations. 
	Second edition. 
	Springer, Berlin. 
 	
\bibitem{revuz-yor}
    Revuz, D.\ and Yor, M.\ (1999).
    Continuous martingales and Brownian motion.
    Third edition.
    Springer, Berlin.

\bibitem{st16}
	Sabot, C.\ and Tarr\`es, P.\ (2016).
	Inverting Ray--Knight identity.
	{\it Probab. Theory Relat. Fields} {\bf 165}, 559--580.
	
 
\bibitem{toth-werner}	
	T\'oth, B. \ and Werner, W. \ (1998).
	The true self-repelling motion.
	{\it Probab.\ Theory Related Fields} {\bf 111}, 375--452. 
	
 \bibitem{walsh84}
    Walsh, J.B.\ (1986).
    An introduction to stochastic partial differential equations.
    {\it \'Ecole d'\'et\'e de probabilit\'es de Saint-Flour XIV 1984}, pp.~265--439, 
    Lecture Notes Math.\ {\bf 1180}, Springer, Berlin.

\bibitem{warren}
	Warren, J.\ (2005).
	A stochastic flow arising in the study of local times.
	{\it Probab. Theorey Relat. Fields}, {\bf 133}, pp. 559--572.


\bibitem{warren-yor}
    Warren, J.\ and Yor, M.\ (1998).
    The Brownian burglar: conditioning Brownian motion by its local time process. 
    {\it S\'em.\ Probab.\ XXXII}, pp.~328--342,
    Lecture Notes Math.\ {\bf 1686}, Springer, Berlin.
 


\bibitem{wernermarkov}
    Werner, W. (2016).
    On the spatial Markov property of soups of unoriented and oriented loops.
    {\it S\'em.\ Probab.\ XLVIII}, pp.~481--503, Lect.\ Notes Math.\ {\bf 2168}, Springer, Berlin.

\bibitem{yor92}
    Yor, M.\ (1992).
    Some aspects of Brownian motion. Part I.
    %{\it Lectures in Mathematics ETH Z\"{u}rich}, 
    Birkh\"{a}user, Basel.

 \end{thebibliography}
\end{document}